\documentclass[11pt]{article}

\usepackage[margin=1in]{geometry}
\usepackage{setspace}
\usepackage[T1]{fontenc}
\usepackage{amsmath,amssymb,amsthm,mathtools}
\RequirePackage{tgtermes}
\RequirePackage{newtxtext}
\RequirePackage{newtxmath}
\RequirePackage{bm}
\RequirePackage{endnotes}

\usepackage{algorithm}
\usepackage{algpseudocode}
\usepackage{graphicx}
\usepackage{subcaption}
\usepackage{tikz}
\usetikzlibrary{positioning,arrows.meta,fit,backgrounds,shapes.geometric}

\usepackage[sort&compress]{natbib}
\bibpunct[, ]{[}{]}{,}{n}{}{,}

\usepackage{xcolor}
\definecolor{citeblue}{RGB}{0,90,160}

\usepackage[
    colorlinks=true,
    linkcolor=black,
    citecolor=citeblue,
    urlcolor=citeblue,
    pdfauthor={\.I.\ Esra B{\"u}y{\"u}ktahtak{\i}n},
    pdftitle={Deep Learning for Sequential Decision Making under Uncertainty: Foundations, Frameworks, and Frontiers}
]{hyperref}

\usepackage[nameinlink,noabbrev]{cleveref}

\theoremstyle{definition}

\title{\textbf{Deep Learning for Sequential Decision Making under Uncertainty: Foundations, Frameworks, and Frontiers}}

\author{
\.I.\ Esra B{\"u}y{\"u}ktahtak{\i}n\\[0.2em]
Grado Department of Industrial and Systems Engineering, Virginia Tech\\
Blacksburg, VA 24061, USA\\
\texttt{esratoy@vt.edu}
}

\date{}

\begin{document}

\maketitle


\vspace{0.75em}

\begin{abstract}
Artificial intelligence (AI) is moving increasingly beyond prediction to support
decisions in complex, uncertain, and dynamic environments. This shift creates a natural intersection with operations research and management science (OR/MS), which has long provided methodological foundations for sequential decision making under uncertainty. At the same time, deep learning advances, including feedforward 
neural networks, recurrent architectures, transformers, large language models (LLMs), 
and deep reinforcement learning, have expanded data-driven modeling for large-scale 
decisions. This tutorial presents an OR/MS-centered perspective on deep learning for 
sequential decision making under uncertainty, bridging neural architectures and OR/MS 
approaches to decision making. Its premise: deep learning complements optimization 
rather than replacing it. Deep learning brings adaptability and scalable approximation, 
whereas OR/MS provides the mathematical rigor to represent constraints, recourse, uncertainty, and 
decision quality. The tutorial reviews key decision making foundations, connects them to the major neural 
architectures in modern AI, and organizes the field around three central themes: 
predict-then-optimize and decision-aware learning, learning-based decision generation 
under constraints for continuous and discrete problems with temporal coupling, and 
deep reinforcement learning for sequential and combinatorial decision making. Impact spans supply chains, service 
systems, healthcare and epidemic response, agriculture, energy, environmental 
sustainability, and autonomous operations. This tutorial frames these developments as part of a 
shift from predictive AI toward decision-capable AI, highlighting OR/MS's role in 
shaping the next generation of integrated learning--optimization systems.
\end{abstract}


\vspace{0.75em}

\section{Introduction and Motivation}

Artificial intelligence (AI) is increasingly moving beyond pattern recognition and prediction toward systems that support, recommend, and sometimes automate decisions. In many domains central to operations research and management science (OR/MS), including healthcare, supply chains, energy, finance, transportation, and public policy, the fundamental challenge is not only to anticipate what may happen, but to determine what should be done as uncertainty unfolds over time. Sequential decision making under uncertainty thus represents a central point of contact between modern AI and OR/MS.  This tutorial adopts a decision-centered lens: how deep learning can support implementable sequential decisions, rather than merely improve forecasts.

Recent advances in deep learning have greatly expanded the reach of data-driven modeling. Feedforward neural networks (FNNs) have proven effective for nonlinear approximation in static settings, recurrent architectures, such as long short-term memory networks (LSTMs), have strengthened learning in temporal environments \citep{goodfellow2016deep,hochreiter1997long}, and transformer architectures now underpin large language models (LLMs), generative AI, and conversational systems such as ChatGPT \citep{vaswani2017attention}. Deep reinforcement learning has further broadened the ability to learn sequential decision policies through interaction in complex dynamic settings \citep{sutton2018reinforcement}. Together, these developments have substantially increased both the modeling capacity and the practical scope of modern learning systems.

At the same time, their growing use in operational settings has made a longstanding distinction impossible to ignore: prediction is not the same as decision making. A highly accurate predictive model does not necessarily produce a high-quality decision, especially when actions are constrained, uncertainty evolves over time, and downstream consequences matter \citep{elmachtoub2022smart}. In sequential settings, the distinction becomes sharper still, because present decisions shape future opportunities, risks, and feasible actions. The relevant objective is, therefore, rarely prediction alone, but the design of policies that remain feasible, adaptive, and effective over time. This is where OR/MS provides an essential perspective and methodology. The field has long developed rigorous frameworks for sequential decision making under uncertainty, including dynamic programming, Markov decision processes, stochastic programming, approximate dynamic programming, and reinforcement learning \citep{bellman1957dynamic,puterman1994markov,birge2011introduction,shapiro2009lectures,sutton2018reinforcement}. These frameworks place decisions at the center of the analysis by explicitly accounting for constraints, recourse, uncertainty, information structure, and system-level trade-offs. They ask not merely whether a model predicts well, but whether a policy is implementable, respects what is known at each stage, and improves operational performance.

However, classical decision making frameworks face serious challenges in contemporary settings. Many sequential problems involve large state spaces, long planning horizons, heterogeneous data streams, complex constraints, and mixed-integer structure, rendering exact dynamic programming intractable and large-scale stochastic programs computationally burdensome. Deep learning brings important strengths to this landscape: neural architectures can approximate complex functional relationships, learn compact representations from large datasets, and in some settings generalize across families of related problem instances rather than solving each new instance entirely from scratch \citep{bengio2021machine}. The central question then is not whether learning should replace optimization, but how learning and optimization can be integrated so that each contributes what it does best.

One useful lens for this integration is the \emph{predict-then-optimize} paradigm, in which uncertain inputs are first predicted and then used in a downstream optimization model. This modular approach is natural in many OR/MS applications, but it also highlights a broader methodological challenge: learning models should often be evaluated not only by statistical prediction accuracy but also by the quality of the decisions they induce. This perspective has motivated several related streams of work, including decision-focused learning, end-to-end task-based learning, prescriptive analytics, differentiable optimization layers, and learning-to-optimize methods \citep{donti2017task,elmachtoub2022smart,bertsimas2020predictive,amos2017optnet,agrawal2019differentiable,bengio2021machine, yilmaz2024expandable}. Related model-based OR/MS work similarly shows how biological, epidemiological, operational, or learned predictive models can serve as components within optimization models, in a \emph{predict-and-optimize} fashion that links system state dynamics to downstream decisions \citep{buyuktahtakin2011dynamic,buyuktahtakin2018ebola,fischetti2018deep,lee2025transformer}.

This tutorial therefore treats deep learning architectures not merely as predictive tools, but as components of sequential decision systems under uncertainty. It focuses on how these architectures can be integrated with OR/MS frameworks in ways that support feasibility, non-anticipativity, recourse, generalization, and decision quality. The central questions are what should be learned in sequential settings, how reinforcement learning and stochastic programming offer complementary perspectives, and how newer architectures such as transformers and large language models can be incorporated into optimization without sacrificing constraints, objectives, and information structure. More broadly, the tutorial argues that the next phase of AI for high-impact operations will depend not only on more powerful learning architectures but also on stronger OR/MS decision frameworks.

This tutorial is intended for researchers, graduate students, and practitioners familiar with optimization, stochastic modeling, decision analysis, or applied machine learning, but not necessarily with deep learning. Rather than providing an exhaustive survey of AI or reinforcement learning, it organizes the learning--optimization landscape around three paradigms: integrating prediction and optimization, learning-based decision generation under constraints for continuous, discrete, and mixed decision problems, and deep reinforcement learning for sequential and combinatorial decision making.

By the end of this tutorial, readers should be able to: (i) distinguish how Markov decision processes (MDPs), dynamic programming, reinforcement learning, and multi-stage stochastic programming formalize sequential decisions under uncertainty; (ii) identify neural architectures suited to static, relational, temporal, and contextual decision structures; (iii) compare predict-then-optimize, decision-aware learning, and learning-to-optimize as learning--optimization integration paradigms; (iv) assess learning-based decision generation under constraints with respect to feasibility, expandability, temporal coupling, non-anticipativity, and recourse; (v) explain the role and limitations of deep reinforcement learning in constrained OR/MS settings; and (vi) recognize representative application domains and emerging research directions for integrated learning--optimization systems. The examples discussed throughout are illustrative rather than exhaustive and are used to highlight recurring modeling principles at the interface of learning and optimization.

\paragraph{Tutorial Overview.}
The remainder of this tutorial is organized as follows. Section~\ref{sec:preliminaries} reviews the foundations of sequential decision making under uncertainty, including Markov decision processes, dynamic programming, reinforcement learning, and multi-stage stochastic programming. Section~\ref{sec:architectures} introduces the major deep learning architectures relevant to decision systems, with emphasis on FNNs, graph neural networks, recurrent architectures such as LSTMs, and transformer attention mechanisms. Section~\ref{sec:integration} discusses key approaches for combining learning and optimization, including predict-then-optimize, decision-aware learning, and learning-to-optimize. Section~\ref{sec:learning_under_constraints} focuses on sequential decision learning under constraints, with particular attention to feasibility, expandability, integrality, temporal coupling, and non-anticipativity. Section~\ref{sec:drl} examines deep reinforcement learning in sequential and combinatorial settings. Section~\ref{sec:applications_impact} highlights representative applications and broader interdisciplinary impact. Section~\ref{sec:challenges} discusses open challenges and research frontiers. Finally, Section~\ref{Conclusion} concludes the tutorial.

\section{Preliminaries: Sequential Decision Making under Uncertainty}
\label{sec:preliminaries}

\subsection{Foundations of Sequential Decision Making}
\label{sec:preliminaries_foundations}

\begin{figure*}[h]
\centering

\begin{subfigure}[t]{0.24\textwidth}
\centering
\begin{tikzpicture}[
    >=Latex,
    every node/.style={font=\scriptsize},
    state/.style={circle,draw,minimum size=7mm,inner sep=1pt},
    action/.style={rectangle,draw,rounded corners,minimum width=9mm,minimum height=6mm},
    lbl/.style={fill=white, inner sep=1pt}
]
\node[state] (s1) at (0,0) {$s_t$};
\node[action] (a1) at (1.35,0) {$a_t$};
\node[state] (s2) at (3.25,0) {$s_{t+1}$};

\draw[->] (s1) -- (a1);
\draw[->] (a1) -- (s2);

\node[lbl] at (2.3,0.53) {$P(\cdot\mid s_t,a_t)$};
\node[align=center] at (1.65,-1.0) {$c(s_t,a_t)$\\Markov state transition};
\end{tikzpicture}
\caption{MDP}
\end{subfigure}
\hfill
\begin{subfigure}[t]{0.24\textwidth}
\centering
\begin{tikzpicture}[
    >=Latex,
    every node/.style={font=\scriptsize},
    box/.style={rectangle,draw,rounded corners,minimum width=1.35cm,minimum height=0.6cm,align=center}
]
\node[box] (vtop) at (0,1.6) {$V_{t-1}(s)$};
\node[box] (vmid) at (0,0.65) {$V_t(s)$};
\node[box] (vbot) at (0,-0.3) {$V_{t+1}(s)$};

\draw[->] (vbot) -- (vmid);
\draw[->] (vmid) -- (vtop);

\node[align=center] at (0,-1.05) {backward recursion\\via Bellman equation};
\end{tikzpicture}

\vspace{-1.5mm}
{\scriptsize
$\displaystyle
V_t(s)=\min_{a\in\mathcal{A}(s)}
\left\{c_t(s,a)+\mathbb{E}\!\left[V_{t+1}(S_{t+1})\,\middle|\,S_t=s,\;A_t=a\right]\right\}
$
}

\vspace{-5mm}
\caption{DP}
\end{subfigure}
\hfill
\begin{subfigure}[t]{0.24\textwidth}
\centering
\hspace*{-2mm}
\begin{tikzpicture}[
    >=Latex,
    every node/.style={font=\scriptsize},
    box/.style={rectangle,draw,rounded corners,minimum width=1.5cm,minimum height=0.7cm,align=center},
    lbl/.style={fill=white, inner sep=1pt}
]
\node[box] (agent) at (0,0) {Agent\\$\pi_\theta(a\mid s)$};
\node[box] (env) at (2.65,0) {Environment};

\draw[->] (agent) -- (env);
\draw[->] (env) to[bend right=20] (agent);

\node[lbl] at (1.3,0.26) {$a_t$};
\node[lbl] at (1.25,-0.42) {$s_{t+1},\,r_t$};
\node[lbl] at (0,0.82) {update $\theta$};
\node[align=center] at (1.3,-1.02) {learn policy/value\\from interaction data};
\end{tikzpicture}
\caption{RL}
\end{subfigure}
\hfill
\begin{subfigure}[t]{0.24\textwidth}
\centering
\begin{tikzpicture}[
    >=Latex,
    every node/.style={font=\scriptsize},
    dec/.style={circle,draw,minimum size=7mm,inner sep=1pt},
    lbl/.style={fill=white, inner sep=1pt}
]
\node[dec] (r)  at (0,1.95) {$x_1$};
\node[lbl] at (0,2.55) {first-stage decision};

\node[dec] (n1) at (-0.75,0.95) {$x_2^1$};
\node[dec] (n2) at ( 0.75,0.95) {$x_2^2$};

\node[dec] (l1) at (-1.45,-0.45) {$x_3^1$};
\node[dec] (l2) at (-0.35,-0.45) {$x_3^2$};
\node[dec] (l3) at ( 0.35,-0.45) {$x_3^3$};
\node[dec] (l4) at ( 1.45,-0.45) {$x_3^4$};

\draw (r) -- (n1);
\draw (r) -- (n2);
\draw (n1) -- (l1);
\draw (n1) -- (l2);
\draw (n2) -- (l3);
\draw (n2) -- (l4);

\node[align=center] at (0,-1.25) {scenario tree\\non-anticipativity and recourse};
\end{tikzpicture}
\caption{MSP}
\end{subfigure}
\caption{Comparison of four sequential decision making frameworks. (a) A Markov decision process (MDP) models stochastic state transitions under actions and stage costs. (b) Dynamic programming (DP) solves such problems through Bellman recursion over value functions. (c) Reinforcement learning (RL) learns policies or value functions from interaction with an environment. (d) Multi-stage stochastic programming (MSP) represents decisions explicitly on a scenario tree with recourse and non-anticipativity constraints.}
\label{fig:sequential_frameworks_comparison}
\end{figure*}
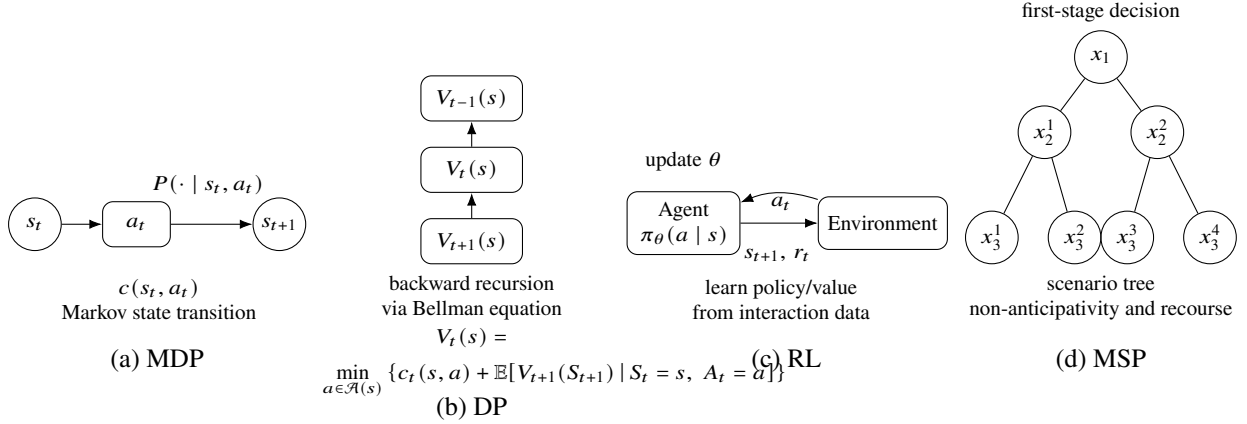

Sequential decision making under uncertainty is a foundational theme in operations research, since decisions are rarely made all at once and instead must adapt over time as information is progressively revealed \citep{bellman1957dynamic,puterman1994markov,birge2011introduction,shapiro2009lectures}.
We consider a finite horizon $t=1,\dots,T$ and a stochastic process $\{\xi_t\}_{t=1}^T$, where $\xi_{[t]}:=(\xi_1,\dots,\xi_t)$ denotes the information available by stage $t$. At each stage, the decision is made after observing the current history, so
\(x_t=x_t(\xi_{[t]})\). Thus, the information pattern is as follows.
\[
\begin{aligned}
&\emph{decide}(x_1)\;\rightarrow\;\emph{observe}(\xi_2)\;\rightarrow\;
\emph{decide}(x_2)\;\rightarrow\;\cdots \\
&\cdots\;\rightarrow\;\emph{observe}(\xi_T)\;\rightarrow\;
\emph{decide}(x_T).
\end{aligned}
\]
Accordingly, each \(x_t\) may depend only on \(\xi_{[t]}\), not on future
realizations; this is the \emph{non-anticipativity} principle
\citep{birge2011introduction,shapiro2009lectures}.
To allow both continuous and integer decisions, let
$x_t(\xi_{[t]}) \in \mathcal{X}_t(\xi_{[t]}) \subseteq \mathbb{R}_+^{\,n_t-q_t}\times\mathbb{Z}_+^{\,q_t}$.
The goal is to minimize the expected total cost over the horizon, where the stage-$t$ term $c_t(\xi_{[t]})^\top x_t(\xi_{[t]})$ represents the immediate cost incurred after observing $\xi_{[t]}$, and subsequent decisions provide recourse as uncertainty unfolds. Suppressing the explicit dependence on $\xi_{[t]}$ for readability, a compact multi-stage formulation is
\begin{equation}
\begin{aligned}
\min_{x}\quad
& \mathbb{E}\!\left[\sum_{t=1}^T c_t^\top x_t\right] \\
\text{s.t.}\quad
& H_1x_1=b_1, \\
& A_t x_{t-1}+H_t x_t=b_t, \qquad t=2,\dots,T, \\
& x_t \in \mathcal{X}_t, \qquad t=1,\dots,T .
\end{aligned}
\label{eq:general_multistage_compact}
\end{equation}
For expositional simplicity, \eqref{eq:general_multistage_compact} is written in a linear recourse form, but the same sequential structure extends more broadly to problems with nonlinear objectives and constraints. Here, for $t\ge2$, the quantities $A_t,H_t,b_t,c_t$, and $\mathcal{X}_t$ depend on the realized history $\xi_{[t]}$. In practice, the uncertainty process is often approximated by a finite scenario tree, where each node corresponds to a realized history, and non-anticipativity requires identical decisions at nodes sharing the same history \citep{birge2011introduction,shapiro2009lectures}. This formulation provides a common foundation for MDPs, dynamic programming, reinforcement learning, and multi-stage stochastic programming, as summarized in Figure~\ref{fig:sequential_frameworks_comparison}.

\subsection{Markov Decision Processes}

The formulation above is history-based: at stage $t$, the decision
$x_t=x_t(\xi_{[t]})$ may depend on all the information revealed so far. However, in many
settings, the information relevant to the decision in the history can be
summarized by a state variable $s_t$, leading to the framework of a Markov
decision process (MDP) \citep{puterman1994markov,bertsekas1995dynamic}.
An MDP is defined by $(\mathcal{S},\mathcal{A},P,c,\gamma,\rho_0)$, where
$\mathcal{S}$ is the state space, $\mathcal{A}(s)$ is the set of feasible
actions in state $s$, $P(s' \mid s,a)$ is the transition probability,
$c(s,a)$ is the one-stage cost, $\gamma \in (0,1]$ is a discount factor, and
$\rho_0$ is the initial state distribution. At stage $t$, the decision maker
observes $s_t$, chooses an action $a_t \in \mathcal{A}(s_t)$, the system
evolves as $s_{t+1} \sim P(\cdot \mid s_t,a_t)$, and incurs cost $c(s_t,a_t)$. Figure~\ref{fig:sequential_frameworks_comparison}(a) illustrates this state--action--transition view of an MDP.

Compared to the notation of the previous subsection, the action $a_t$ plays
the role of the stage decision, while the state $s_t$ summarizes the
 information relevant to the decision contained in the history $\xi_{[t]}$. More
precisely, when there exists a mapping $s_t=\phi_t(\xi_{[t]})$ such that
feasible actions, one-stage costs, and conditional distribution of the next
state depend only on history through $s_t$, one may restrict attention to
policies of the form $a_t=\pi_t(s_t)$ without loss of optimality
\citep{puterman1994markov,bertsekas1995dynamic}. The key assumption is the
Markov property:
\begin{equation}
\mathbb{P}(s_{t+1}\mid s_1,a_1,\dots,s_t,a_t)
=
\mathbb{P}(s_{t+1}\mid s_t,a_t),
\label{eq:markov_property}
\end{equation}
which states that, given the current state and action, the next state is
independent of the earlier history.

The objective is to find a policy $\pi$ mapping states to actions to minimize the
 expected cumulative discounted cost.
\begin{equation}
J(\pi)=\mathbb{E}_\pi\!\left[\sum_{t=1}^T \gamma^{t-1} c(s_t,a_t)\right].
\end{equation}
For a finite horizon, the optimality recursion is naturally time-indexed, as discussed in
Section~\ref{sec:dp}. In the commonly used infinite-horizon discounted stationary
specialization, the optimal value function satisfies the Bellman optimality equation
\begin{equation}
\begin{aligned}
V^*(s)=\min_{a \in \mathcal{A}(s)} \Bigl\{
c(s,a)
+\gamma \sum_{s' \in \mathcal{S}} P(s' \mid s,a)\,V^*(s')
\Bigr\}.
\end{aligned}
\end{equation}

Thus, an MDP may be viewed as a special case of a general sequential model
when the decision-relevant history can be represented by a Markov state
without loss of optimality. This abstraction is powerful and computationally
convenient, but can be restrictive in operations research applications with
complex feasibility constraints, integer decisions, and explicit
scenario-based coupling, which are often represented more naturally in
optimization-based formulations.

\subsection{Dynamic Programming}
\label{sec:dp}

Given the MDP formulation above, dynamic programming (DP) provides a recursive framework for characterizing optimal sequential decisions through value functions \citep{bellman1957dynamic,bertsekas1995dynamic}. Let $V_t(s)$ denote the optimal cost-to-go when the system is in state $s$ at stage $t$. Then the Bellman recursion is
\begin{equation}
\begin{aligned}
V_t(s)=\min_{a \in \mathcal{A}(s)} \Bigl\{ \,
&c_t(s,a) \\
&+\mathbb{E}\!\left[V_{t+1}(s_{t+1}) \mid s_t=s,\; a_t=a\right]
\Bigr\},
\end{aligned}
\end{equation}
with terminal condition \(V_{T+1}(\cdot)=0\). Figure~\ref{fig:sequential_frameworks_comparison}(b) depicts this backward-recursion view of dynamic programming.

This recursion captures the central logic of sequential decision making: the optimal value at the current stage equals the immediate cost plus the expected future cost induced by the current action. In this sense, DP provides the main analytical bridge between state-based MDP representation and optimal decision making over time. Its main limitation is computational. As the state space increases in dimension or granularity, the exact evaluation and storage of the value functions quickly become intractable, leading to the well-known curse of dimensionality \citep{bellman1957dynamic}. This challenge has motivated approximate dynamic programming, reinforcement learning, and other learning-based methods designed to scale sequential decision making to high-dimensional settings \citep{bertsekas1996neuro,powell2007approximate, sutton2018reinforcement}.

\subsection{Reinforcement Learning}
Reinforcement learning (RL) is based on the same MDP framework but addresses settings in which the transition law, the one-stage cost, or both are not explicitly known and must instead be learned from data or interaction with the environment \citep{sutton2018reinforcement,watkins1992q,konda2000actor}. In the standard RL view, an agent observes the current state, selects an action according to a policy, receives feedback through a cost or reward, and uses this interaction to improve future decisions. A policy $\pi(a\mid s)$ specifies how actions are selected in each state, and the goal is to learn a policy that minimizes the expected cumulative cost. Figure~\ref{fig:sequential_frameworks_comparison}(c) summarizes this agent--environment interaction.

In general, RL methods can be grouped into three main classes. Value-based methods estimate quantities such as the value function $V_t(s)$ or the action-value function $Q_t(s,a)$ and make decisions from these estimates. Policy-based methods instead parameterize the policy directly and optimize it from the data. Actor--critic methods combine these two ideas by learning both a policy representation and a value-based performance signal.

Relative to classical DP, RL replaces exact model-based recursion with data-driven approximation, making it attractive for high-dimensional problems and environments that are difficult to model analytically. At the same time, standard RL formulations do not naturally enforce the rich feasibility constraints, integrality requirements, and explicit non-anticipativity conditions that arise in many operations research applications. Thus, while RL offers scalable tools for learning sequential decisions, additional modeling structure is often needed in constrained OR settings. In many contemporary applications, these ideas are often combined with deep neural networks to approximate value functions or policies, giving rise to deep reinforcement learning, which we discuss in later sections.

\subsection{Multi-Stage Stochastic Programming}

While MDPs, DP, and RL adopt a state-based view, multi-stage stochastic programming (MSP)
remains closer to the original history-based formulation in \eqref{eq:general_multistage_compact},
where decisions adapt to revealed information and are linked across stages through recourse
constraints \citep{birge2011introduction,shapiro2009lectures}. Representative algorithmic
and decomposition advances for multi-stage stochastic integer programming include
\citep{zou2019sddip,lulli2004branch,buyuktahtakin2022stage,daryalal2024lagrangian,
boland2018combining}. In MSP, this structure is represented explicitly over a finite set of scenarios, each corresponding to a realization of the uncertainty process.

In a scenario-based extensive formulation, stage-wise decisions are indexed by both time and scenario, say \(x_t^\omega\), and the objective minimizes expected total cost, typically of the form \(\sum_{\omega\in\Omega} p^\omega \sum_{t=1}^T c_t^{\omega\top}x_t^\omega\), subject to scenario-wise feasibility and non-anticipativity. The latter requires \(x_t^\omega=x_t^{\omega'}\) whenever scenarios \(\omega\) and \(\omega'\) share the same history up to stage \(t\). Figure~\ref{fig:sequential_frameworks_comparison}(d) illustrates this scenario-tree view, where branching captures the revelation of uncertainty and downstream decisions provide recourse.

Compared with the MDP framework, MSP places greater emphasis on explicit decision variables and constraint systems. This is especially valuable in OR/MS applications involving multi-scale temporal or spatial dynamics, network flows, capacity limits, logical conditions, budget restrictions, or mixed-integer decisions, where structural fidelity is often more important than a compact state representation. The trade-off is computational: scenario-based formulations grow rapidly with the number of stages and uncertainty realizations, and integer decisions make the resulting extensive forms even more challenging, motivating decomposition and scenario-based acceleration methods such as stage-\(t\) scenario dominance \citep{buyuktahtakin2022stage}.

\subsection{Perspective and Implications}

The frameworks described above provide complementary lenses for sequential decision making under uncertainty. The general multi-stage formulation in \eqref{eq:general_multistage_compact} starts from the broad history-based view: uncertainty is revealed over time, decisions adapt to available information, and recourse and non-anticipativity are explicit. MDPs and dynamic programming obtain analytical power by summarizing history through a state variable and exploiting Bellman recursion; reinforcement learning replaces explicit model knowledge with data-driven estimation of value functions or policies; and multi-stage stochastic programming preserves scenario-dependent decisions and constraints more explicitly. Recent work further bridges these views by incorporating MDP-style state transitions, decision-dependent uncertainty, and limited statistical learning into multi-stage stochastic programs through policy-graph representations \citep{morton2026mdp}.

From an OR/MS perspective, the key trade-off is between modeling fidelity and computational tractability. MDPs and DPs provide strong optimality principles, but require a manageable state representation. RL scales more naturally to complex and high-dimensional problems, but standard formulations do not automatically enforce feasibility, integrality, or scenario-tree consistency. MSP captures these features directly and is especially valuable in constrained, multi-scale temporal or spatial decision problems, though often at substantial computational cost. This tension between structure and scalability helps explain the growing interest in integrating learning and optimization. It also raises a natural next question: which neural architectures are best suited to represent sequential decisions, value functions, and policies in complex OR/MS settings?

\section{Neural Deep Learning Architectures from a Decision Making Lens}
\label{sec:architectures}

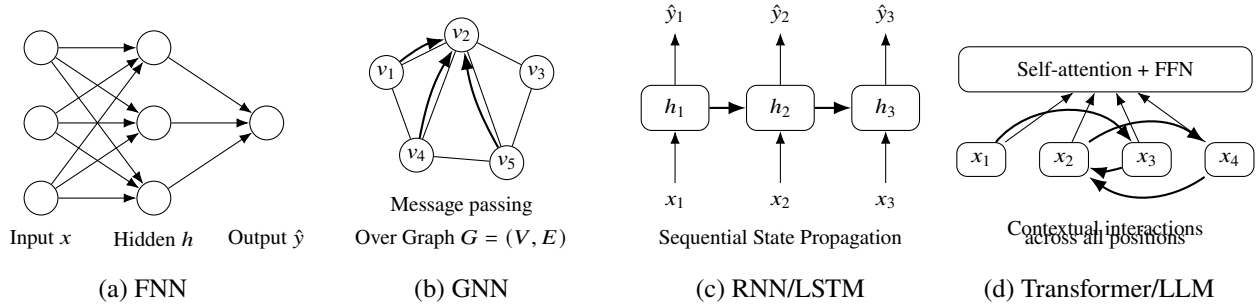
\begin{figure*}[h!]
\centering

\begin{subfigure}[t]{0.23\textwidth}
\centering
\begin{tikzpicture}[
    x=1cm,y=1cm,
    neuron/.style={circle, draw, minimum size=0.45cm, inner sep=0pt},
    >=Latex
]
\node[neuron] (i1) at (0,1.0) {};
\node[neuron] (i2) at (0,0.0) {};
\node[neuron] (i3) at (0,-1.0) {};

\node[neuron] (h1) at (1.5,1.0) {};
\node[neuron] (h2) at (1.5,0.0) {};
\node[neuron] (h3) at (1.5,-1.0) {};

\node[neuron] (o1) at (3.0,0.0) {};

\foreach \i in {i1,i2,i3}
    \foreach \h in {h1,h2,h3}
        \draw[->,thin] (\i) -- (\h);

\foreach \h in {h1,h2,h3}
    \draw[->,thin] (\h) -- (o1);

\node at (0,-1.55) {\scriptsize Input \(x\)};
\node at (1.5,-1.55) {\scriptsize Hidden \(h\)};
\node at (3.0,-1.55) {\scriptsize Output \(\hat y\)};
\end{tikzpicture}
\caption{FNN}
\end{subfigure}
\hfill
\begin{subfigure}[t]{0.23\textwidth}
\centering
\begin{tikzpicture}[
    x=1cm,y=1cm,
    nodeg/.style={circle, draw, minimum size=0.45cm, inner sep=0pt},
    >=Latex
]
\node[nodeg] (v1) at (0.2,0.9) {\scriptsize \(v_1\)};
\node[nodeg] (v2) at (1.2,1.4) {\scriptsize \(v_2\)};
\node[nodeg] (v3) at (2.2,0.9) {\scriptsize \(v_3\)};
\node[nodeg] (v4) at (0.6,-0.2) {\scriptsize \(v_4\)};
\node[nodeg] (v5) at (1.8,-0.3) {\scriptsize \(v_5\)};

\draw[-,thin] (v1) -- (v2);
\draw[-,thin] (v2) -- (v3);
\draw[-,thin] (v1) -- (v4);
\draw[-,thin] (v2) -- (v4);
\draw[-,thin] (v2) -- (v5);
\draw[-,thin] (v3) -- (v5);
\draw[-,thin] (v4) -- (v5);

\draw[->,thick] (v1) to[bend left=12] (v2);
\draw[->,thick] (v4) to[bend left=10] (v2);
\draw[->,thick] (v5) to[bend left=10] (v2);

\node[align=center] at (1.2,-1.1) {\scriptsize Message passing\\[-1pt] \scriptsize over graph \(G=(V,E)\)};
\end{tikzpicture}
\caption{GNN}
\end{subfigure}
\hfill
\begin{subfigure}[t]{0.23\textwidth}
\centering
\begin{tikzpicture}[
    x=1cm,y=1cm,
    block/.style={draw, rounded corners, minimum width=0.9cm, minimum height=0.6cm},
    >=Latex
]
\node[block] (r1) at (0.4,0) {\scriptsize \(h_1\)};
\node[block] (r2) at (1.8,0) {\scriptsize \(h_2\)};
\node[block] (r3) at (3.2,0) {\scriptsize \(h_3\)};

\draw[->,thick] (r1) -- (r2);
\draw[->,thick] (r2) -- (r3);

\draw[->,thin] (0.4,-1.0) -- (r1);
\draw[->,thin] (1.8,-1.0) -- (r2);
\draw[->,thin] (3.2,-1.0) -- (r3);

\draw[->,thin] (r1) -- (0.4,1.0);
\draw[->,thin] (r2) -- (1.8,1.0);
\draw[->,thin] (r3) -- (3.2,1.0);

\node at (0.4,-1.25) {\scriptsize \(x_1\)};
\node at (1.8,-1.25) {\scriptsize \(x_2\)};
\node at (3.2,-1.25) {\scriptsize \(x_3\)};

\node at (0.4,1.25) {\scriptsize \(\hat y_1\)};
\node at (1.8,1.25) {\scriptsize \(\hat y_2\)};
\node at (3.2,1.25) {\scriptsize \(\hat y_3\)};

\node[align=center] at (1.8,-1.75) {\scriptsize Sequential state propagation};
\end{tikzpicture}
\caption{RNN/LSTM}
\end{subfigure}
\hfill
\begin{subfigure}[t]{0.23\textwidth}
\centering
\begin{tikzpicture}[
    x=1cm,y=1cm,
    tok/.style={draw, rounded corners, minimum width=0.65cm, minimum height=0.45cm},
    >=Latex
]
\node[tok] (t1) at (0.2,0) {\scriptsize \(x_1\)};
\node[tok] (t2) at (1.3,0) {\scriptsize \(x_2\)};
\node[tok] (t3) at (2.4,0) {\scriptsize \(x_3\)};
\node[tok] (t4) at (3.5,0) {\scriptsize \(x_4\)};

\draw[->,thick] (t2) to[bend left=35] (t4);
\draw[->,thick] (t4) to[bend left=35] (t2);
\draw[->,thick] (t1) to[bend left=50] (t3);
\draw[->,thick] (t3) to[bend left=18] (t2);

\node[draw, rounded corners, minimum width=3.9cm, minimum height=0.6cm] (attn) at (1.85,1.2) {\scriptsize Self-attention + FFN};

\foreach \x in {t1,t2,t3,t4}
    \draw[->,thin] (\x) -- (attn);

\node[align=center] at (1.85,-1.0) {\scriptsize Contextual interactions\\[-1pt] \scriptsize across all positions};
\end{tikzpicture}
\caption{Transformer/LLM}
\end{subfigure}

\caption{Comparison of FNN, GNN, RNN/LSTM, and transformer architectures for learning-based decision systems.}
\label{fig:architecture_comparison}
\end{figure*}

This section summarizes four major families of neural architectures used in learning-based decision systems: feedforward neural networks (FNNs), graph neural networks (GNNs), recurrent architectures including recurrent neural networks (RNNs) and long short-term memory networks (LSTMs), and transformer/LLM architectures. These models are based on the foundational work in multilayer neural networks and representation learning \citep{rumelhart1986learning,lecun2015deep,goodfellow2016deep}. Our focus is not on neural-network implementation details, but on the type of decision structure each architecture naturally represents: static feature mappings, relational dependencies, temporal dynamics, and long-range contextual interactions. Throughout, \(x\) denotes an input, \(h\) a hidden representation, \(\hat y\) a model output or predicted decision, and \(\theta\) the set of trainable parameters. Figure~\ref{fig:architecture_comparison} provides a visual comparison of the architectures discussed in this section.

For OR/MS readers, a useful way to compare neural architectures is through their \emph{inductive bias}: the type of decision structure the architecture represents naturally. As summarized in Table~\ref{tab:architecture_choice}, FNNs, GNNs, recurrent models, and transformers differ in whether they primarily encode static, relational, temporal, or contextual dependencies. These architectures can be used to approximate optimization-relevant objects, including decisions, recourse values, value functions, and policies. A common source of these inductive biases is parameter sharing: FNNs apply a fixed nonlinear map to problem features; GNNs reuse update rules across nodes and edges; recurrent models reuse transition rules across time; and transformers reuse attention projections across positions. Architecture choice in learning-based decision systems is therefore guided not only by prediction accuracy, but also by how well the architecture matches the structure of the underlying optimization problem.

\begin{table}[h!]
\centering
\small
\caption{Neural architecture choice for OR/MS decision models.}
\label{tab:architecture_choice}
\begin{tabular}{p{0.18\textwidth}p{0.34\textwidth}p{0.36\textwidth}}
\hline
\textbf{Architecture} & \textbf{Natural inductive bias} & \textbf{Typical OR/MS use} \\
\hline
FNN & Static feature-to-output mappings
& Predicting costs, recourse values, decisions, or value-function approximations from fixed problem features \citep{elmachtoub2022smart,mnih2015human} \\
GNN & Relational structure among nodes, arcs, or interacting entities
& Routing, network design, supply chains, facility location, and graph-structured combinatorial optimization or solver guidance \citep{gasse2019learn2branch, khalil2022mipgnn} \\
RNN/LSTM & Temporal dependence and evolving hidden state
& Inventory control, lot sizing, dynamic resource allocation, sequential decision learning, and recurrent value or policy approximation under partial observability \citep{yilmaz2023lstm,hausknecht2015drqn} \\
Transformer & Long-range contextual interactions across sequences, sets or decisions
& Attention-based routing, long-horizon planning, and transformer-based MIP/production-planning solution prediction \citep{kool2019attention,buyuktahtakin2024transformers} \\
\hline
\end{tabular}
\end{table}

\subsection{Feedforward Neural Networks (FNNs)}

Feedforward neural networks are the canonical architecture for \emph{static} decision mappings. They are most appropriate when the output depends on a fixed feature vector rather than on temporal history, graph structure, or broader contextual interactions. This makes them a natural starting point for predict--then--optimize pipelines, surrogate optimization, and approximation of static decision rules \citep{rumelhart1986learning,hornik1989multilayer,goodfellow2016deep}.

An \(L\)-layer FNN maps an input \(x\in\mathbb{R}^{d_0}\) to an output \(\hat y\in\mathbb{R}^{d_L}\) through repeated affine transformations and nonlinear activations. Starting from \(h^{(0)}=x\), each layer computes \(h^{(\ell+1)}=\sigma^{(\ell)}(W^{(\ell)}h^{(\ell)}+b^{(\ell)})\) for \(\ell=0,\ldots,L-1\), with \(\hat y=h^{(L)}\). The trainable parameters are the weights and biases \(\theta=\{W^{(\ell)},b^{(\ell)}\}_{\ell=0}^{L-1}\). The key mechanism is repeated nonlinear feature transformation, allowing the network to construct progressively richer latent representations of the input.

In OR/MS applications, FNNs are useful for predicting optimal or near-optimal decisions, approximating objective values or recourse functions, learning dispatching or ranking rules, and embedding fast surrogate models inside larger optimization pipelines. Their main limitation is that they are memoryless: they do not naturally encode sequential dependence, relational structure, or long-range contextual interactions. These limitations motivate richer architectures such as GNNs, recurrent models, and transformers.

\subsection{Graph Neural Networks (GNNs)}

Many OR/MS problems are naturally defined on graphs, including transportation networks, supply chains, routing systems, facility-location models, and precedence-constrained scheduling problems. In these settings, the input is not simply a feature vector, but a graph \(G=(V,E)\) in which the nodes and edges encode relationships between entities. Graph neural networks are designed for such data by combining node features with graph structure through message passing \citep{scarselli2009graph,kipf2017semi,hamilton2017inductive,wu2021comprehensive}.

Let \(h_v^{(\ell)}\) denote the representation of node \(v\in V\) at layer \(\ell\), with \(h_v^{(0)}\) as the initial node feature vector. A generic message-passing layer updates each node by aggregating information from its neighbors:
\begin{equation}
h_v^{(\ell+1)}
=
\psi^{(\ell)}\!\left(
h_v^{(\ell)},
\mathrm{AGG}\bigl(\{h_u^{(\ell)}:u\in\mathcal{N}(v)\}\bigr)
\right),
\end{equation}
where \(\mathcal{N}(v)\) is the neighborhood of node \(v\), \(\mathrm{AGG}(\cdot)\) is a permutation-invariant aggregation operator such as sum, mean, or max, and \(\psi^{(\ell)}(\cdot)\) is a learnable update function. In OR/MS applications, edge features, such as costs, capacities, distances, travel times, or precedence relations, can also be incorporated into the message or update functions. Repeating the process allows each node representation to incorporate information from progressively larger graph neighborhoods.

For decision making, the appeal of GNNs is that they explicitly exploit relational structure. Depending on the task, the final representations can be used to produce node-level, edge-level, or graph-level predictions. This makes GNNs well aligned with classical OR/MS formulations in which topology, adjacency, flow, precedence, or neighborhood effects determine solution quality. Typical applications include routing, network design, facility location, scheduling with precedence relations, learning heuristics for graph-structured combinatorial optimization, and solver guidance for mixed-integer programming \citep{khalil2017learning,gasse2019learn2branch,khalil2022mipgnn,bengio2021machine,wu2021comprehensive}.

\subsection{Sequential Architectures: RNNs, LSTMs, and Temporal Convolutions}

Many OR/MS problems are inherently sequential: the current decision depends not only on current data, but also on past states, observations, or actions. Recurrent neural networks are designed for this setting by maintaining a hidden state that evolves over time and summarizes relevant history \citep{elman1990finding,rumelhart1986learning,goodfellow2016deep}.

For a sequence \(x_{1:T}=(x_1,\ldots,x_T)\), an RNN maintains a hidden state \(h_t\) and produces an output \(\hat y_t\) according to
\begin{equation}
h_t=\phi_\theta(x_t,h_{t-1}),\qquad
\hat y_t=\psi_\theta(h_t),
\quad t=1,\ldots,T.
\end{equation}
Thus, the architecture shifts from a static map \(x\mapsto\hat y\) to a dynamic map in which the prediction at time \(t\) depends on both the current input and the accumulated hidden state. For OR/MS readers, this is useful because many systems are sequential, having evolving inventories, queues, resource levels, demand histories, or operational states.

Standard RNNs can struggle with long-range dependence because gradients may vanish or explode over many recurrent steps \citep{bengio1994learning}. Long short-term memory networks address this limitation by introducing a memory cell and gates that regulate what information is retained, updated, and exposed through the hidden state \citep{hochreiter1997long,gers2000learning}. The key modeling advantage is that LSTMs can preserve information over longer horizons while selectively updating memory as new inputs arrive. This makes them useful in long-horizon or delayed-effect settings, such as inventory control, lot sizing, dynamic allocation, and sequential decision learning \citep{hochreiter1997long,goodfellow2016deep,yilmaz2023lstm}.

Temporal convolutional networks provide another way to model sequential data by applying one-dimensional filters across time. With causal or dilated convolutions, they can extract local and multi-scale temporal patterns in a parallelizable way \citep{vandenoord2016wavenet,bai2018empirical}. However, in this tutorial, the emphasis is on recurrent and attention-based architectures because they connect most directly to the sequential decision-learning frameworks discussed later.

For decision making applications, the central strength of recurrent architectures is history dependence. Their main limitation is that information is processed sequentially, which can restrict parallelization and make very long-context reasoning difficult. These limitations motivate attention-based architectures.

\subsection{Transformer/LLM Architectures and Attention}

Transformer architectures are a central framework for modern sequence modeling and form the basis of most large language models \citep{vaswani2017attention,brown2020language}. Their defining innovation is self-attention, which allows each element of a sequence to directly incorporate information from other relevant elements rather than relying on a recurrent hidden state. This makes transformers useful when decisions depend on long-range interactions, a variable-length context, heterogeneous input, or rich dependencies across many entities or time periods.

Let \(x_{1:n}=(x_1,\ldots,x_n)\) denote an input sequence. Each element is mapped to an embedding and combined with positional information, producing initial representations \(h_i^{(0)}\) for \(i=1,\ldots,n\). In OR/MS applications, positional or structural encodings should respect the symmetry of the problem. Sequence position is meaningful in time-indexed planning, where period order matters. In contrast, arbitrary item or facility orderings should not affect predictions in permutation-invariant problems such as facility location, set selection, or some routing formulations; otherwise, the model may incorrectly treat an arbitrary input order as meaningful, instead of learning from the actual decision structure of the problem. At layer \(\ell\), the sequence representations are stacked in \(H^{(\ell)}\), and self-attention forms query, key, and value matrices \(Q\), \(K\), and \(V\). The attention operation is
\begin{equation}
\mathrm{Attn}(Q,K,V)
=
\mathrm{softmax}\!\left(\frac{QK^\top}{\sqrt{d_k}}\right)V.
\end{equation}
The interpretation is straightforward: each position compares itself with other positions, assigns larger weights to the most relevant ones, and updates its representation using a weighted combination of their value vectors. Unlike recurrent models, information does not need to travel step by step through time, which helps transformers to capture long-range dependencies more directly \citep{vaswani2017attention}.

\begin{figure}[h!]
\centering
\resizebox{0.8\textwidth}{!}{%
\begin{tikzpicture}[
    x=1cm,y=1cm,
    >=Latex,
    token/.style={draw, rounded corners, minimum width=0.9cm, minimum height=0.55cm, align=center},
    mat/.style={draw, rounded corners, minimum width=1.2cm, minimum height=0.7cm, align=center, fill=gray!8},
    op/.style={draw, rounded corners, minimum width=1.8cm, minimum height=0.8cm, align=center, fill=blue!5},
    note/.style={align=center, font=\small}
]

\node[token] (x1) at (0.4,0) {$h_1^{(\ell)}$};
\node[token] (x2) at (1.8,0) {$h_2^{(\ell)}$};
\node[token] (x3) at (3.4,0) {$\cdots$};
\node[token] (x4) at (4.6,0) {$h_n^{(\ell)}$};

\node[note] at (2.1,-0.9) {Input sequence representations\\at layer \(\ell\)};

\node[mat] (Q) at (0.4,2.0) {$Q$};
\node[mat] (K) at (2.5,2.0) {$K$};
\node[mat] (V) at (4.3,2.0) {$V$};

\node[note] at (1.4,2.8) {\(Q=H^{(\ell)}W_Q\)};
\node[note] at (3.4,2.8) {\(K=H^{(\ell)}W_K\)};
\node[note] at (5.4,2.8) {\(V=H^{(\ell)}W_V\)};

\foreach \x in {x1,x2,x3,x4} {
    \draw[->] (\x.north) -- (Q.south);
    \draw[->] (\x.north) -- (K.south);
    \draw[->] (\x.north) -- (V.south);
}

\node[op] (score) at (2.1,4.1) {Similarity scores\\ \(QK^\top / \sqrt{d_k}\)};
\draw[->] (Q.north) -- (score.south west);
\draw[->] (K.north) -- (score.south);

\node[op] (soft) at (2.1,5.6) {Attention weights\\ \(\mathrm{softmax}\!\left(QK^\top/\sqrt{d_k}\right)\)};
\draw[->] (score.north) -- (soft.south);

\node[op] (attn) at (2.1,7.2) {Contextual aggregation\\ \(\mathrm{Attn}(Q,K,V)\)};
\draw[->] (soft.north) -- (attn.south);
\draw[->] (V.north) to[out=90,in=0] (attn.east);

\node[token] (y1) at (0,9.0) {$h_1^{(\ell+1)}$};
\node[token] (y2) at (1.4,9.0) {$h_2^{(\ell+1)}$};
\node[token] (y3) at (2.8,9.0) {$\cdots$};
\node[token] (y4) at (4.2,9.0) {$h_n^{(\ell+1)}$};

\draw[->] (attn.north) -- (2.1,8.2);
\draw[->] (2.1,8.2) -| (y1.south);
\draw[->] (2.1,8.2) -- (y2.south);
\draw[->] (2.1,8.2) -- (y3.south);
\draw[->] (2.1,8.2) -| (y4.south);

\node[note] at (2.1,9.8) {Updated sequence representations};

\node[draw, rounded corners, align=left, font=\small, fill=green!5] (insight) at (9.5,5.2) {%
\textbf{Attention mechanism:}\\
Each position \(i\)\\
1. forms a query from \(h_i^{(\ell)}\),\\
2. compares it with all keys,\\
3. assigns larger weights to more relevant positions,\\
4. updates its representation using a weighted\\
combination of the values.};

\end{tikzpicture}
}
\caption{Illustration of the self-attention mechanism. The layer-\(\ell\) token representations \(\{h_i^{(\ell)}\}_{i=1}^n\) are stacked into the matrix \(H^{(\ell)}\), from which the query, key, and value matrices are formed. Pairwise similarity scores are then computed, normalized into attention weights, and used to produce updated contextual representations \(\{h_i^{(\ell+1)}\}_{i=1}^n\). Adapted from \citet{vaswani2017attention}.}
\label{fig:self_attention_notation}
\end{figure}
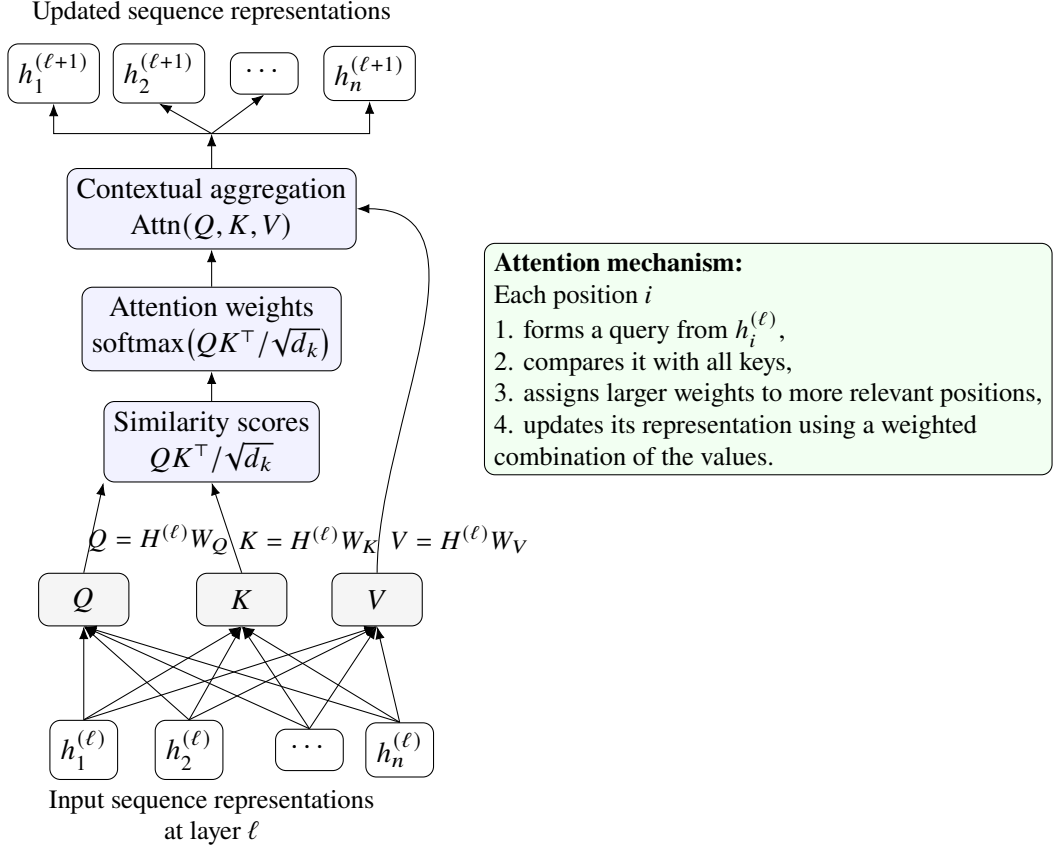

Figure~\ref{fig:self_attention_notation} summarizes this self-attention mechanism, showing how relevance scores across positions are converted into attention weights and used to form updated contextual representations. In practice, transformer layers use multihead attention, feedforward sublayers, residual connections, and normalization to capture several types of relationships and stabilize training \citep{vaswani2017attention}. Transformers appear in encoder-only, encoder--decoder, and decoder-only forms. Encoder-only models learn contextualized representations, encoder--decoder models map input sequences to output sequences, and decoder-only models generate outputs autoregressively and form the basis of modern LLMs \citep{vaswani2017attention,brown2020language}.

For OR/MS decision making, the main advantage of transformers is their ability to model rich contextual interactions without compressing all prior information into a single recurrent state. This is especially useful when decisions depend on long-range dependencies, heterogeneous entities, structured text, or complex constraints embedded in sequential context. Compared to recurrent models, transformers are more parallelizable and often more effective in capturing global nonlinear structure, making them increasingly relevant for learning-based optimization and decision support \citep{bengio2021machine,kool2019attention,buyuktahtakin2024transformers}.

\subsection{Design Considerations for OR/MS Applications}

In OR/MS applications, architecture choice matters when it reflects the structure of
the decision problem. Returning to the question posed at the end of Section 2 — which architectures are best suited to represent sequential decisions, value functions, and policies — Table 1 summarizes the answer at the level of inductive bias. FNNs provide a natural starting point for static decision mappings and
value-function approximation, as in Deep Q-Networks \citep{mnih2015human}; GNNs extend
this role to relational and graph-structured states and actions, as in graph-based
reinforcement learning for combinatorial optimization \citep{khalil2017learning}; RNNs
and LSTMs, including Deep Recurrent Q-Networks \citep{hausknecht2015drqn}, represent
value functions and policies that must condition on an evolving history over a planning
horizon; and transformers \citep{chen2021decision,buyuktahtakin2024transformers}
support policies and decision representations with long-range contextual dependencies.
Learned models often require optimization layers, constraint-aware design, feasibility
repair, warm starts, or solver-based refinement because feasibility and
implementability are not automatic \citep{amos2017optnet,agrawal2019differentiable,
yilmaz2024expandable}. Generalization across horizons, scenarios, and problem scales
is also difficult, particularly in sequential settings where errors propagate over
time. Deep learning is therefore most valuable when it captures recurring structure
across related instances, temporal patterns, network interactions, or contextual
dependencies, while optimization preserves feasibility, recourse, and decision
quality. This perspective motivates the learning--optimization integration paradigms
discussed next.

\section{Learning--Optimization Integration: Paradigms}
\label{sec:integration}

Learning--optimization integration can be organized by the role learning plays in the decision pipeline. As shown in Figure~\ref{fig:learning_optimization_paradigms}, predict-then-optimize (Figure~\ref{fig:learning_optimization_paradigms}(a)) estimates uncertain inputs before optimization; decision-aware learning (Figure~\ref{fig:learning_optimization_paradigms}(b)) trains models using downstream decision quality; and learning-to-optimize (Figure~\ref{fig:learning_optimization_paradigms}(c)) leverages solved instances, optimization traces, or related problem data to generate, accelerate, or adapt optimization decisions. A separate predict-then-optimize pipeline may suffice when forecasts are reliable and the optimization model is stable, whereas tighter integration becomes valuable when prediction errors change decisions, similar instances are solved repeatedly, or computational tractability is central. Broader perspectives are provided by recent surveys on machine learning for combinatorial optimization, end-to-end constrained optimization learning, and optimization under uncertainty \citep{bengio2021machine,kotary2021survey,sadana2025survey}.

\begin{figure*}[h!]
\centering
\resizebox{0.92\textwidth}{!}{
\begin{tikzpicture}[
    x=1cm, y=1cm,
    >=Latex,
    font=\small,
    modelbox/.style={
        rectangle, draw, rounded corners,
        minimum width=3.2cm, minimum height=0.95cm,
        align=center, very thick, fill=blue!12
    },
    optbox/.style={
        rectangle, draw, rounded corners,
        minimum width=3.2cm, minimum height=0.95cm,
        align=center, very thick, fill=green!14
    },
    outbox/.style={
        rectangle, draw, rounded corners,
        minimum width=3.2cm, minimum height=0.95cm,
        align=center, very thick, fill=blue!7
    },
    dashedbox/.style={
        rectangle, draw, dashed, rounded corners,
        minimum width=2.8cm, minimum height=0.9cm,
        align=center, very thick
    },
    fadedbox/.style={
        rectangle, draw=gray!55, dashed, rounded corners,
        minimum width=2.6cm, minimum height=0.85cm,
        align=center, very thick, fill=gray!10, text=gray!60
    },
    title/.style={font=\bfseries},
    arrow/.style={->, very thick},
    fbarrow/.style={->, line width=1.1pt, dashed, draw=orange!80!black},
    note/.style={font=\scriptsize, align=center}
]
\node[title] at (0,5.6)    {(a) Predict-then-Optimize};
\node[title] at (5.8,5.6)  {(b) Decision-Aware Learning};
\node[title] at (11.3,5.6) {(c) Learning-to-Optimize};
\node[dashedbox] (a0) at (0,4.3)   {Contextual\\ Features};
\node[modelbox]  (a1) at (0,2.75)  {Predictive Model};
\node[dashedbox] (a2) at (0,1.2)   {$\hat{\xi}$};
\node[optbox]    (a3) at (0,-0.35) {Optimization};
\node[dashedbox] (a4) at (0,-1.9)  {Decision};
\draw[arrow] (a0) -- (a1);
\draw[arrow] (a1) -- (a2);
\draw[arrow] (a2) -- (a3);
\draw[arrow] (a3) -- (a4);
\node[dashedbox] (b0) at (5.8,4.3)   {Contextual\\ Features};
\node[modelbox]  (b1) at (5.8,2.75)  {Learning Model};
\node[dashedbox] (b2) at (5.8,1.2)   {$\hat{\xi}$};
\node[optbox]    (b3) at (5.8,-0.35) {Optimization};
\node[dashedbox] (b4) at (5.8,-1.9)  {Decision};
\draw[arrow] (b0) -- (b1);
\draw[arrow] (b1) -- (b2);
\draw[arrow] (b2) -- (b3);
\draw[arrow] (b3) -- (b4);
\draw[fbarrow]
    (b3.east) .. controls +(2.1,0.42) and +(2.1,-0.42) .. (b1.east);
\node[note, text=orange!80!black, anchor=center] at (8.1,1.2)
    {Decision loss\\ (gradient)};
\node[dashedbox] (c0)  at (11.3,4.3)   {Problem Instance\\ $\alpha$};
\node[modelbox]  (c1)  at (11.3,2.75)  {Learning Model};
\node[fadedbox]  (cxi) at (11.3,1.2)   {$\hat{\xi}$};
\node[outbox]    (c2)  at (11.3,-0.35) {Candidate Decision\\ or Solver Guidance};
\node[dashedbox] (c3)  at (11.3,-1.9)  {Decision};
\draw[gray!60, thick] (cxi.south west) -- (cxi.north east);
\node[note, text=gray!60, anchor=west] at (12.75,1.2) {(bypassed)};
\draw[arrow] (c0) -- (c1);
\draw[arrow] (c1.east) .. controls +(1.55,-0.25) and +(1.55,0.25) .. (c2.east);
\draw[arrow] (c2) -- (c3);
\node[modelbox, minimum width=0.55cm, minimum height=0.38cm, label=right:{\scriptsize Learned component}]
    at (2.1,-3.5) {};
\node[optbox, minimum width=0.55cm, minimum height=0.38cm, label=right:{\scriptsize Optimization}]
    at (6.4,-3.5) {};
\node[dashedbox, minimum width=0.55cm, minimum height=0.38cm, label=right:{\scriptsize Data / decision}]
    at (10.0,-3.5) {};
\end{tikzpicture}}
\caption{Comparison of learning--optimization paradigms. All three pipelines begin from a problem input and end in a decision; color denotes the role of each block (learned, optimization, or data/decision). (a) Predict-then-optimize separates prediction and optimization: a model produces predicted parameters $\hat{\xi}$ that are then optimized over. (b) Decision-aware learning keeps optimization as a separate block but trains the model using gradients of a downstream decision loss (dashed feedback loop). (c) Learning-to-optimize maps a problem instance directly to a candidate decision or solver guidance, which occupies the optimization stage and thereby subsumes (direct decision maps) or accelerates (solver guidance) explicit optimization; the parameter-prediction step $\hat{\xi}$ is bypassed (struck through).}
\label{fig:learning_optimization_paradigms}
\end{figure*}

\subsection{Predict-Then-Optimize}

The predict-then-optimize paradigm remains the dominant approach in many applications. In this framework, observed covariates or contextual information \(z\), such as historical demand, prices, weather, customer features, or current system descriptors, are first used to predict uncertain parameters \(\xi\). The predicted parameters are then passed to an optimization model:
\begin{equation}
    \hat{\xi} = f_{\theta}(z), \quad
    x^*(\hat{\xi}) = \arg\min_{x \in \mathcal{X}(\hat{\xi})} c(x, \hat{\xi}).
\end{equation}
This paradigm aligns with classical OR/MS workflows, where forecasts or scenarios are generated and then used within optimization models \citep{shapiro2009lectures,birge2011introduction}.

Despite its simplicity and modularity, predict-then-optimize can suffer from a misalignment between predictive accuracy and decision quality. Predictive models are typically trained by minimizing statistical loss functions such as mean squared error,
\begin{equation}
    \min_{\theta} \mathbb{E}\left[\|f_{\theta}(z) - \xi\|^2\right],
\end{equation}
which does not necessarily minimize the downstream decision loss,
\begin{equation}
    \mathbb{E}\left[c(x^*(f_{\theta}(z)), \xi)\right].
\end{equation}
This issue is central to decision-focused learning and related predict-and-optimize approaches, where small prediction errors can sometimes induce large decision errors in constrained or combinatorial optimization problems \citep{elmachtoub2022smart}.

Moreover, predict-then-optimize treats prediction and optimization as separate steps, so it may not account for how prediction errors affect the final decision. In constrained or combinatorial problems, even small errors can change the optimizer's solution, reflecting the instability that can arise in parametric optimization \citep{rockafellar2009variational}. As a result, the approach can produce suboptimal decisions and, when feasibility depends on uncertain parameters, decisions that may prove infeasible under realized conditions.

\subsection{Decision-Aware Learning}

Decision-aware learning addresses the misalignment between prediction accuracy and decision quality by incorporating the optimization problem directly into the learning objective. Instead of minimizing prediction error alone, the model is trained to minimize expected decision loss:
\begin{equation}
    \min_{\theta} \; \mathbb{E}_{s,\xi} \left[c\left(x^*(f_{\theta}(s)), \xi\right)\right].
\end{equation}

This paradigm often relies on differentiating through the optimization problem or constructing surrogate loss functions that approximate the decision objective. The foundational contributions include \citet{donti2017task}, who introduce task-based learning for stochastic optimization, and \citet{elmachtoub2022smart}, who propose the SPO+ loss for structured prediction, where SPO stands for Smart Predict-then-Optimize. The SPO+ loss is a convex surrogate for the decision loss induced by using predicted parameters in a downstream optimization problem.

More broadly, this approach is closely related to prescriptive analytics \citep{bertsimas2020predictive}, where predictive models are explicitly designed to support decision making. Recent advances in differentiable optimization \citep{amos2017optnet,agrawal2019differentiable} have further enabled the end-to-end integration of optimization layers within neural networks.

However, decision-aware learning faces significant challenges in sequential and combinatorial settings. In multi-stage stochastic programming, scenario-indexed decisions must satisfy non-anticipativity constraints:
\begin{equation}
    x_t^{\omega} = x_t^{\omega'}
    \quad \text{whenever} \quad
    \xi_{[t]}^{\omega} = \xi_{[t]}^{\omega'}.
\end{equation}
That is, two scenarios with the same information history up to stage \(t\) must receive the same stage \(t\) decision. Such requirements are difficult to handle with standard differentiable optimization layers, which are most natural for continuous and differentiable problem classes. The challenge becomes more pronounced when the underlying optimization model includes integer decisions, logical constraints, recourse decisions, or nonconvex feasible regions. These limitations motivate alternative paradigms that learn the decision structure while relying on optimization-based screening, repair, or refinement to recover feasible and implementable solutions.

\subsection{Learning-to-Optimize}

Learning-to-optimize (L2O) leverages solved instances, solver information, or related problem data to improve how optimization problems are solved or approximated.  In this tutorial, we organize L2O methods into three broad roles: learning to generate  high-quality decisions directly, learning to accelerate optimization algorithms, and learning to adapt or embed surrogate optimization models. This organization is broadly consistent with the recent overview of \citet{chen2024learning}.

The first paradigm aims to learn the \emph{solution operator} itself. Given an optimization instance \(\alpha\), the model returns an optimal decision
\begin{equation}
    x^*(\alpha)=\arg\min_{x\in\mathcal{X}(\alpha)} c(x,\alpha),
\end{equation}
where \(\mathcal{X}(\alpha)\) is the feasible set and \(c(x,\alpha)\) is the objective function. Rather than solving this problem from scratch for every new instance, one trains a model
\begin{equation}
    \hat{x}(\alpha)=f_\theta(\alpha)
\end{equation}
to map instance features directly to optimal or near-optimal decisions. This setting is particularly attractive in OR/MS applications where closely related problems must be solved repeatedly under changing data, such as online, parametric, and real-time decision environments \citep{bengio2021machine,yilmaz2023lstm,yilmaz2024expandable}.

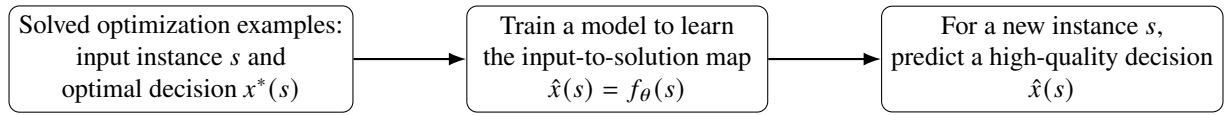
\begin{figure}[h!]
\centering
\begin{tikzpicture}[
    >=Latex,
    font=\small,
    node distance=1.5cm,
    box/.style={
        draw,
        rounded corners,
        align=center,
        minimum height=1.2cm,
        minimum width=4.0cm
    }
]

\node[box] (data) {Solved optimization examples:\\input instance \(\alpha\) and\\optimal decision \(x^*(\alpha)\)};

\node[box, right=of data] (train) {Train a model to learn\\the input-to-solution map\\\(\hat{x}(\alpha)=f_\theta(\alpha)\)};

\node[box, right=of train] (predict) {For a new instance \(\alpha\),\\predict a high-quality decision\\\(\hat{x}(\alpha)\)};

\draw[->, thick] (data) -- (train);
\draw[->, thick] (train) -- (predict);

\end{tikzpicture}
\caption{High-level learning-to-optimize pipeline for direct solution generation. Solved optimization instances are used to train a model that maps a problem input directly to an optimal or near-optimal decision, enabling fast prediction for new instances.}
\label{fig:l2o_simple}
\end{figure}

As illustrated in Figure~\ref{fig:l2o_simple}, the workflow is canonical: solved instances \((\alpha,x^*(\alpha))\) are generated offline, these input--solution pairs are used to train a predictive model, and the trained model is then applied to a new instance \(\alpha\) to produce a high-quality decision \(\hat{x}(\alpha)\). From an OR/MS perspective, part of the optimization effort is shifted offline so that good decisions can be generated more quickly for future related instances.

The second paradigm uses learning to accelerate, rather than replace, optimization. Here, the solver remains central, but the selected components are learned from the data to improve performance. Examples include learning branching rules, node selection strategies, cutting-plane policies, primal heuristics, warm starts, and decomposition components such as column generation or pricing
in combinatorial optimization \citep{khalil2016branch,lodi2017branching,bengio2021machine,khalil2017learning, hijazi2026improving}. This paradigm is especially appealing in OR/MS because the learned component guides the search without using machine learning as a proxy for the optimization model itself; the feasibility, constraints, and objective evaluation remain governed by the underlying solver. In this way, learning can improve practical scalability while preserving the structure and reliability of classical optimization algorithms.

The third paradigm, which we refer to as \emph{surrogate-embedded optimization}, adapts the optimization model itself by incorporating learned or predicted components. In many real-world applications, the nominal optimization formulation is only an approximation of the underlying decision environment. Learning or predictive modeling can therefore be used to construct surrogate objectives, approximate value or recourse functions, or embed predictive models directly within the optimization problem:
\begin{equation}
    \min_{x\in\mathcal{X}} \; \hat{c}_\theta(x)
    \qquad \text{or} \qquad
    \min_{x\in\mathcal{X}} \; c(x)+\hat{Q}_\theta(x),
\end{equation}
where \(\hat{c}_\theta(x)\) may approximate an unknown or expensive objective and \(\hat{Q}_\theta(x)\) may approximate a downstream value, recourse, or continuation function. More generally, a learned predictive model can also be embedded through auxiliary variables and constraints, for example,
\begin{equation}
    y=f_\theta(x), \qquad g(x,y)\le 0, \qquad x\in\mathcal{X},
\end{equation}
so that the prediction \(y\) influences the feasibility, cost, or performance of the downstream system. For OR/MS, this paradigm is especially important when data can make optimization models more realistic, adaptive, and decision-relevant \citep{bertsimas2020predictive,fischetti2018deep,anderson2020strong,misic2020treeensembles, lee2025transformer}. Recent work on two-stage stochastic programming follows this same logic by approximating the expected second-stage value function with a neural network and embedding the resulting surrogate within a tractable optimization model \citep{patel2022neur2sp}. Related work on optimization over trained neural networks further strengthens this perspective by studying how learned predictive models can be embedded and optimized directly within mathematical programming formulations \citep{tong2024ottnn}. More recently, optimization-based formulations of modern learning architectures, including transformer attention mechanisms represented within mixed-integer nonlinear programming models, further extend this paradigm by bringing the predictive structure directly into the optimization layer \citep{lee2025transformer}.

Related model-based optimization traditions have long embedded domain dynamics directly within decision models, including epidemic--logistics optimization for Ebola \citep{buyuktahtakin2018ebola}, invasive-species control and surveillance models that couple logistic biological growth, compartmental models, and spatial spread with optimization \citep{buyuktahtakin2011dynamic,kibis2021eab}, and cancer treatment-planning models that integrate Gompertz tumor-growth dynamics with multimodal treatment optimization \citep{kibis2019cancer}. Although these examples are not deep learning models, they illustrate the broader OR/MS principle that mechanistic models can predict system behavior within an optimization framework, thereby improving decision relevance.

Although these paradigms are analytically distinct, they often overlap in practice: a learned decision predictor may be followed by feasibility repair, learned components may guide exact algorithms, and surrogate models may be embedded within larger optimization pipelines. The next part of the tutorial focuses on learning-based decision generation for sequential settings, where mappings from problem instances to decisions must remain feasible, implementable, and consistent with the information available over time. Section~\ref{sec:learning_under_constraints} develops this perspective through recurrent, expandable, and non-anticipative learning--optimization architectures, while Section~\ref{sec:drl} examines deep reinforcement learning as a policy-learning approach for sequential and combinatorial decision problems.

\section{Learning Optimal Decisions under Constraints}
\label{sec:learning_under_constraints}

This section focuses on learning-based decision generation for sequential OR/MS problems with discrete, continuous, or mixed decisions, where temporal coupling, recourse, and implementability are central. Differentiable optimization layers, decision-focused and predict-and-optimize methods, neural combinatorial optimization, and neural approaches to stochastic programming provide important related perspectives \citep{amos2017optnet,agrawal2019differentiable,donti2017task,elmachtoub2022smart,wilder2019decision,mandi2024decision,vinyals2015pointer,kool2019attention,NazariEtAl2018,patel2022neur2sp}. However, once integer variables, logical constraints, recourse, and scenario coupling are introduced, the resulting optimization problems are generally nonconvex or nondifferentiable, so standard gradient-based end-to-end methods are difficult to apply directly. The setting emphasized here is computationally demanding: the deterministic mixed-integer core is already NP-hard, and the multi-stage extension layers temporal coupling, recourse, and non-anticipativity onto it. Broader surveys situate these developments within the machine learning literature for combinatorial optimization and end-to-end constrained optimization \citep{bengio2021machine,kotary2021survey,sadana2025survey}.

The frameworks below follow a learning--optimization route: they learn mappings from problem data to high-quality decisions using solved optimization instances, and then rely on validation, repair, warm-starting, partial fixing, or solver-based refinement to recover implementable solutions. This design is useful because constructive neural heuristics for routing and sequencing do not directly extend to multi-stage formulations with recourse, coupling constraints, and scenario-dependent information. We examine three architectures: an LSTM-based learner for temporal coupling \citep{yilmaz2023lstm}; an expandable encoder--decoder framework for generalization across horizons and dimensions without retraining \citep{yilmaz2024expandable}; and a non-anticipative learning--optimization framework for scenario-tree information structure in multi-stage stochastic programming \citep{yilmaz2025nonanticipative}. Together, they illustrate how deep learning can support sequential optimization while preserving the central role of OR/MS models and solvers.

\paragraph{Guiding example.}
Consider a capacitated lot-sizing problem over periods \(t=1,\ldots,T\). In each period, the decision maker chooses a binary setup decision \(u_t\), production quantity \(q_t\), and inventory level \(I_t\), forming the stage decision \(x_t=(u_t,q_t,I_t)\). Constraints such as \(q_t\le C_tu_t\) and \(I_t=I_{t-1}+q_t-d_t\) illustrate the mixed-integer and temporally coupled structure of the problem: current setup and production choices affect future inventories, feasibility, and capacity use. Under demand uncertainty, the model becomes a multi-stage stochastic program, where scenarios with the same information history \(\xi_{[t]}\) available at stage \(t\) must receive the same stage-\(t\) decision. This is the non-anticipativity requirement introduced in Section~\ref{sec:preliminaries}.
The learning frameworks below approximate the decision
mappings discussed in Sections~\ref{sec:architectures} and~\ref{sec:integration}; see Figure~\ref{fig:learning_optimization_paradigms}(c). The training and prediction process can be summarized as follows. \emph{Training data:} previously solved optimization instances provide pairs of input sequences and target decisions. For a horizon of length \(T\), a training instance consists of period-level inputs \((z_1,\ldots,z_T)\) and target decisions \((x_1^*,\ldots,x_T^*)\), where \(x_t^*\) is obtained from an optimal or near-optimal solution of that instance. In capacitated lot sizing, for example, a period-level input---the stage-dependent data \(\psi_t\)---may include production cost \(\rho_t=2\), setup cost \(f_t=100\), capacity \(C_t=40\), demand \(d_t=30\), and initial inventory \(I_0=10\); the solved instance may provide the corresponding optimal setup label \(u_t^*=1\). \emph{Input representation:} at stage \(t\), the neural input is \(z_t = (\xi_{[t]}, \psi_t, \hat{x}_{t-1})\). Here, \(\xi_{[t]}\) denotes the information history observed up to stage \(t\), \(\psi_t\) denotes stage-dependent optimization data such as costs, capacities, right-hand sides, and feasible-set information, and \(\hat{x}_{t-1}\) denotes the previous predicted decision when such feedback is used. \emph{Model training:} the network parameters \(\theta\) are fitted so that \(f_\theta(z_t)\) approximates \(x_t^*\), or a relevant component of it. In LSTM-Opt for capacitated lot sizing, the primary learned output is the binary setup decision, and the model is trained using the binary cross-entropy loss
\begin{equation}
-\frac{1}{T}\sum_{t=1}^{T}\left[u_t^*\log(\hat{u}_t)+(1-u_t^*)\log(1-\hat{u}_t)\right],
\label{eq:binary_cross_entropy_lstm}
\end{equation}
shown here for a single instance, with the empirical objective averaging this loss over the training set. Here, \(\hat{u}_t=[f_\theta(z_t)]_u\in[0,1]\) is the probability of the predicted setup decision. \emph{Prediction and optimization recovery:} at prediction time, the trained model maps \((z_1,\ldots,z_T)\) to predicted decisions \((\hat{x}_1,\ldots,\hat{x}_T)\), or to selected decision components. For binary setup decisions, the predicted probability is converted into a binary value using a threshold, for example \(\tilde{u}_t=1\) if \(\hat{u}_t\ge 0.5\) and \(\tilde{u}_t=0\) otherwise. Production and inventory decisions are then recovered by fixing the predicted binary setup variables in the MIP and solving the resulting linear program. Learning only binary setup decisions, rather than all components of \(x_t\), limits the training objective to a single-output task; jointly predicting setups, production quantities, and inventories couples heterogeneous targets into a multi-objective loss that can degrade prediction quality~\citep{yilmaz2023lstm}. Validation, repair, warm-starting, partial fixing, or solver-based refinement then help recover feasible, implementable solutions.

\subsection{Sequential Decision Learning via LSTM Frameworks}

A natural first instantiation of direct decision learning arises when the optimal solution itself is a sequence. In such settings, the task is not to predict a single static quantity but to learn a mapping from time-indexed problem data to a time-indexed sequence of decisions. This perspective is especially relevant in the field of OR/MS, where decisions are linked across time through inventories, capacities, budgets, setups, and uncertainty realizations.

LSTM architectures are well suited for this sequence-to-sequence setting because their recurrent memory mechanism allows the prediction at stage \(t\) to depend on both current input data and information propagated through the sequence \citep{hochreiter1997long}. In multi-stage optimization, dependence arises not only through the observed data, such as demands, costs, capacities, or realized uncertainty, but also through the decisions themselves. In capacitated lot sizing, for example, a setup or production decision in one period affects future inventories, capacity use, and subsequent production choices. Through its recurrent state, an LSTM can learn such temporal relationships so that the predicted sequence \(\hat{x}_1,\ldots,\hat{x}_T\) reflects both the evolution of the input data and the internal logic induced by constraints and the objective function \citep{yilmaz2023lstm}.

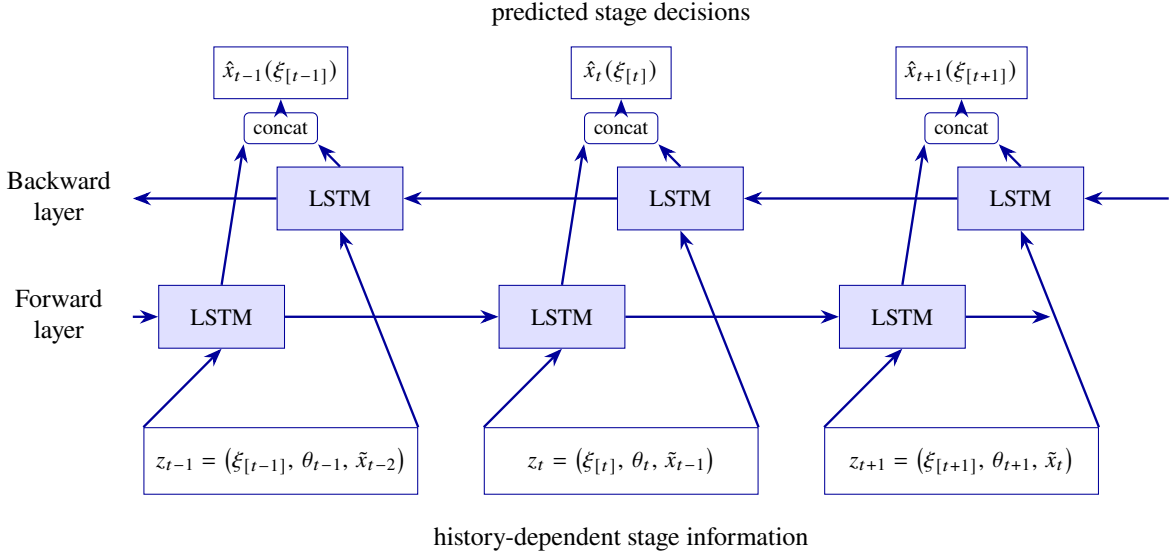
\begin{figure}[h!]
\centering
\resizebox{0.95\textwidth}{!}{
\begin{tikzpicture}[
    >=Stealth,
    font=\footnotesize,
    line/.style={draw=blue!60!black, line width=0.9pt},
    box/.style={draw=blue!60!black, fill=blue!12, rectangle,
                minimum width=1.7cm, minimum height=0.85cm, align=center},
    outbox/.style={draw=blue!60!black, fill=white, rectangle,
                   minimum width=1.35cm, minimum height=0.7cm, align=center},
    inobox/.style={draw=blue!60!black, fill=white, rectangle,
                   minimum width=3.7cm, minimum height=0.9cm, align=center},
    comb/.style={draw=blue!60!black, fill=white, rounded corners=2pt,
                 minimum width=0.8cm, minimum height=0.38cm, align=center, font=\scriptsize},
    lab/.style={font=\small}
]

\def\s{4.6}

\node[box] (f1) at (0,0) {LSTM};
\node[box] (f2) at (\s,0) {LSTM};
\node[box] (f3) at ({2*\s},0) {LSTM};

\node[box] (b1) at (1.6,1.6) {LSTM};
\node[box] (b2) at ({1.6+\s},1.6) {LSTM};
\node[box] (b3) at ({1.6+2*\s},1.6) {LSTM};

\node[comb] (c1) at (0.8,2.55) {concat};
\node[comb] (c2) at ({0.8+\s},2.55) {concat};
\node[comb] (c3) at ({0.8+2*\s},2.55) {concat};

\node[outbox] (x1) at (0.8,3.3) {$\hat{x}_{t-1}(\xi_{[t-1]})$};
\node[outbox] (x2) at ({0.8+\s},3.3) {$\hat{x}_{t}(\xi_{[t]})$};
\node[outbox] (x3) at ({0.8+2*\s},3.3) {$\hat{x}_{t+1}(\xi_{[t+1]})$};

\node[inobox] (z1) at (0.8,-1.95)
{$z_{t-1}=\big(\xi_{[t-1]},\,\psi_{t-1},\,\hat{x}_{t-2}\big)$};

\node[inobox] (z2) at ({0.8+\s},-1.95)
{$z_t=\big(\xi_{[t]},\,\psi_t,\,\hat{x}_{t-1}\big)$};

\node[inobox] (z3) at ({0.8+2*\s},-1.95)
{$z_{t+1}=\big(\xi_{[t+1]},\,\psi_{t+1},\,\hat{x}_{t}\big)$};

\draw[line,->] (-1.2,0) -- (f1.west);
\draw[line,->] (f1.east) -- (f2.west);
\draw[line,->] (f2.east) -- (f3.west);
\draw[line,->] (f3.east) -- ({2*\s+2.0},0);

\draw[line,->] ({1.6+2*\s+2.0},1.6) -- (b3.east);
\draw[line,->] (b3.west) -- (b2.east);
\draw[line,->] (b2.west) -- (b1.east);
\draw[line,->] (b1.west) -- (-1.2,1.6);

\draw[line,->] (z1.north west) -- (f1.south);
\draw[line,->] (z1.north east) -- (b1.south);

\draw[line,->] (z2.north west) -- (f2.south);
\draw[line,->] (z2.north east) -- (b2.south);

\draw[line,->] (z3.north west) -- (f3.south);
\draw[line,->] (z3.north east) -- (b3.south);

\draw[line,->] (f1.north) -- (c1.south west);
\draw[line,->] (b1.north) -- (c1.south east);
\draw[line,->] (c1.north) -- (x1.south);

\draw[line,->] (f2.north) -- (c2.south west);
\draw[line,->] (b2.north) -- (c2.south east);
\draw[line,->] (c2.north) -- (x2.south);

\draw[line,->] (f3.north) -- (c3.south west);
\draw[line,->] (b3.north) -- (c3.south east);
\draw[line,->] (c3.north) -- (x3.south);

\node[lab] at ({0.8+\s},4.1) {predicted stage decisions};
\node[lab] at ({0.8+\s},-3.0) {history-dependent stage information};

\node[align=center, lab] at (-2.2,1.6) {Backward\\layer};
\node[align=center, lab] at (-2.2,0) {Forward\\layer};

\end{tikzpicture}%
}
\caption{Bidirectional LSTM architecture for sequential decision learning in multi-stage optimization. At stage \(t\), the input \(z_t=(\xi_{[t]},\psi_t,\hat{x}_{t-1})\) combines the revealed information history, stage-dependent optimization data, and the previous decision prediction.}
\label{fig:bilstm_msp_decisions}
\end{figure}

As illustrated in Figure~\ref{fig:bilstm_msp_decisions}, the bidirectional LSTM architecture maps the time-indexed stage information to a structured sequence of stage decisions rather than isolated pointwise predictions \citep{yilmaz2023lstm}. At stage \(t\), the input vector \(z_t=(\xi_{[t]},\psi_t,\hat{x}_{t-1})\) summarizes the revealed history \(\xi_{[t]}\), stage-dependent optimization data \(\psi_t\), and previous decision prediction \(\hat{x}_{t-1}\). Recurrent forward and backward layers process the sequence, their hidden representations are concatenated, and the combined representation is mapped to the stage-wise prediction \(\hat{x}_t(\xi_{[t]})\). Because the architecture is bidirectional, it should be interpreted as an offline solution-learning model for fully specified planning instances; an online implementation would restrict the learned decision rule to information available at the decision epoch.

From an OR/MS perspective, the model is not merely fitting correlations in sequential data. Rather, it learns recurring patterns induced by the underlying optimization structure, including temporal coupling across stages, dependencies among decisions, and regularities shaped by inventory balance equations, resource limits, and the objective function itself \citep{bengio2021machine, yilmaz2024expandable}. Although feasibility is not enforced explicitly within the network, these structural relationships can be reflected implicitly in the learned representation. The output is therefore best viewed as a \emph{sequential decision trajectory}, in which each stage decision is informed by both evolving input data and its role in the broader optimization sequence.

Evaluation of learning-based decision frameworks requires a broader set of metrics than classical optimization alone, including objective quality, feasibility, infeasibility rate, computational gain, early incumbent quality, and robustness in downstream decision making \citep{yilmaz2023lstm,yilmaz2024expandable}. These properties make LSTM-based learning--optimization frameworks attractive for multi-stage applications such as production planning, inventory control, energy scheduling, and epidemic response. In such settings, the learned decision sequence can serve as either a fast approximate decision rule or as guidance for downstream optimization through fixing, partial fixing, or warm-starting. Related work on capacitated lot sizing suggests that the same learn-to-decision logic can also be pursued with temporal CNNs and transformer-based architectures for dynamic mixed-integer optimization \citep{choi2024tcnn,buyuktahtakin2024transformers}.

\subsection{Expandable Learning--Optimization Architectures}

A central limitation of many learning-based optimization methods is that they are tied to the instance sizes seen during training. For OR/MS applications, practical value requires more than high predictive accuracy on a fixed benchmark; it requires learned models that can transfer across longer planning horizons and larger problem dimensions without rebuilding the full pipeline. PredOpt \citep{yilmaz2024expandable} addresses this challenge through an expandable learning--optimization architecture that combines a sequence-to-sequence predictor with time-wise and item-wise expansion mechanisms.

PredOpt connects the neural architectures discussed in Section~\ref{sec:architectures} with the decision frameworks in Section~\ref{sec:preliminaries}. Built on a sequence-to-sequence encoder--decoder learning architecture with attention \citep{sutskever2014sequence,bahdanau2015neural,luong2015effective}, it learns reusable stage-wise decision patterns while the optimization model preserves feasibility and solution quality, capturing temporal dependencies across stages while also learning recurring structural patterns induced by stage data, coupled decisions, constraints, and the objective. As shown in Figure~\ref{fig:predopt_compact_detailed}, PredOpt presents a modular predict--screen--optimize pipeline, providing a concrete instantiation of the learning-to-optimize paradigm summarized in Figure~\ref{fig:learning_optimization_paradigms}(c). The predictor identifies promising decision patterns, infeasible assignments are screened out, and the resulting partial solution is passed to the optimization solver through fixing or warm-starting.

PredOpt is best viewed as a hybrid decision pipeline rather than a stand-alone predictor. The learned model provides fast guidance, while feasibility screening of the predicted solution and solver-based refinement preserve the integrity of both the optimization model and sequential decision. This complementarity is critical in sequential decision problems, where scalability and implementability must coexist.

\begin{figure}[t]
\centering
\resizebox{\textwidth}{!}{%
\begin{tikzpicture}[
    >=Latex,
    font=\small,
    node distance=6mm,
    every node/.style={align=center},
    block/.style={
        draw,
        rounded corners=3pt,
        line width=0.9pt,
        minimum height=1.25cm,
        text width=3.0cm,
        inner sep=4pt,
        fill=gray!8
    },
    finalblock/.style={
        draw,
        rounded corners=3pt,
        line width=0.9pt,
        minimum height=1.25cm,
        text width=2.8cm,
        inner sep=4pt,
        fill=gray!15
    },
    arrow/.style={->, line width=0.95pt}
]

\node[block] (input) {{\itshape Test instance}\\longer horizon / larger dimension};

\node[block, right=8mm of input] (predict)
{{\itshape Expandable encoder--decoder}\\time-wise and item-wise prediction};

\node[block, right=8mm of predict] (screen)
{{\itshape Feasibility screening}\\remove infeasible assignments};

\node[block, right=8mm of screen] (solver)
{{\itshape Optimization solver}\\variable fixing or warm-start};

\node[finalblock, right=8mm of solver] (solution)
{{\itshape Refined solution}\\feasible, high-quality output};

\draw[arrow] (input) -- (predict);
\draw[arrow] (predict) -- (screen);
\draw[arrow] (screen) -- (solver);
\draw[arrow] (solver) -- (solution);

\end{tikzpicture}%
}
\caption{Compact view of the expandable PredOpt architecture: prediction, feasibility screening, and solver-based refinement.}
\label{fig:predopt_compact_detailed}
\end{figure}

\subsubsection{Generalization Across Horizons and Incremental Expansion}

The first source of expandability in PredOpt is generalization across planning horizons. Because the encoder--decoder with attention is not tied to a fixed output length, a model trained on shorter horizons can be applied to longer ones while preserving temporal dependence across stages \citep{yilmaz2024expandable, bahdanau2015neural,luong2015effective}. This is particularly relevant in OR/MS settings such as rolling-horizon planning, seasonal operations, and adaptive scheduling, where the horizon length often changes over time. Importantly, this expansion is incremental rather than all-at-once: the trained model transfers across horizon segments without retraining, reusing the learned temporal structure on longer sequences. This makes the framework broadly applicable to families of related sequential problems, not just a single fixed-scale benchmark.

\subsubsection{Item-Wise Expansion and Higher-Dimensional Generalization}

The second source of expandability is generalization across problem dimension. In many OR/MS problems, the mathematical structure remains the same, while the number of items, products, or decision components changes. PredOpt addresses this challenge through an item-wise expansion mechanism that allows a model trained on smaller instances to be reused on larger ones without retraining \citep{yilmaz2024expandable}. As in Figure~\ref{fig:predopt_compact_detailed}, the predictor remains embedded in the same predict--screen--optimize pipeline, but its scope is extended to a larger set of elements through repeated subset-based prediction. Algorithm~\ref{alg:dimension_wise_generalization} summarizes the mechanism. The trained predictor is treated as a reusable local decision module: it is applied to overlapping subsets whose size matches the training dimension, the repeated item-level predictions are aggregated, and the resulting larger-dimensional partial solution is passed, when needed, to a downstream optimizer for feasibility restoration and refinement.

\begin{algorithm}[h!]
\caption{Subset-Based Item-Wise Expansion for Larger-Dimensional Instances \cite{yilmaz2024expandable}}
\label{alg:dimension_wise_generalization}
\small
\textbf{Input:} trained predictor \(M\) built on item dimension \(d^M\); test-instance data \(\alpha\); full item set \(\mathcal{D}\); minimum coverage threshold \(\delta\); aggregation rule \(\mathcal{A}(\cdot)\).

\textbf{Output:} aggregated prediction \(\hat{x}_{j,1:T}\) for each item \(j \in \mathcal{D}\).

\begin{enumerate}
    \item Initialize cumulative predictions \(\bar{x}_{j,1:T} \leftarrow 0\) and prediction counts \(\gamma_j \leftarrow 0\) for all \(j \in \mathcal{D}\).
    \item \textbf{while} there exists an item \(j \in \mathcal{D}\) with \(\gamma_j < \delta\) \textbf{do}
    \begin{enumerate}
        \item Sample a subset \(S \subseteq \mathcal{D}\) such that \(|S| = d^M\).
        \item Construct the restricted input \(\alpha_S\) for subset \(S\).
        \item Make a forward pass with the trained predictor:
        \[
        \hat{x}_{S,1:T} = M(\alpha_S).
        \]
        \item For each \(j \in S\), update
        \[
        \bar{x}_{j,1:T} \leftarrow \bar{x}_{j,1:T} + \hat{x}_{j,1:T},
        \qquad
        \gamma_j \leftarrow \gamma_j + 1.
        \]
    \end{enumerate}
    \item For each \(j \in \mathcal{D}\), compute the final aggregated prediction
    \[
    \hat{x}_{j,1:T} \leftarrow \mathcal{A}\!\left(\bar{x}_{j,1:T},\gamma_j\right).
    \]
    \item Optionally pass the aggregated prediction \(\hat{x}\) to a downstream optimization model for feasibility restoration and solution refinement.
\end{enumerate}
\end{algorithm}

Returning to the lot-sizing example, suppose that the trained predictor \(M\) was built on \(d^M=3\) products, while the test instance has five products, \(\mathcal{D}=\{1,2,3,4,5\}\). PredOpt samples overlapping subsets of size three, such as \(S_1=\{1,2,3\}\), \(S_2=\{1,3,4\}\), and \(S_3=\{2,4,5\}\), along with additional subsets as needed until the desired coverage threshold is reached. For each subset \(S\), the restricted input \(\alpha_S\) contains the costs, demands, capacities, and initial inventories for the products in \(S\), and the forward pass \(M(\alpha_S)\) predicts local setup trajectories, for example \(\hat{u}_{j,1:T}\), for the products in that subset.

Shared resources are scaled so that each restricted instance reflects the resource tightness
of the full problem. In multi-item capacitated lot sizing, where products share period capacity \(\sum_{j\in\mathcal{D}}q_{jt}\le C_t\), one may set \(C_{t,S}=C_t\sum_{j\in S}d_{jt}/\sum_{j\in\mathcal{D}}d_{jt}\), when the denominator is positive. For instance, if \(C_t=100\), the total demand at period \(t\) is \(200\), and subset \(S\) accounts for demand \(50\), then \(C_{t,S}=25\). Because each product appears in several subsets, its repeated predictions are aggregated. For example, after several subset evaluations, if product \(1\) receives setup predictions \((1,1,0,1)\) for period \(t\), averaging gives \(0.75\), which can be used as a threshold to obtain \(\hat{u}_{1t}=1\). The resulting larger-dimensional setup prediction is then passed to the downstream optimizer, which restores global feasibility across all products, capacities, inventory balances, and integrality constraints.

\subsubsection{Computational Evidence and Practical Scalability}

The computational results in \citep{yilmaz2024expandable} show that expandability leads to substantial practical gains. For the Multi-item Capacitated Lot Size Problem (MCLSP) with 12 items, PredOpt achieves about \(99\%\) prediction accuracy and only \(0.04\%\)--\(0.11\%\) optimality gap across horizons from \(T=30\) to \(T=150\), while requiring only \(4.6\) seconds at \(T=150\) compared to \(25{,}085\) seconds for CPLEX. In Multi-dimensional Knapsack (MSMK), PredOpt remains competitive as dimension grows, with roughly \(91.6\%\)--\(93.8\%\) prediction accuracy, small optimality gaps, and solution times of only seconds on instances for which direct CPLEX runtimes become orders of magnitude larger. These results suggest that expandability is not merely a modeling convenience; it provides significant computational leverage on larger and longer sequential decision problems.

In sum, expandable learning--optimization architectures address a central OR/MS challenge: how to transfer learned decision structure beyond the training scale while preserving feasibility and solution quality. PredOpt illustrates how deep learning captures recurring temporal and structural patterns across instances, while optimization enforces global consistency, delivering computational leverage at scale without sacrificing solution quality. The result is not just a predictor, but a scalable decision pipeline anchored in the logic of mathematical optimization \citep{yilmaz2024expandable}.

\subsection{Non-Anticipative Learning--Optimization Frameworks}

Sequential decision making under uncertainty requires decisions to remain implementable as information is revealed over time. In multi-stage stochastic programming, this is captured by \emph{non-anticipativity}: if two scenarios share the same realization history up to stage \(t\), then the stage-\(t\) decision must also be the same \citep{shapiro2009lectures}. Standard sequence models do not enforce this requirement and may therefore exploit distinctions across scenarios that are not yet observable, producing predictions that are statistically plausible but operationally invalid. To address this challenge, \citet{yilmaz2025nonanticipative} propose the \emph{Non-anticipative Encoder--Decoder with Attention} (NEDA), which extends the attention-based encoder--decoder architecture of \citet{luong2015effective} to the scenario-tree structure of multi-stage stochastic programs.

The NEDA framework is intuitive: at stage \(t\), scenarios that are still indistinguishable to the decision maker should also be indistinguishable to the neural network. Let \(h^e_{t,s}=[\overrightarrow{h}^e_{t,s},\overleftarrow{h}^e_{t,s}]\) denote the hidden state of the encoder for scenario \(s\) at stage \(t\), and let \(\Omega_{t,s}\) be the set of scenarios that share the same information history up to stage \(t\) as in scenario \(s\). NEDA enforces non-anticipativity by replacing the scenario-specific encoder state with a scenario-group average:
\begin{equation}
\bar{h}^e_{t,s}
=
\frac{\sum_{s' \in \Omega_{t,s}} h^e_{t,s'}}{|\Omega_{t,s}|}.
\label{eq:neda_avg_hidden}
\end{equation}

Continuing the lot-sizing example, suppose a three-stage scenario tree has four terminal demand scenarios. At stage \(t=2\), the decision maker has not yet observed which terminal scenario will occur. Instead, scenarios \(s_1\) and \(s_2\) pass through the same stage-2 information node, while scenarios \(s_3\) and \(s_4\) pass through a different stage-2 information node. Non-anticipativity therefore requires \(x_2^{s_1}=x_2^{s_2}\) and \(x_2^{s_3}=x_2^{s_4}\), even though the scenarios may diverge later. NEDA enforces this logic at the representation level. For example, if \(h^e_{2,s_1}=(0.2,0.8)\) and \(h^e_{2,s_2}=(0.4,0.6)\), both scenarios use the averaged representation \(\bar{h}^e_{2,s_1}=\bar{h}^e_{2,s_2}=(0.3,0.7)\) before decoding the stage-2 decision. Thus, scenarios that are indistinguishable at stage 2 are passed to the decoder through a common information-restricted representation. Once the scenarios separate at a later stage, their representations and decisions may differ.

Because all scenarios in the same information set use the same averaged hidden state, the decoder receives a common representation for indistinguishable scenarios and can produce the same stage-\(t\) prediction by construction. In effect, the model is prevented from encoding differences that the decision maker has not yet observed.

This is the key innovation over deterministic frameworks such as PredOpt. NEDA does not merely learn the temporal structure between stages; it also aligns hidden-state representations throughout the scenario tree so that the learned decision map respects the information structure of the stochastic program. When the averaged hidden states are highly consistent, the model can make confident common predictions; when they are not, the downstream ScenPredOpt framework relies more heavily on feasibility repair, LP-based heuristics, and solver refinement to recover an implementable high-quality solution \citep{yilmaz2025nonanticipative}.

\begin{figure}[h!]
\centering
\includegraphics[width=0.92\linewidth]{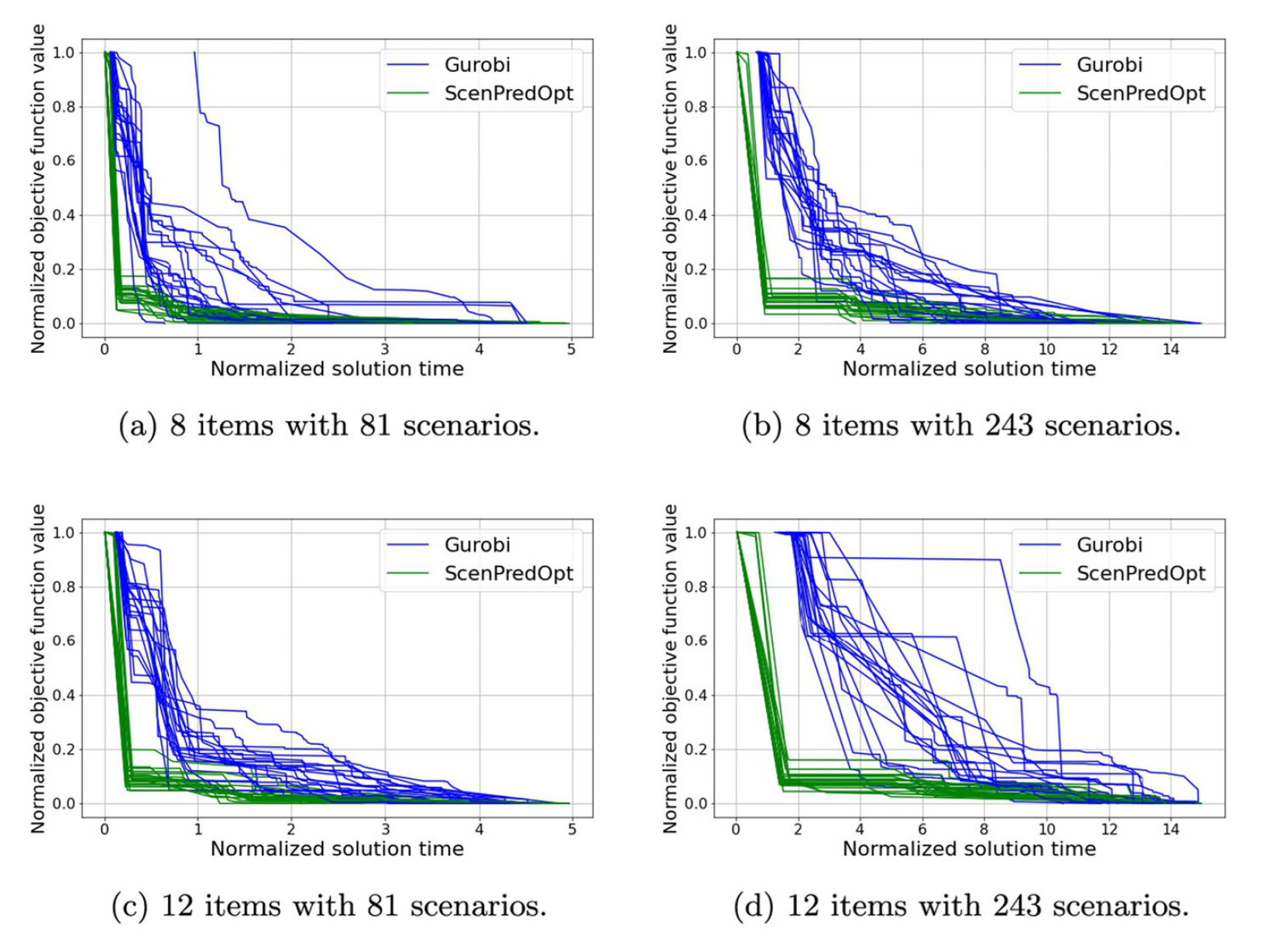}
\caption{Normalized objective-value trajectories for Gurobi and ScenPredOpt across representative instances. ScenPredOpt attains stronger early incumbents and faster practical improvement during the initial stages of the solution process. Adapted from \citet{yilmaz2025nonanticipative} under the Creative Commons Attribution 4.0 International License (CC BY 4.0).}
\label{fig:scenpredopt_progress}
\end{figure}

\paragraph{Computational Performance.}
Figure~\ref{fig:scenpredopt_progress} shows that ScenPredOpt reduces the normalized objective value faster and more consistently than Gurobi across representative instance classes in a cost minimization problem. The gap widens as uncertainty and problem size increase, with ScenPredOpt reaching near-optimal solutions in a small fraction of the total runtime while Gurobi exhibits slower and more variable convergence. This behavior highlights the practical value of integrating non-anticipativity and feasibility within the learning--optimization pipeline.

Across the reported stochastic test sets, ScenPredOpt delivers a strong speed--quality trade-off relative to direct Gurobi solves and benchmark methods, including Progressive Hedging (PH) \citep{watson2011progressive}, SDDiP \citep{ding2019sddip,zou2019sddip} when applicable, and tailored heuristics \citep{absi2019,bertsimas2002}. In Stochastic Multi-dimensional Knapsack (SMSMK) experiments, it reduces Gurobi runtimes from \(1155\)--\(7200\) seconds to \(1\)--\(69\) seconds while maintaining small optimality gaps of \(0.77\%\)--\(1.31\%\) and prediction accuracies above \(91\%\) \citep{yilmaz2025nonanticipative}. PH can sometimes achieve smaller gaps, but it remains substantially slower; tailored heuristics are faster than PH but generally less accurate. Thus, the main computational advantage of ScenPredOpt is its ability to identify strong incumbents quickly while preserving near-optimal solution quality, especially in larger, scenario-rich settings where exact- and decomposition-based methods become increasingly expensive.

A related large-scale direction embeds transformer-based learning within a Benders-type decomposition framework \citep{benders1962partitioning} for stochastic combinatorial optimization, illustrating how learned models can accelerate structured exact methods rather than replace them \citep{choi2025benderstransformers}.

\section{Deep Reinforcement Learning for Sequential Decision Making}
\label{sec:drl}

This section examines deep reinforcement learning (DRL) as a policy-approximation framework for sequential OR/MS problems in which uncertainty unfolds over time, and explicit stochastic optimization models may be difficult to specify or solve at scale. We focus on three OR/MS-oriented uses: simulation-integrated policy learning, approximation of stochastic-programming recourse and stage-wise decisions, and optimization-informed policy learning for combinatorial optimization. The goal is not to survey the broader DRL literature, including model-based RL and world models \citep{moerland2023model}, offline RL \citep{levine2020offline}, safe and constrained RL \citep{garcia2015comprehensive,achiam2017constrained}, or multi-agent RL \citep{zhang2021multi}. Instead, we emphasize settings where feasibility, integrality, temporal coupling, recourse, and domain-specific state/action design determine whether DRL can serve as a useful component of OR/MS decision systems.

\subsection{DRL as Data-Driven Policy Approximation for Stochastic Decision Problems}

Building on the MDP, DP, and RL foundations introduced in Section~\ref{sec:preliminaries}, we now focus on DRL as a scalable policy-learning framework for sequential decision problems under uncertainty. Whereas classical DP characterizes optimality through Bellman recursion \citep{bellman1957dynamic}, DRL replaces exact value-function computation and full state enumeration with sample-based learning and deep function approximation. This makes it especially relevant for OR/MS settings with large state spaces, sequential uncertainty, and complex dynamics, where exact dynamic programming or stochastic optimization becomes computationally impractical.

In value-based methods under a cost-minimization convention, the action-value function \(Q^\pi(s,a)\) represents an expected cost-to-go and is updated using temporal-difference learning, as in Q-learning \citep{watkins1989learning}. With one-step target
\begin{equation}
y_t = c_t + \gamma \min_{a'} Q_t(s_{t+1},a'),
\end{equation}
where \(y_t\) denotes the temporal-difference target and \(c_t\) denotes the one-step cost, the update is
\begin{equation}
Q_{t+1}(s_t,a_t)
\leftarrow
Q_t(s_t,a_t)
+ \alpha \bigl[y_t - Q_t(s_t,a_t)\bigr].
\end{equation}
Deep Q-Networks (DQN) \citep{mnih2015human} replace tabular value functions with neural approximators \(Q(s,a;\theta)\), enabling learning in large state spaces. Policy-based methods instead parameterize the decision rule directly through \(\pi_\theta(a \mid s)\) and optimize expected return or cost with gradient-based methods \citep{williams1992simple}, while actor--critic approaches combine value approximation with policy optimization to improve learning stability \citep{konda2000actor}.

The transition from RL to DRL occurs when deep architectures are used to represent value functions or policies. When the full state of the system is not directly observable, the natural extension is a partially observable Markov decision process (POMDP), in which decisions are based on observations or belief states rather than fully observed states \citep{kaelbling1998pomdp}. In this setting, Deep Recurrent RL methods such as Deep Recurrent Q-Networks (DRQN) provide a natural bridge by using memory to summarize observation histories and approximate latent state information relevant to decision making \citep{hausknecht2015drqn}. For OR/MS applications, this matters because many sequential decision problems involve long-range temporal dependence, partial observability, and high-dimensional state descriptions. Recurrent neural networks, LSTMs \citep{hochreiter1997long}, and attention-based transformers \citep{vaswani2017attention} provide the ability to capture evolving state information, intertemporal dependencies, and longer-horizon effects that are often central to planning, control, and resource allocation problems.

From an OR/MS standpoint, DRL may be interpreted as a simulation-based approximation to stochastic dynamic programming and, more broadly, to multi-stage stochastic programming (MSP). In MSP, decisions are optimized over scenario trees subject to non-anticipativity and recourse structure. DRL replaces this explicit scenario-tree representation with a parameterized policy
\begin{equation}
x_t = \pi_\theta(s_t),
\label{eq:policy_state_mapping}
\end{equation}
where the action \(x_t\) is the stage-\(t\) decision and the state \(s_t\) summarizes the information available at stage \(t\). In this sense, non-anticipativity holds automatically: because \(x_t = \pi_\theta(s_t)\) is a fixed function evaluated at the realized state, the decision depends only on information available by stage \(t\) and never on future realizations, provided \(s_t\) is measurable with respect to that information, that is, a function of history \(\xi_{[t]}\) rather than of future outcomes. Non-anticipativity is thus inherited from the policy form itself, and the challenge shifts from enforcing it to ensuring that the state, action space, and learning environment together produce decisions that remain feasible and implementable.

This viewpoint clarifies both the promise and limitations of DRL for OR/MS. Its strength lies in scalability: it avoids the exponential growth of scenario trees and can be learned in rich simulation environments where exact modeling is impractical. Its limitation is that feasibility, integrality, and structural consistency are not guaranteed by default. As a result, effective uses of DRL often combine policy learning with optimization structure, simulation models, or feasibility-enforcing mechanisms, so that the flexibility of learning is complemented by the rigor of mathematical optimization.

\subsection{DRL Integration with Simulation}

A particularly important role for deep reinforcement learning in operations research arises when the decision environment can be simulated with high fidelity but cannot be represented through a tractable analytical transition model. In such settings, the main difficulty is not the absence of system knowledge, but the inability to express the system dynamics in a form amenable to classical dynamic programming, stochastic programming, or other policy-based optimization methods. Simulation-integrated DRL addresses this challenge by treating a simulation model as the learning environment and training a policy through repeated interaction with that environment \citep{BushajEtAl2023}. This creates a practical bridge between mechanistic system modeling and sequential decision optimization, especially in settings that are nonlinear, path-dependent, partially observed, or shaped by heterogeneous interacting agents \citep{sutton2018reinforcement,mnih2015human,Powell2019,Powell2022}.

From a decision making perspective, this paradigm expands the scope of sequential optimization beyond problems that can be summarized by compact state transitions or scenario trees. Rather than requiring a closed-form stochastic model, one can rely on a simulation engine that already captures operational rules, feedback effects, and domain-specific complexity. The DRL agent then learns a policy directly from the simulated evolution of the system, with reward or cost signals encoding the underlying objectives and trade-offs. In this sense, simulation-integrated DRL can be viewed as a policy-learning layer built on top of a computational model of the system \citep{Powell2019,Powell2022,sutton2018reinforcement, BushajEtAl2023}.

\begin{figure}[h!]
\centering
\begin{subfigure}[t]{0.82\textwidth}
\centering
\begin{tikzpicture}[
    >=Latex,
    font=\small,
    box/.style={rectangle, draw, rounded corners=2pt, minimum width=3.2cm, minimum height=1.0cm, align=center},
    line/.style={->, thick}
]
\node[box] (state) {Current state \(s_t\)};
\node[box, below=0.9cm of state] (agent) {DRL agent};
\node[box, below=0.9cm of agent] (sim) {Simulation environment};
\node[box, below=0.9cm of sim] (next) {Reward \(r_t\), next state \(s_{t+1}\)};
\node[box, right=1.8cm of next] (update) {Policy update};
\draw[line] (state) -- (agent);
\draw[line] (agent) -- node[right] {\(a_t\)} (sim);
\draw[line] (sim) -- (next);
\draw[line] (next.east) -- ++(0.9,0) |- (update.west);
\draw[line] (update.north) |- (agent.east);
\end{tikzpicture}
\caption{High-level workflow}
\label{fig:simdrl_panel_a}
\end{subfigure}

\vspace{1.5em}

\begin{subfigure}[t]{0.92\textwidth}
\centering
\begin{tikzpicture}[
    >=Latex,
    font=\small,
    box/.style={rectangle, draw, rounded corners=2pt, minimum width=2.4cm, minimum height=0.95cm, align=center},
    tbox/.style={rectangle, draw, rounded corners=2pt, minimum width=1.1cm, minimum height=0.55cm, align=center},
    line/.style={->, thick}
]
\node[tbox] (t1) at (0,2.7) {\(t=1\)};
\node[tbox] (t2) at (4.3,2.7) {\(t=2\)};
\node[tbox] (tT) at (9.6,2.7) {\(t=T-1\)};
\node[box] (sim1) at (0,0) {Simulation};
\node[box] (agent2) at (4.3,0) {DRL agent};
\node[box] (sim2) at (4.3,1.5) {Simulation};
\node[box] (agentT) at (9.6,1.5) {DRL agent};
\node[box] (simT) at (9.6,0) {Simulation};
\draw[line] (sim1.east) -- node[above]{\(s_1\)} (agent2.west);
\draw[line] (agent2.north) -- node[right]{\(a_1\)} (sim2.south);
\draw[line] (sim2.east) -- node[above]{\(r_1,\ s_2\)} ++(0.9,0);
\draw[line] ([xshift=-0.8cm]agentT.west) -- (agentT.west);
\draw[line] (agentT.south) -- node[right]{\(a_{T-1}\)} (simT.north);
\node at (6.95,1.72) {\Large \(\cdots\)};
\node at (6.95,0.22) {\Large \(\cdots\)};
\node[right=0.1cm of simT] {\(r_{T-1},\ s_T\)};
\end{tikzpicture}
\caption{Sequential time-indexed view}
\label{fig:simdrl_panel_b}
\end{subfigure}
\caption{Two complementary views of simulation-integrated deep reinforcement learning. Panel~(a) presents the high-level closed-loop workflow linking the DRL agent, the simulation environment, and policy updates. Panel~(b) presents the same framework in time-indexed form, emphasizing the sequential evolution of states, actions, and rewards across decision epochs. Adapted and generalized from \citet{BushajEtAl2023}.}
\label{fig:sim_drl_combined_vertical}
\end{figure}
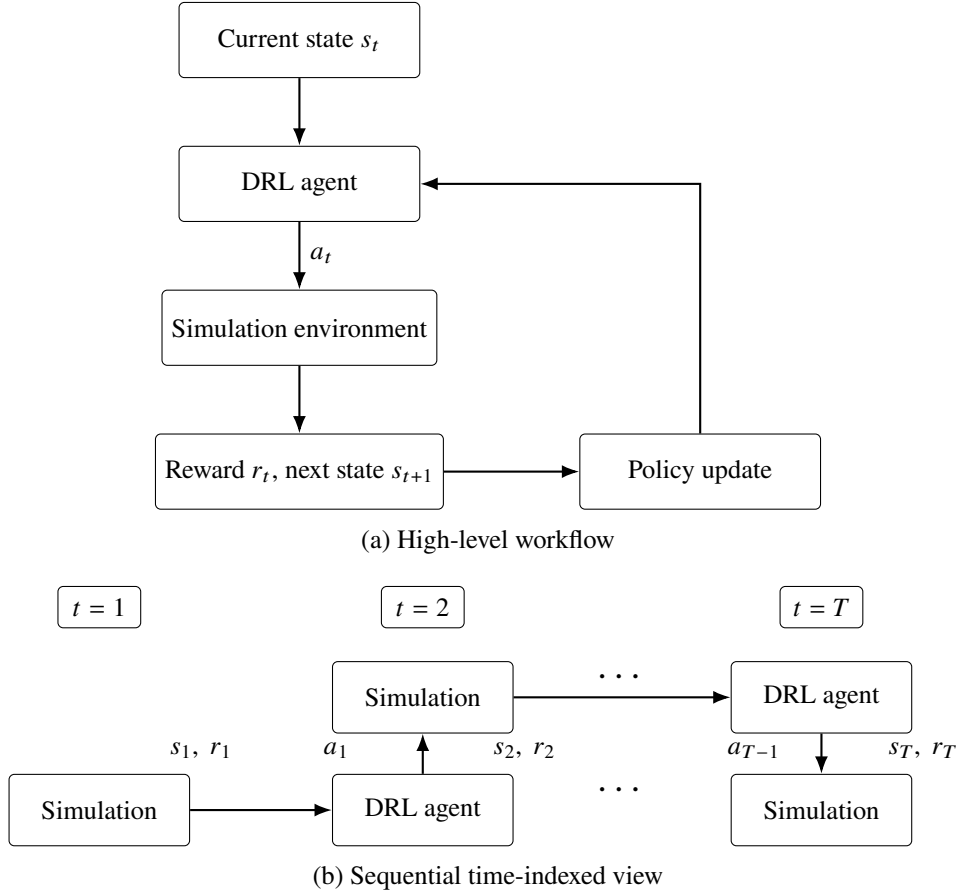

Figure~\ref{fig:sim_drl_combined_vertical} provides two complementary views of this integration and makes clear why simulation becomes part of the policy-learning loop \citep{BushajEtAl2023}. Panel~(a) presents the closed-loop high-level workflow that connects the current state, the DRL agent, the simulation environment, the resulting reward and the next state, and the subsequent policy update. Panel~(b) makes the sequential structure explicit over time. At each decision epoch \(t\), the simulator provides a state description \(s_t\), the DRL agent selects an action \(a_t\), and the simulator propagates the consequences of that action to generate a reward \(r_t\) together with the information needed to construct the next state \(s_{t+1}\). Thus, learning proceeds by repeatedly applying the current policy within the simulator:
\begin{equation}
\begin{aligned}
a_t &\sim \pi_\theta(\cdot \mid s_t),\\
(r_t,s_{t+1}) &= \mathcal{S}(s_t,a_t),
\qquad t=1,\ldots,T-1.
\end{aligned}
\label{eq:simulator_interaction}
\end{equation}
Here, \(\mathcal{S}\) denotes the simulation environment: instead of specifying an explicit analytical transition kernel, the simulator returns the reward or cost feedback \(r_t\) and the next state \(s_{t+1}\) induced by action \(a_t\) \citep{BushajEtAl2023,sutton2018reinforcement,Powell2022}.

Epidemic control provides a compelling example. In such problems, the spread of the disease depends on rich interactions between contact behavior, intervention timing, vaccination, compliance, and population heterogeneity. These effects are difficult to encode through explicit transition equations, yet they can often be represented naturally through agent-based simulation. Bushaj et al.~\citep{BushajEtAl2023} operationalize this idea through a simulation-deep reinforcement learning framework in which a DRL agent interacts with the Covasim agent-based epidemic simulator \citep{KerrEtAl2021} to learn adaptive intervention policies under multi-objective reward trade-offs. This modeling logic is also consistent with related DRL-based epidemic-control studies that learn mitigation policies in large-scale simulated environments \citep{LibinEtAl2021}.

More broadly, simulation-integrated DRL uses complex simulators as environments for learning, testing, and refining sequential decision policies. Rather than using simulation only for post hoc evaluation, it uses simulation as part of the optimization loop to learn responsive, adaptive, and operationally meaningful policies in environments where explicit stochastic formulations are unavailable or too restrictive for practical solution \citep{Powell2019,Powell2022,BushajEtAl2023}.

\subsection{DRL for Stochastic Programs}

Another important direction is the use of DRL to approximate solutions to stochastic programs. A leading example is \citet{YilmazBuyuktahtakin2024}, who recast a two-stage stochastic program as a multi-agent DRL framework. The central idea is to replace the repeated solution of a large scenario-based optimization model with learned stage-wise decision policies. This is especially appealing in two-stage stochastic programming, where the problem naturally separates into here-and-now decisions and recourse decisions taken after uncertainty is revealed.

The key OR/MS insight is that the second-stage recourse decision can be represented not only as a nested optimization problem but also as a learnable decision process. In the framework of \citet{YilmazBuyuktahtakin2024}, Agent~1 learns the first-stage decision \(x\), while Agent~2 learns the second-stage recourse decision \(y\) conditional on the first-stage decision and the realized uncertainty. In this way, the logic of stochastic programming is preserved, but repeated optimization is replaced by learned decision rules.

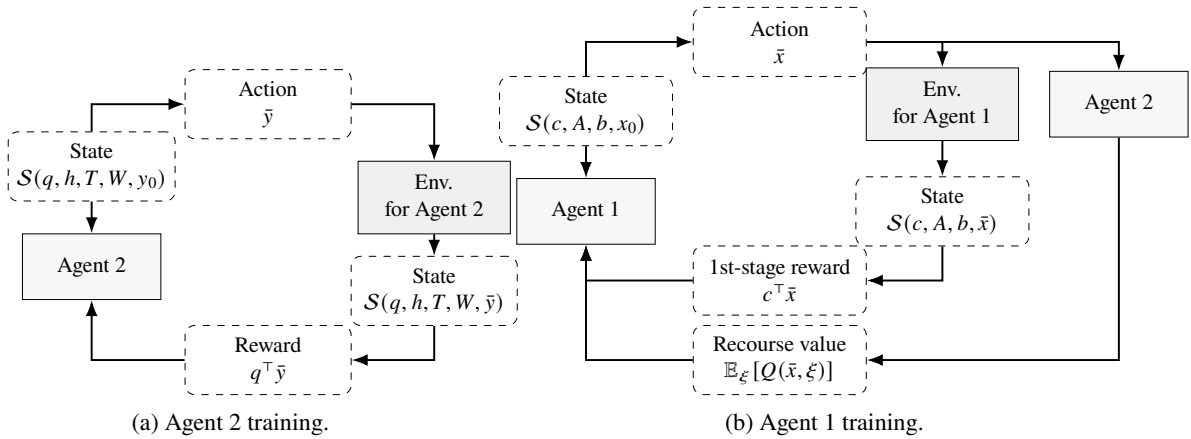
\begin{figure}[h!]
\centering
\resizebox{0.92\textwidth}{!}{%
\begin{minipage}{\textwidth}
\centering

\begin{subfigure}[t]{0.39\textwidth}
\centering
\begin{tikzpicture}[
    >=Latex,
    font=\footnotesize,
    agent/.style={
        rectangle, draw,
        minimum width=2.0cm, minimum height=0.95cm,
        align=center, fill=gray!6
    },
    env/.style={
        rectangle, draw,
        minimum width=2.2cm, minimum height=1.05cm,
        align=center, fill=gray!12
    },
    info/.style={
        rectangle, draw, dashed, rounded corners=4pt,
        minimum width=2.4cm, minimum height=1.0cm,
        align=center
    },
    line/.style={->, thick}
]

\node[agent] (agent2) at (0,0) {Agent 2};

\node[info] (state0) at (0,1.45)
{State\\\(\mathbf{s}(q,h,T,W,y_0)\)};

\node[info] (action) at (2.55,2.35)
{Action\\\(\bar y\)};

\node[env] (env2) at (4.95,1.0)
{Env.\\for Agent 2};

\node[info] (state1) at (4.95,-0.35)
{State\\\(\mathbf{s}(q,h,T,W,\bar y)\)};

\node[info] (signal) at (2.55,-1.35)
{Learning signal\\\(q^\top \bar y\)};

\draw[line] (state0.south) -- (agent2.north);
\draw[line] (state0.north) |- (action.west);
\draw[line] (action.east) -| (env2.north);
\draw[line] (env2.south) -- (state1.north);
\draw[line] (state1.south) |- (signal.east);
\draw[line] (signal.west) -| (agent2.south);

\end{tikzpicture}
\caption{Agent 2 training.}
\label{fig:2srl_agent2}
\end{subfigure}
\hfill
\begin{subfigure}[t]{0.57\textwidth}
\centering
\begin{tikzpicture}[
    >=Latex,
    font=\footnotesize,
    agent/.style={
        rectangle, draw,
        minimum width=2.0cm, minimum height=0.95cm,
        align=center, fill=gray!6
    },
    env/.style={
        rectangle, draw,
        minimum width=2.2cm, minimum height=1.05cm,
        align=center, fill=gray!12
    },
    info/.style={
        rectangle, draw, dashed, rounded corners=4pt,
        minimum width=2.5cm, minimum height=1.0cm,
        align=center
    },
    line/.style={->, thick}
]

\node[agent] (agent1) at (0,0) {Agent 1};

\node[info] (stateinit) at (0,1.45)
{State\\\(\mathbf{s}(c,A,b,x_0)\)};

\node[info] (actionx) at (2.8,2.45)
{Action\\\(\bar x\)};

\node[env] (env1) at (5.15,1.55)
{Env.\\for Agent 1};

\node[agent] (agent2b) at (7.7,1.55) {Agent 2};

\node[info] (statex) at (5.15,0.0)
{State\\\(\mathbf{s}(c,A,b,\bar x)\)};

\node[info] (signal1) at (2.8,-1.0)
{1st-stage signal\\\(c^\top \bar x\)};

\node[info] (signal2) at (2.8,-2.15)
{Recourse value\\\(\mathbb{E}_{\xi}[Q(\bar x,\xi)]\)};

\draw[line] (stateinit.south) -- (agent1.north);
\draw[line] (stateinit.north) |- (actionx.west);

\draw[line] (actionx.east) -| (env1.north);
\draw[line] (actionx.east) -| (agent2b.north);

\draw[line] (env1.south) -- (statex.north);
\draw[line] (statex.south) |- (signal1.east);

\draw[line] (agent2b.south) |- (signal2.east);

\draw[line] (signal1.west) -| (agent1.south);
\draw[line] (signal2.west) -| (agent1.south);

\end{tikzpicture}
\caption{Agent 1 training.}
\label{fig:2srl_agent1}
\end{subfigure}

\caption{Two-stage DRL training overview. Agent~2 learns recourse actions after scenario realization, while Agent~1 learns first-stage decisions using both immediate objective contribution and downstream recourse feedback. Adapted from \citet{YilmazBuyuktahtakin2024}.}
\label{fig:two_stage_drl}
\end{minipage}%
}
\end{figure}

Figure~\ref{fig:two_stage_drl} illustrates the training structure. Panel~(a) shows the training loop for Agent~2, which learns a recourse policy after a scenario is realized. Starting from a second-stage state, schematically denoted by \(\mathbf{s}(q,h,T,W,y_0)\), Agent~2 chooses a recourse action \(\bar y\), transitions to the updated state \(\mathbf{s}(q,h,T,W,\bar y)\), and receives a learning signal based on the second-stage value, represented here by \(q^\top \bar y\). Thus, Agent~2 learns how to map the realized second-stage uncertainty into corrective recourse actions. Panel~(b) highlights the more distinctive feature of the framework: Agent~1 is trained not only on the immediate first-stage contribution \(c^\top \bar x\), but also on the downstream value induced by that decision, represented by \(\mathbb{E}_{\xi}[Q(\bar x,\xi)]\). From the first-stage state \(\mathbf{s}(c,A,b,x_0)\), Agent~1 selects a first-stage action \(\bar x\), which is then passed to both the first-stage environment and Agent~2. The learning signal therefore reflects both immediate performance and future recourse consequences. In OR/MS terms, Agent~1 learns a here-and-now policy that internalizes the expected downstream recourse value.

A one-period example illustrates the logic. Suppose Agent~1 chooses whether to reserve additional capacity, \(x\in\{0,1\}\), before demand is observed, at first-stage cost \(20x\). The demand is high with probability \(0.3\) and low otherwise. If high demand occurs and no capacity has been reserved, i.e., \(x=0\), the system incurs a second-stage shortage penalty of \(100\); if capacity has been reserved, i.e., \(x=1\), Agent~2 can use it as a recourse action, reducing the penalty to \(10\). In the low-demand scenario, no second-stage penalty is incurred. Thus, the expected cost is \(0.3(100)=30\) for \(x=0\) and \(20+0.3(10)=23\) for \(x=1\), so reserving capacity is optimal despite its immediate cost. In the two-agent DRL view, Agent~2 learns the scenario-dependent recourse decision after demand is observed, while Agent~1 learns the first-stage action that anticipates the expected value of downstream recourse. Since the underlying optimization problem is written as a cost-minimization model, the corresponding DRL reward can be interpreted as the negative of the total cost.

This interaction between Agent~1 and Agent~2 is what makes the framework especially compelling. Agent~2 learns to respond once uncertainty is observed, while Agent~1 learns to anticipate that future response when making the initial decision. The result is a policy-based approximation to the recourse mapping: rather than evaluating many scenarios through repeated optimization, the framework learns how first-stage choices shape second-stage performance and improves both policies through interaction over time. In this sense, DRL provides a complementary computational lens for stochastic programming: rather than solving a scenario-based optimization model repeatedly at deployment, it can learn stage-wise policy approximations that reflect inter-stage dependence and recourse effects \citep{YilmazBuyuktahtakin2024,Powell2019,Powell2022,sutton2018reinforcement}.

The computational motivation is a speed--quality trade-off. The two-stage DRL framework reduces the solution times from hours to fractions of seconds with reported optimality gaps of roughly \(6\%\)--\(10\%\) \citep{YilmazBuyuktahtakin2024}. These gaps are meaningful and should not be interpreted as matching the guarantees of exact stochastic programming, sample-average approximation, or problem-specific approximation algorithms. Instead, DRL is best viewed as a complementary approximation paradigm for real-time or large-scale settings, where exact solves are prohibitive, or learned policies can provide fast approximations, screening rules, or warm starts.

\subsection{DRL for Combinatorial Optimization}

Deep reinforcement learning has emerged as a promising approximation framework for combinatorial optimization, where decisions are discrete, constraints are tight, and the exhaustive search quickly becomes impractical. Early studies showed that neural policies can learn strong solution patterns for structured problems such as routing and sequencing \citep{bello2017neural,NazariEtAl2018,kool2019attention}. However, the main lesson for optimization is that performance depends less on the generic DRL algorithm than on how well the learning environment reflects the structure of the underlying model, including feasibility, constraint interactions, and the organization of the decision space \citep{bengio2021machine,MazyavkinaEtAl2021,BushajBuyuktahtakin2024}.

A particularly relevant example is \citet{BushajBuyuktahtakin2024}, who develop a problem-structured DRL framework for the multi-dimensional knapsack problem (MKP). The contribution is not simply the use of DRL on a binary integer program but the redesign of the learning environment so that it encodes the combinatorial structure of the model. The framework combines three elements: a K-means-based procedure for constructing a feasible reduced model, a heuristic variable-reordering step, and a reordered two-dimensional state representation on which several DRL algorithms---Advantage Actor--Critic (A2C), Actor--Critic using Kronecker-Factored Trust Region (ACKTR), Deep Q-Network (DQN), and Proximal Policy Optimization (PPO2)---are trained and compared. The broader point is that initialization, state design, and feasibility support can improve policy learning when they are embedded directly into the pipeline rather than left to post-processing.

The K-means component constructs a reduced but informative relaxation of the MKP by clustering constraints according to their coefficient and right-hand side patterns. Representative constraints are selected and solved with CPLEX; if the resulting solution violates omitted constraints, the most violated constraints are added back recursively until feasibility for the original model is restored. Thus, the K-means procedure is not used to produce the final solution directly but to create a feasible and informative initialization for the DRL environment.

\begin{figure}[h!]
\centering
\resizebox{0.92\linewidth}{!}{%
\begin{tikzpicture}[
    >=Latex,
    font=\footnotesize,
    block/.style={
        rectangle, draw, rounded corners=2pt,
        minimum width=2.2cm, minimum height=0.9cm, align=center
    },
    alg/.style={
        rectangle, draw, rounded corners=2pt,
        minimum width=1.7cm, minimum height=0.8cm, align=center
    },
    io/.style={
        trapezium, trapezium left angle=70, trapezium right angle=110,
        draw, minimum width=2.2cm, minimum height=1.0cm, align=center
    },
    decision/.style={
        diamond, draw, aspect=1.6, align=center, inner sep=1pt
    },
    merge/.style={
        circle, draw, minimum size=0.5cm, inner sep=0pt
    },
    line/.style={->, thick}
]

\node[io] (instances) at (0,0) {Multiple\\training\\instances};

\node[block] (kmeans)  at (3.2,1.35) {K-means\\clustering};
\node[block] (cplex)   at (6.0,1.35) {Solve reduced\\model / CPLEX};
\node[decision] (feas) at (8.7,1.35) {Feasible?};
\node[block] (viol)    at (8.7,-0.95) {Add violated\\constraints};

\node[block] (heur)    at (3.2,-1.35) {Heuristic};
\node[block] (reord)   at (6.0,-1.35) {Reorder\\variables};

\node[merge] (plus) at (10.9,0) {\(+\)};

\node[block, minimum width=3.8cm, minimum height=1.0cm] (env) at (14.4,0)
{RL environment\\(reordered initial 2D state)};

\node[alg] (a2cu)   at (12.6,2.65) {A2C\\(untrained)};
\node[alg] (acktru) at (14.4,2.65) {ACKTR\\(untrained)};
\node[alg] (dqnu)   at (16.2,2.65) {DQN\\(untrained)};
\node[alg] (ppo2u)  at (18.0,2.65) {PPO2\\(untrained)};

\node[alg] (a2ct)   at (12.6,-2.65) {A2C\\(trained)};
\node[alg] (acktrt) at (14.4,-2.65) {ACKTR\\(trained)};
\node[alg] (dqnt)   at (16.2,-2.65) {DQN\\(trained)};
\node[alg] (ppo2t)  at (18.0,-2.65) {PPO2\\(trained)};

\draw[line] (instances.east) -- ++(0.9,0) |- (kmeans.west);
\draw[line] (instances.east) -- ++(0.9,0) |- (heur.west);

\draw[line] (kmeans.east) -- (cplex.west);
\draw[line] (cplex.east) -- (feas.west);

\draw[line] (feas.south) -- node[right] {No} (viol.north);
\draw[line] (viol.north) -- ++(0,0.45) -| (cplex.south);

\draw[line] (feas.east) -- ++(0.9,0) -| node[pos=0.15, above] {Yes} (plus.north);

\draw[line] (heur.east) -- (reord.west);
\draw[line] (reord.south) -- ++(0,-0.55) -| (plus.south);

\draw[line] (plus.east) -- (env.west);

\draw[line] (a2cu.south)   -- ++(0,-0.9) -| ([xshift=-1.1cm]env.north);
\draw[line] (acktru.south) -- ++(0,-0.7) -- ([xshift=-0.35cm]env.north);
\draw[line] (dqnu.south)   -- ++(0,-0.7) -- ([xshift=0.35cm]env.north);
\draw[line] (ppo2u.south) -- ++(0,-1.15) -| ([xshift=1.75cm]env.north);

\draw[line] ([xshift=-1.1cm]env.south) -- ++(0,-0.9) -| (a2ct.north);
\draw[line] ([xshift=-0.35cm]env.south) -- ++(0,-0.7) -- (acktrt.north);
\draw[line] ([xshift=0.35cm]env.south) -- ++(0,-0.7) -- (dqnt.north);
\draw[line] ([xshift=1.1cm]env.south) -- ++(0,-0.9) -| (ppo2t.north);

\end{tikzpicture}%
}
\caption{Problem-structured DRL pipeline for the multi-dimensional knapsack problem. Training instances are processed through a K-means-based feasible relaxation path and a heuristic variable-reordering path, then transformed into a reordered two-dimensional RL environment on which candidate DRL algorithms are trained. Adapted and simplified from \citet{BushajBuyuktahtakin2024}.}
\label{fig:mkp_drl_flowchart}
\end{figure}
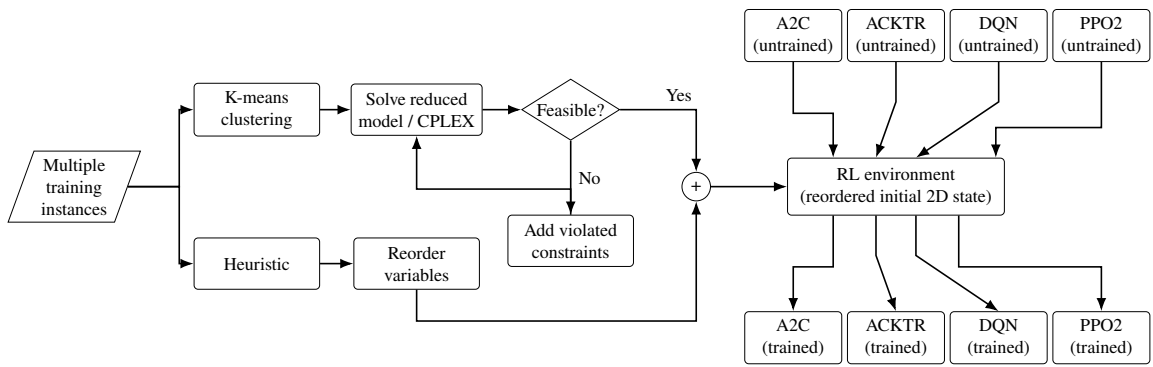

Figure~\ref{fig:mkp_drl_flowchart} summarizes the pipeline. Training instances are processed through two coordinated paths: a feasible initialization path based on K-means and a heuristic variable-reordering path. Once a feasible reduced solution is obtained, it is merged with the reordered representation to form a two-dimensional RL environment. Candidate DRL algorithms are then trained in this environment. In effect, the MKP is not presented as an unstructured sequence of independent \(0\)--\(1\) choices; instead, the environment exposes patterns in the constraint system and arranges the decision variables so that the agent can learn how local item selections affect global feasibility and objective value.

An important modeling insight is that the strongest initial heuristic is not always the best learning signal. The authors report that higher-quality classical heuristics can reduce exploration and encourage premature convergence. The K-means initialization therefore serves a more balanced role: it supplies a feasible and informative benchmark without overly constraining the agent's search behavior.

The computational payoff is substantial. \citet{BushajBuyuktahtakin2024} show that this design yields major speedups over CPLEX on medium and large MKP instances while maintaining very small optimality gaps and strong agreement with high-quality item-selection patterns. More broadly, the framework illustrates a useful principle for combinatorial optimization: DRL is most effective when used as a problem-structured policy layer, supported by informed initialization, feasibility-preserving design, and hybrid interaction with optimization logic. This points to promising directions such as DRL-guided primal heuristics, warm-start generation, feasibility-aware local improvement policies, and tighter integration with exact methods including branch-and-bound, cutting planes, and decomposition \citep{BushajBuyuktahtakin2024,bengio2021machine,MazyavkinaEtAl2021}. Across simulation-based systems, stochastic programs, and combinatorial optimization, DRL tends to be most useful as a problem-structured policy layer, with feasibility and bounds still provided by optimization.


\section{Applications and Interdisciplinary Impact}
\label{sec:applications_impact}

The practical importance of the learning--optimization interface is already evident in application domains where decisions are repeated, information unfolds over time, and operational constraints cannot be ignored. Recent work spans predict-then-optimize, contextual prescriptive analytics, decision-focused learning, neural combinatorial optimization, and reinforcement-learning-based decision systems \citep{bertsimas2020predictive,elmachtoub2022smart,elbalghiti2023generalization,kallus2023stochastic,bengio2021machine,mandi2024decision}. Despite their methodological differences, these approaches share a common principle: learning enables the transfer of experience between related instances and adaptation to changing environments, while optimization remains indispensable for enforcing feasibility, coordinating decisions, and preserving system-level performance. This perspective is also consistent with the OR/MS tradition of approximate dynamic programming and sequential decision analytics, where learned or approximate policies are ultimately evaluated by the quality of the decisions they induce \citep{Powell2019}.

\subsection{Supply Chains, Service Systems, and Procurement}

Supply chains continue to provide one of the clearest testbeds for these ideas because decisions recur over time, propagate across networks, and must be made under uncertainty. In routing and related combinatorial problems, pointer- and attention-based architectures have demonstrated that learned policies can generate high-quality solutions across broad classes of instances \citep{vinyals2015pointer,bello2017neural,kool2019attention}. In inventory and replenishment settings, deep reinforcement learning has shown considerable promise when a classical exploitable structure is difficult to leverage directly, while recent hybrid methods integrate policy learning with integer programming to manage large constrained action spaces \citep{gijsbrechts2022inventory,oroojlooyjadid2022beergame,harsha2025deeppolicy,cibaku2026vaccine}. Related advances in queueing, procurement, and sequential production--inventory systems further indicate that learned surrogates and integrated learning--optimization frameworks can support repeated operational decisions without requiring each new instance to be solved from scratch \citep{dai2022queueing,mandl2023commodity,busch2023commoditydeep,yilmaz2023lstm,yilmaz2024expandable,yilmaz2025nonanticipative}. Recent review evidence identifies supply chains as a major domain for integration of AI--optimization \citep{galande2025supplychainreview}, while related work on supply chain viability, resilience, and sustainability, including intertwined supply networks and AI-enabled intelligent digital twins, further underscores the breadth of this integration \citep{ivanov2020viability,ivanov2023idt, jozani2025proactive}.

\subsection{Healthcare and Epidemic Response}

The value of integrating learning with decision structure becomes especially evident in healthcare and response to epidemics. In clinical decision support, reinforcement learning and related sequential models must contend with partial observability, confounding, and delayed outcomes, so MDP-based formulations provide an important foundation, but not a complete solution \citep{komorowski2018ai,liu2022mlepomdp,saghafian2024adtr}. Related work on HIV prevention and treatment, epidemic resource allocation with equity considerations, integrated epidemic--supply chain planning, vaccine distribution, and other healthcare logistics settings further illustrate the breadth of this interface, spanning approximate dynamic programming, multi-stage stochastic optimization, risk-averse planning, and hybrid clustering and reinforcement-learning-based methods \citep{cosgun2018hiv,yin2021equity,yin2023ventilator,cibaku2026vaccine}. At the population level, epidemic control introduces spatial interdependence, limited resources, and the need for coordinated intervention over time, making learning more effective when it is embedded within simulation and optimization rather than used as a stand-alone policy engine \citep{BushajEtAl2023}. Related examples of integrated prediction, simulation, and optimization appear in healthcare and epidemic applications, including Ebola logistics, cancer treatment planning under spatio-temporal tumor growth, fairness-driven epidemic allocation, and data-driven COVID-19 ventilator allocation under uncertainty and risk \citep{buyuktahtakin2018ebola,kibis2019cancer,jozani2025proactive,yin2023ventilator}.

\subsection{Energy, Agriculture, and Environmental Systems}

Energy and infrastructure systems benefit similarly from hybrid approaches that couple rapid approximation with physical and operational discipline. Learning-based surrogates and data-driven decision models are increasingly used to accelerate repeated planning and control tasks in power systems and microgrids, including optimal power flow and real-time energy management applications \citep{fioretto2020predicting,pan2021deepopf,nakabi2021microgrid}. Agriculture offers an equally instructive example: irrigation, fertilization, treatment, and harvest decisions are inherently dynamic, resource-constrained, and state-dependent, and are now being explored through reinforcement-learning-based decision support \citep{gautron2022crop}. Comparable opportunities arise in environmental management, including wildfire response planning, where decisions must adapt over time while remaining consistent with ecological processes, spatial interactions, and budget limitations, reflecting a broader push toward sustainability-oriented decision making across these interconnected domains \citep{ntaimo2013simulation}.

\subsection{Robotics, UAVs, and Other Emerging Applications}

Robotics, unmanned aerial vehicles (UAVs), and autonomous delivery systems extend these questions into real-time operational environments. In such settings, learned policies must react quickly while respecting energy, synchronization, routing, and service constraints. Recent work on UAV routing and truck--drone coordination illustrates the promise of reinforcement learning and learning-based heuristics in exactly this regime \citep{fan2023uav,wu2023truckdrone}. The significance of these applications lies not in replacing optimization, but in broadening the reach of sequential decision models to environments where speed, uncertainty, and combinatorial structure interact in more immediate and consequential ways. Together, these advances show that high-impact decision making will be shaped not by learning or optimization alone but by their principled integration. They also highlight the inherently interdisciplinary nature of this interface and its growing role in designing intelligent, scalable, and implementable decision systems for complex societal challenges.



\section{Open Challenges and Research Frontiers}
\label{sec:challenges}

\subsection{Generalization, Scale, and Evaluation}

A central challenge in learning-based decision making is generalization across instances, horizons, and scales. Models may perform well on the distributions used for training, but degrade on larger, longer, or structurally different problems. This issue is especially important for learned decision-prediction maps: high in-distribution prediction accuracy does not necessarily imply reliable out-of-distribution decision quality when new instances are generated from a different problem or uncertainty distribution. This is especially pronounced in combinatorial and sequential settings, where feasibility, solution quality, and operational structure must be preserved as size, uncertainty, and planning horizons change. Unlike prediction tasks, where generalization is often measured by approximation error, decision making also requires robustness of induced decisions. Architectures such as LSTMs and transformers can capture temporal and structural dependence, but do not, by themselves, guarantee consistent performance across scales or distribution shifts. Sampling-based stochastic optimization also provides a classical foundation for estimating, comparing, and improving decisions when uncertainty must be represented through finite samples \citep{homemdemello2014monte}. Promising directions include expandable architectures, representations based on symmetry, training objectives aligned with the quality of downstream decisions, distributionally robust formulations that target uniform performance under distribution shift, and out-of-distribution stress testing combined with optimization-based repair or refinement when learned decisions deteriorate \citep{yilmaz2023lstm,yilmaz2024expandable,yilmaz2025nonanticipative,vaswani2017attention,bengio2021machine,donti2017task,wilder2019decision,duchi2021uniform, delage2010dro}. Another open direction is the design of learning proxies that jointly represent mixed continuous and discrete decisions, as coupling heterogeneous targets in a single predictive model remains difficult and can degrade prediction quality \citep{yilmaz2024expandable}.

Progress also requires a more rigorous evaluation. Learning--optimization methods are still often tested on limited or highly synthetic datasets, making meaningful comparison difficult. More informative benchmarks should assess solution quality, feasibility, computational efficiency, generalization, calibration, and robustness to distribution shift across realistic application domains \citep{bengio2021machine,donti2017task,wilder2019decision}.

\subsection{Feasibility, Reliability, and Trust}

Feasibility and reliability remain central in the deployment of learned decision models, especially in high-stakes settings. Unlike optimization models, which explicitly enforce constraints, learning-based systems may produce infeasible decisions unless constraints are carefully incorporated into the architecture or training process. This is particularly challenging when decisions involve discreteness, combinatorial structure, recourse, or explicit scenario-tree consistency. Penalty methods, projection layers, and differentiable surrogates are useful, but they may introduce approximation error or an additional computational burden. Therefore, a natural direction is to combine learned candidate solutions with optimization-based repair, certification, or refinement, so that feasibility is handled structurally rather than heuristically \citep{amos2017optnet,agrawal2019differentiable,fischetti2018deep,yilmaz2025nonanticipative}. More broadly, stronger guarantees are still needed for feasibility, robustness, risk-aware performance, and suboptimality under uncertainty. This is especially important for models that control expected tail losses through conditional value-at-risk (CVaR) and related risk-averse optimization frameworks \citep{rockafellar2000optimization,dentcheva2024risk}. Promising directions include building constraint satisfaction directly into learning architectures, using conformal prediction and chance constraints to provide probabilistic reliability or feasibility certificates for learned policies, and designing hybrid pipelines in which optimization-based repair or refinement provides post-hoc feasibility and solution-quality guarantees, especially under distributional shift \citep{angelopoulos2023conformal,bertsimas2018data,jiang2022chance,yilmaz2025nonanticipative}.

Trust is equally important. Decision makers must understand not only what is recommended, but also why it is reasonable and under what conditions it may fail. Classical optimization models often offer a degree of transparency through explicit objectives, constraints, and sensitivity information. Deep models, by contrast, may function as opaque approximators. Bridging this gap will require a combination of post-hoc explanation tools and more traditional sensitivity, scenario, and policy analysis \citep{ribeiro2016why,lundberg2017unified,bengio2021machine}. In sequential settings, trust also depends on how learned decisions behave under perturbations, tighter constraints, and alternative uncertainty realizations. Human--AI collaboration is also likely to matter for trust in practice, as evidence from field settings and large-scale meta-analysis suggests that human involvement can improve reliance and decision quality, but only when task design and division of labor align with the complementary strengths of humans and AI \citep{yang2025advisor,vaccaro2024useful}.

\subsection{Integrated Learning, Optimization, and Hybrid Decision Frameworks}

A recurring theme throughout this tutorial is that learning and optimization are often most effective when treated as interacting components rather than separate stages. The main question is not whether one should replace the other, but how the two should interact.

In predict-then-optimize settings, learning supplies inputs to an optimization model; in decision-aware learning, the optimization problem shapes the training objective. More tightly integrated approaches embed optimization layers within neural architectures or use learning to guide branching, decomposition, and search \citep{donti2017task,wilder2019decision,amos2017optnet,agrawal2019differentiable,khalil2016branch,gasse2019learn2branch,khalil2022mipgnn,lodi2017branching,bengio2021machine}. The trade-off is clear: end-to-end training can become computationally demanding, while loosely coupled pipelines may miss important feedback between prediction and decision quality. A useful path forward is a modular design in which learning provides fast approximations or structural guidance, while optimization preserves feasibility and improves the quality of the final solution \citep{bengio2021machine,yilmaz2024expandable,yilmaz2025nonanticipative}.

The same logic extends to deep reinforcement learning. Dynamic programming, approximate dynamic programming, and multi-stage stochastic programming derive decisions from explicit models of system dynamics, constraints, and uncertainty \citep{bellman1957dynamic,powell2007approximate,shapiro2009lectures}. By contrast, DRL learns policies through interaction and is particularly appealing when system dynamics are complex, simulation-based, or only partially known \citep{sutton2018reinforcement,mnih2015human,williams1992simple,konda2000actor}. A properly defined MDP or POMDP policy is non-anticipative when its state or belief state contains only information available at the decision epoch; the main OR/MS challenge is that DRL does not automatically enforce feasibility, integrality, explicit scenario-tree consistency, or recourse structure. This limitation has motivated hybrid frameworks in which learning provides candidate actions, value approximations, or policy guidance, while optimization maintains structural consistency and enhances reliability \citep{bengio2021machine,fischetti2018deep,khalil2016branch,lodi2017branching,yilmaz2025nonanticipative}.

\subsection{Uncertainty, Robustness, and Risk in Decision Making}

Point predictions are rarely sufficient for sequential decisions under uncertainty. The
learned estimates of demands, transitions, scenario probabilities, or recourse behavior can be augmented
with uncertainty information and incorporated into downstream decisions. Bayesian approximations,
including dropout-based methods, as well as deep ensembles and conformal prediction, offer practical routes to
quantify this uncertainty \citep{gal2016dropout,lakshminarayanan2017simple,angelopoulos2023conformal,
kendall2017uncertainties}. These estimates become decision-relevant once translated into
scenarios, uncertainty sets, or ambiguity sets within stochastic, robust, and distributionally
robust optimization (DRO) models \citep{bental1998robust,bertsimas2018data,delage2010dro,
wiesemann2014dro,esfahani2018wasserstein}, including distributionally robust Markov decision
processes that embed ambiguity directly in the dynamic model \citep{xu2012drmdp}.

A complementary view treats uncertainty as a belief updated over the course of the decision,
not as a quantity estimated once. Bayesian sequential decision making updates beliefs as
information arrives and uses them to guide future actions, as in Bayesian optimization under
costly evaluations \citep{jones1998efficient,frazier2018tutorial} and Bayesian graph traversal
for routing under uncertainty \citep{caballero2025bayesian}. Adaptive dynamic programming
(ADP) in the control and reinforcement learning literature offers a related perspective on
sequential learning, in which approximate value functions and policies are learned from data
or interaction when system dynamics are only partially known \citep{liu2021adaptive}. Although
related to the OR/MS tradition of approximate dynamic programming \citep{powell2007approximate},
this literature places particular emphasis on online learning, stability, and policy approximation.
Bridging these perspectives with large-scale stochastic, constrained, and mixed-integer OR/MS
models remains an important open direction. Beyond learning policies or value functions, adaptivity can also be embedded directly in the structure of the optimization model.

A further frontier lies in unifying the structural notions of adaptivity that are now emerging in stochastic
programming: optimizing \emph{when} a decision is revised \citep{basciftci2024adaptive} and
adapting the \emph{stage structure itself} to a randomly evolving number of stages
\citep{siddig2025multistage}. Combining these ideas with learning-enhanced value-function approximation could yield sequential decision systems that adapt their policies, their commitment timing, and their stage structure as information unfolds.

Related themes also appear in CS and AI, where distributional robustness is connected to adversarial
training, worst-group performance, and distribution shift \citep{sinha2018certifying,
duchi2021uniform,sagawa2020groupdro}, with recent work extending it to heterogeneous
subpopulations and latent covariate shift \citep{duchi2023latentdro}. A natural next step is
connecting uncertainty-aware learning more directly with risk-sensitive and distributionally
robust decision making, so learned models inform reliability under model error, rare events,
and distribution shift---not only expected performance.

\subsection{Emerging AI Frontiers in OR/MS Decision Making}

Some of the most promising opportunities now lie across disciplinary boundaries. Many important decision problems combine sequential uncertainty, heterogeneous data, human behavior, physical and institutional constraints, and multiscale dynamics in ways that no single methodology can adequately capture. Examples include resilient supply networks, adaptive healthcare and epidemic planning, precision agriculture and food systems, ecological and environmental management, energy-grid resilience, autonomous logistics, and decision support in biological and biomedical systems. These interfaces also extend naturally to strategic and public-policy settings: recent work on principal--agent games, for example, integrates optimization, machine learning, and game theory to develop interpretable heuristics that support fast and explainable policy decision making in forest health management \citep{baswapuram2026principalagent}.

A closely related frontier concerns how recent AI models can expand the way optimization is formulated, explored, solved, and ultimately used. Beyond direct learning of decision policies, large language models and related generative approaches are beginning to support natural-language optimization modeling, heuristic generation, algorithmic discovery for hard optimization problems, interactive planning, and model diagnosis or repair \citep{ahmaditeshnizi2024optimus,cetinkaya2026llmheuristics,romeraparedes2024funsearch,venkatachalam2025llmnetwork,ao2026optirepair,simchilevi2025llm}. FunSearch illustrates this direction by leveraging an LLM to propose executable algorithmic components that are evaluated by an external scoring procedure and iteratively improved \citep{romeraparedes2024funsearch}. Recent LLM-assisted heuristic discovery for mixed-integer single-machine scheduling provides a closely related OR/MS example: human--LLM collaboration produced scalable heuristics that improve upon simple human-designed dispatching rules for the NP-hard single-machine total tardiness problem, with performance evaluated through optimality gaps, solution times, and comparisons against MIP, dynamic programming, and classical heuristics across instances with up to 500 jobs \citep{cetinkaya2026llmheuristics}. These results suggest that LLMs may help discover scalable heuristic rules for hard constrained optimization problems, but they also highlight the need for rigorous verification, benchmarking, and solver-based evaluation rather than accepting generated rules, formulations, or explanations as valid by construction, especially because LLM outputs may contain hallucinated assumptions or algorithmic steps. This points to a broader role for AI: not merely as a predictive engine, but as an innovation tool and an interface between data, models, algorithms, and decision makers.

Other emerging interfaces reinforce the same theme. Progress at the intersection of game theory and deep learning shows how deep reinforcement learning and empirical game-theoretic analysis can be combined to learn robust policies in sequential multi-agent environments \citep{lanctot2017unified}. Such developments are especially relevant when strategic interaction, adaptation, and learning unfold together over time. Recent work on incorporating safety constraints into multi-agent DRL for naval search and defense applications further illustrates how current DRL approaches can be extended to enforce feasibility explicitly through hard constraints \citep{choi2026safety}. Another emerging interface links learning and optimization with quantum computing. Hybrid quantum--classical methods, especially variational approaches such as the Quantum Approximate Optimization Algorithm (QAOA), address hard combinatorial optimization problems through parameterized quantum circuits whose parameters are iteratively tuned by a classical optimizer. In this sense, they sit naturally at the boundary of learning and optimization, although scalability, trainability, and practical advantage remain open questions \citep{blekos2024qaoa}. The greatest promise lies in hybrid decision platforms that integrate learning, optimization, generative AI, and, over time, quantum-enhanced methods within an integrated decision framework.

\section{Conclusion}
\label{Conclusion}

This tutorial presents an OR/MS-centered view of deep learning for sequential decision making under uncertainty. Its central message is straightforward: prediction alone is not enough; what matters is decision quality under constraints, recourse, uncertainty, and evolving information \citep{elmachtoub2022smart,donti2017task,wilder2019decision,bengio2021machine,yilmaz2024expandable}. From that perspective, the tutorial has brought together several strands of work that are often studied separately. At the foundational level, dynamic programming, reinforcement learning, and multi-stage stochastic programming provide complementary ways to model sequential decisions under uncertainty \citep{bellman1957dynamic,shapiro2009lectures,birge2011introduction,sutton2018reinforcement}. At the architectural level, modern neural models, including feedforward networks, recurrent models, transformers, and deep reinforcement learning, expand the range of mappings, policies, and value functions that can be learned from rich data \citep{goodfellow2016deep,vaswani2017attention,sutton2018reinforcement,BushajEtAl2023,YilmazBuyuktahtakin2024}. At the methodological level, decision-aware learning, learning-to-optimize, expandable architectures, non-anticipative learning, and hybrid DRL--optimization approaches suggest a broader design principle: learning and optimization are often most useful as interacting components of an integrated decision pipeline \citep{yilmaz2023lstm,yilmaz2024expandable,yilmaz2025nonanticipative,donti2017task,wilder2019decision}.

A consistent theme throughout the tutorial is complementarity. Deep learning contributes to flexible approximation, adaptability, and scalability. Optimization contributes the mathematical structure needed to model constraints, recourse, and uncertainty, together with decision criteria tied to system-level objectives. This complementarity is especially important in large-scale settings, where classical methods provide theoretical grounding but can be computationally demanding, while purely learned policies may be flexible but not yet reliable enough for direct deployment \citep{bellman1957dynamic,shapiro2009lectures,birge2011introduction,mnih2015human}. Viewed in this way, the role of learning is often not to replace optimization, but to strengthen it: by learning useful structure, accelerating repeated solution processes, improving time-to-good decisions, and supporting decision rules that remain compatible with solver-based refinement \citep{bengio2021machine,khalil2016branch,lodi2017branching,yilmaz2023lstm,yilmaz2024expandable,yilmaz2025nonanticipative, baswapuram2026principalagent}.

The applications reviewed here reinforce the breadth of this opportunity. Across supply chains, service systems, healthcare and epidemic response, agriculture and environmental systems, energy, robotics, and autonomous mobility, decisions must be made sequentially, under uncertainty, and under operational constraints. These are precisely the settings in which OR/MS decision models and modern learning architectures can be combined most productively. The opportunity is not simply to predict better, but to build adaptive decision systems that learn from data while remaining feasible.

At the same time, the research frontier remains substantial. Important directions include stronger generalization across instances and scales, tighter integration of uncertainty quantification and information acquisition with downstream decision making, more reliable guarantees for feasibility and robustness, and benchmark environments that better reflect real sequential settings \citep{bengio2021machine,angelopoulos2023conformal,bertsimas2018data,yilmaz2025nonanticipative}. Emerging interfaces that involve large language models, generative AI, algorithmic discovery, and formal reasoning tools may also expand the way optimization models are formulated, explored, repaired, and communicated \citep{brown2020language,ahmaditeshnizi2024optimus,cetinkaya2026llmheuristics}.

More broadly, the field appears to be moving from predictive AI toward decision-capable AI. In that transition, OR/MS has much to contribute: not only optimization methods, but also rigorous analytical standards for defining and evaluating high-quality decisions in terms of optimality, feasibility, robustness, and implementability. From this perspective, deep learning for sequential decision making under uncertainty is best understood not as a departure from OR/MS, but as a natural extension of its core mission: to design decision systems that are analytically grounded, computationally effective, and operationally useful.

\section*{Acknowledgments}

The author sincerely thanks the four anonymous reviewers and the three INFORMS TutORials editors for their thoughtful comments and constructive guidance, which improved the clarity, scope, and exposition of this tutorial. The author also gratefully acknowledges doctoral students, collaborators, and the broader scholarly community whose work helped inform the examples and perspectives discussed in this tutorial.

\bibliographystyle{informs2014}
\bibliography{references}

@article{absi2019,
  author  = {Absi, Nabil and van den Heuvel, Wilco},
  title   = {Worst-case analysis of relax and fix heuristics for
             lot-sizing problems},
  journal = {European Journal of Operational Research},
  volume  = {279},
  number  = {2},
  pages   = {449--458},
  year    = {2019},
  doi     = {10.1016/j.ejor.2019.06.010}
}

@inproceedings{agrawal2019differentiable,
  author    = {Akshay Agrawal and Brandon Amos and Shane Barratt and
               Stephen Boyd and Steven Diamond and J. Zico Kolter},
  title     = {Differentiable Convex Optimization Layers},
  booktitle = {Advances in Neural Information Processing Systems},
  volume    = {32},
  pages     = {1--11},
  year      = {2019},
  publisher = {Curran Associates, Inc.},
  url       = {https://proceedings.neurips.cc/paper/2019/hash/
               9d16f2d4-Abstract.html},
  doi       = {10.48550/arXiv.1910.12430}
}

@inproceedings{ahmaditeshnizi2024optimus,
  author    = {Ali Ahmaditeshnizi and Wenzhi Gao and Madeleine Udell},
  title     = {{OptiMUS}: Scalable Optimization Modeling with
               ({MI}){LP} Solvers and Large Language Models},
  booktitle = {Proceedings of the 41st International Conference on
               Machine Learning},
  series    = {Proceedings of Machine Learning Research},
  volume    = {235},
  pages     = {577--596},
  year      = {2024},
  publisher = {PMLR},
  url       = {https://proceedings.mlr.press/v235/ahmaditeshnizi24a.html}
}

@inproceedings{amos2017optnet,
  author    = {Brandon Amos and J. Zico Kolter},
  title     = {{OptNet}: Differentiable Optimization as a Layer in
               Neural Networks},
  booktitle = {Proceedings of the 34th International Conference on
               Machine Learning},
  series    = {Proceedings of Machine Learning Research},
  volume    = {70},
  pages     = {136--145},
  year      = {2017},
  publisher = {PMLR},
  doi       = {10.48550/arXiv.1703.00443},
  url       = {https://arxiv.org/abs/1703.00443}
}

@misc{galande2025supplychainreview,
  author       = {Galande, Nipun and Jozani, Kimiya Mohammadi and
                  B{\"u}y{\"u}ktahtak{\i}n, {\.I}. Esra},
  title        = {Artificial Intelligence in Supply Chain
                  Optimization: A Systematic Review of Machine
                  Learning Models, Methods, and Applications},
  howpublished = {Optimization Online},
  pages        = {1--66},
  year         = {2025},
  month        = dec,
  note         = {Published online December 8, 2025},
  url          = { https://optimization-online.org/?p=32864}
}

@article{anderson2020strong,
  author  = {Anderson, Ross and Huchette, Joey and Ma, Will and
             Tjandraatmadja, Christian and Vielma, Juan Pablo},
  title   = {Strong Mixed-Integer Programming Formulations for
             Trained Neural Networks},
  journal = {Mathematical Programming},
  volume  = {183},
  number  = {1},
  pages   = {3--39},
  year    = {2020},
  doi     = {10.1007/s10107-020-01474-5}
}

@article{angelopoulos2023conformal,
  author  = {Anastasios N. Angelopoulos and Stephen Bates},
  title   = {Conformal Prediction: A Gentle Introduction},
  journal = {Foundations and Trends in Machine Learning},
  volume  = {16},
  number  = {4},
  pages   = {494--591},
  year    = {2023},
  doi     = {10.1561/2200000101}
}

@misc{ao2026optirepair,
  author       = {Ruicheng Ao and David Simchi-Levi and
                  Xinshang Wang},
  title        = {{OptiRepair}: Closed-Loop Diagnosis and Repair of
                  Supply Chain Optimization Models with {LLM} Agents},
  howpublished = {arXiv preprint arXiv:2602.19439},
  year         = {2026},
  doi          = {10.48550/arXiv.2602.19439},
  url          = {https://arxiv.org/abs/2602.19439}
}

@inproceedings{bahdanau2015neural,
  author    = {Bahdanau, Dzmitry and Cho, Kyunghyun and
               Bengio, Yoshua},
  title     = {Neural Machine Translation by Jointly Learning to
               Align and Translate},
  booktitle = {International Conference on Learning Representations
               ({ICLR} 2015)},
  year      = {2015},
  doi       = {10.48550/arXiv.1409.0473},
  url       = {https://arxiv.org/abs/1409.0473}
}

@book{bellman1957dynamic,
  author    = {Richard Bellman},
  title     = {Dynamic Programming},
  year      = {1957},
  publisher = {Princeton University Press},
  address   = {Princeton, NJ},
  isbn      = {9780691146683},
  url       = {https://press.princeton.edu/books/paperback/9780691146683/dynamic-programming}
}

@inproceedings{bello2017neural,
  author    = {Bello, Irwan and Pham, Hieu and Le, Quoc V. and
               Norouzi, Mohammad and Bengio, Samy},
  title     = {Neural Combinatorial Optimization with Reinforcement
               Learning},
  booktitle = {International Conference on Learning Representations
               ({ICLR} 2017)},
  year      = {2017},
  doi       = {10.48550/arXiv.1611.09940},
  url       = {https://openreview.net/forum?id=Bk9mxlSFx}
}

@article{bengio1994learning,
  author  = {Bengio, Yoshua and Simard, Patrice and Frasconi, Paolo},
  title   = {Learning long-term dependencies with gradient descent
             is difficult},
  journal = {IEEE Transactions on Neural Networks},
  volume  = {5},
  number  = {2},
  pages   = {157--166},
  year    = {1994},
  doi     = {10.1109/72.279181}
}

@article{bengio2021machine,
  author  = {Yoshua Bengio and Andrea Lodi and Antoine Prouvost},
  title   = {Machine Learning for Combinatorial Optimization:
             A Methodological Tour d'Horizon},
  journal = {European Journal of Operational Research},
  volume  = {290},
  number  = {2},
  pages   = {405--421},
  year    = {2021},
  doi     = {10.1016/j.ejor.2020.07.063}
}

@article{bental1998robust,
  author  = {Aharon Ben-Tal and Arkadi Nemirovski},
  title   = {Robust Convex Optimization},
  journal = {Mathematics of Operations Research},
  volume  = {23},
  number  = {4},
  pages   = {769--805},
  year    = {1998},
  doi     = {10.1287/moor.23.4.769}
}

@book{bertsekas1995dynamic,
  author    = {Dimitri P. Bertsekas},
  title     = {Dynamic Programming and Optimal Control},
  year      = {1995},
  publisher = {Athena Scientific},
  address   = {Belmont, MA},
  volume    = {1},
  isbn      = {9781886529434},
  url       = {https://www.athenasc.com/dpcontents.html}
}

@book{bertsekas1996neuro,
  author    = {Bertsekas, Dimitri P. and Tsitsiklis, John N.},
  title     = {Neuro-Dynamic Programming},
  publisher = {Athena Scientific},
  address   = {Belmont, MA},
  year      = {1996},
  isbn      = {978-1-886529-10-6},
  url       = {https://athenasc.com/ndpbook.html}
}

@article{bertsimas2002,
  author  = {Bertsimas, Dimitris and Demir, Ramazan},
  title   = {An Approximate Dynamic Programming Approach to
             Multidimensional Knapsack Problems},
  journal = {Management Science},
  volume  = {48},
  number  = {4},
  pages   = {550--565},
  year    = {2002},
  doi     = {10.1287/mnsc.48.4.550.208}
}

@article{bertsimas2018data,
  author  = {Dimitris Bertsimas and Vishal Gupta and Nathan Kallus},
  title   = {Data-Driven Robust Optimization},
  journal = {Mathematical Programming},
  volume  = {167},
  number  = {2},
  pages   = {235--292},
  year    = {2018},
  doi     = {10.1007/s10107-017-1125-8}
}

@article{bertsimas2020predictive,
  author  = {Bertsimas, Dimitris and Kallus, Nathan},
  title   = {From Predictive to Prescriptive Analytics},
  journal = {Management Science},
  volume  = {66},
  number  = {3},
  pages   = {1025--1044},
  year    = {2020},
  doi     = {10.1287/mnsc.2018.3253}
}

@book{birge2011introduction,
  author    = {John R. Birge and Fran{\c{c}}ois Louveaux},
  title     = {Introduction to Stochastic Programming},
  year      = {2011},
  publisher = {Springer},
  address   = {New York, NY},
  edition   = {2},
  doi       = {10.1007/978-1-4614-0237-4}
}

@inproceedings{brown2020language,
  author    = {Tom B. Brown and Benjamin Mann and Nick Ryder and
               Melanie Subbiah and Jared D. Kaplan and
               Prafulla Dhariwal and Arvind Neelakantan and
               Pranav Shyam and Girish Sastry and Amanda Askell and
               Sandhini Agarwal and Ariel Herbert-Voss and
               Gretchen Krueger and Tom Henighan and Rewon Child and
               Aditya Ramesh and Daniel M. Ziegler and Jeffrey Wu and
               Clemens Winter and Christopher Hesse and Mark Chen and
               Eric Sigler and Mateusz Litwin and Scott Gray and
               Benjamin Chess and Jack Clark and Christopher Berner and
               Sam McCandlish and Alec Radford and Ilya Sutskever and
               Dario Amodei},
  title     = {Language Models are Few-Shot Learners},
  booktitle = {Advances in Neural Information Processing Systems},
  volume    = {33},
  pages     = {1877--1901},
  year      = {2020},
  publisher = {Curran Associates, Inc.},
  doi       = {10.48550/arXiv.2005.14165},
  url       = {https://proceedings.neurips.cc/paper/2020/hash/
               1457c0d6bfcb4967418bfb8ac142f64a-Abstract.html}
}

@article{BushajBuyuktahtakin2024,
  author  = {Sabah Bushaj and {\.I}. Esra
             B{\"u}y{\"u}ktahtak{\i}n},
  title   = {A {K}-means Supported Reinforcement Learning Framework
             to Multi-dimensional Knapsack},
  journal = {Journal of Global Optimization},
  volume  = {89},
  number  = {3},
  pages   = {655--685},
  year    = {2024},
  doi     = {10.1007/s10898-024-01364-6}
}

@article{BushajEtAl2023,
  author  = {Sabah Bushaj and Xuecheng Yin and Arjeta Beqiri and
             Donald Andrews and {\.I}. Esra
             B{\"u}y{\"u}ktahtak{\i}n},
  title   = {A Simulation-Deep Reinforcement Learning ({SiRL})
             Approach for Epidemic Control Optimization},
  journal = {Annals of Operations Research},
  volume  = {328},
  number  = {1},
  pages   = {245--277},
  year    = {2023},
  doi     = {10.1007/s10479-022-04926-7}
}

@article{buyuktahtakin2022stage,
  author  = {B{\"u}y{\"u}ktahtak{\i}n, {\.I}. Esra},
  title   = {Stage-$t$ Scenario Dominance for Risk-Averse
             Multi-Stage Stochastic Mixed-Integer Programs},
  journal = {Annals of Operations Research},
  volume  = {309},
  pages   = {1--35},
  year    = {2022},
  doi     = {10.1007/s10479-021-04388-3}
}

@article{cetinkaya2026llmheuristics,
  author  = {{\c{C}}etinkaya, {\.{I}}brahim O{\u{g}}uz and
             B{\"u}y{\"u}ktahtak{\i}n, {\.{I}}. Esra and
             Shojaee, Parshin and Reddy, Chandan K.},
  title   = {Discovering Heuristics with Large Language Models
             ({LLMs}) for Mixed-Integer Programs:
             Single-Machine Scheduling},
  journal = {Computers \& Operations Research},
  volume  = {186},
  pages   = {107325},
  year    = {2026},
  doi     = {10.1016/j.cor.2025.107325}
}

@article{chen2024learning,
  author  = {Chen, Xiaohan and Liu, Jialin and Yin, Wotao},
  title   = {Learning to Optimize: A Tutorial for Continuous and
             Mixed-Integer Optimization},
  journal = {Science China Mathematics},
  volume  = {67},
  number  = {6},
  pages   = {1191--1262},
  year    = {2024},
  doi     = {10.1007/s11425-023-2293-3}
}

@article{dai2022queueing,
  author  = {J. G. Dai and Mark Gluzman},
  title   = {Queueing Network Controls via Deep Reinforcement
             Learning},
  journal = {Stochastic Systems},
  volume  = {12},
  number  = {1},
  pages   = {30--67},
  year    = {2022},
  doi     = {10.1287/stsy.2021.0081}
}

@misc{ding2019sddip,
  author       = {Ding, Lian and Ahmed, Shabbir and Shapiro, Alexander},
  title        = {A {Python} Package for Multi-Stage Stochastic
                  Programming},
  howpublished = {Optimization Online},
  year         = {2019},
  url          = {https://optimization-online.org/2019/05/7199/},
  note         = {Technical report, pp.\ 1--41}
}

@inproceedings{donti2017task,
  author    = {Priya L. Donti and Brandon Amos and J. Zico Kolter},
  title     = {Task-based End-to-end Model Learning in Stochastic
               Optimization},
  booktitle = {Advances in Neural Information Processing Systems},
  volume    = {30},
  pages     = {5484--5494},
  year      = {2017},
  publisher = {Curran Associates, Inc.},
  doi       = {10.48550/arXiv.1703.04529}
}

@article{elmachtoub2022smart,
  author  = {Adam N. Elmachtoub and Paul Grigas},
  title   = {Smart ``Predict, then Optimize''},
  journal = {Management Science},
  volume  = {68},
  number  = {1},
  pages   = {9--26},
  year    = {2022},
  doi     = {10.1287/mnsc.2020.3922}
}

@article{elman1990finding,
  author  = {Elman, Jeffrey L.},
  title   = {Finding Structure in Time},
  journal = {Cognitive Science},
  volume  = {14},
  number  = {2},
  pages   = {179--211},
  year    = {1990},
  doi     = {10.1016/0364-0213(90)90002-E}
}

@article{fan2023uav,
  author  = {Fan, Mingfeng and Wu, Yaoxin and Liao, Tianjun and
             Cao, Zhiguang and Guo, Hongliang and
             Sartoretti, Guillaume and Wu, Guohua},
  title   = {Deep Reinforcement Learning for {UAV} Routing in the
             Presence of Multiple Charging Stations},
  journal = {IEEE Transactions on Vehicular Technology},
  volume  = {72},
  number  = {5},
  pages   = {5732--5746},
  year    = {2023},
  doi     = {10.1109/TVT.2022.3232607}
}

@article{fischetti2018deep,
  author  = {Matteo Fischetti and Jason Jo},
  title   = {Deep Neural Networks and Mixed Integer Linear
             Optimization},
  journal = {Constraints},
  volume  = {23},
  number  = {3},
  pages   = {296--309},
  year    = {2018},
  doi     = {10.1007/s10601-018-9285-6}
}

@inproceedings{gal2016dropout,
  author    = {Yarin Gal and Zoubin Ghahramani},
  title     = {Dropout as a {Bayesian} Approximation: Representing
               Model Uncertainty in Deep Learning},
  booktitle = {Proceedings of the 33rd International Conference on
               Machine Learning},
  series    = {Proceedings of Machine Learning Research},
  volume    = {48},
  pages     = {1050--1059},
  year      = {2016},
  publisher = {PMLR},
  doi       = {10.48550/arXiv.1506.02142},
  url       = {https://arxiv.org/abs/1506.02142}
}

@article{gers2000learning,
  author  = {Gers, Felix A. and Schmidhuber, J{\"u}rgen and
             Cummins, Fred},
  title   = {Learning to Forget: Continual Prediction with {LSTM}},
  journal = {Neural Computation},
  volume  = {12},
  number  = {10},
  pages   = {2451--2471},
  year    = {2000},
  doi     = {10.1162/089976600300015015}
}

@article{gijsbrechts2022inventory,
  author  = {Gijsbrechts, Joren and Boute, Robert N. and
             Van Mieghem, Jan A. and Zhang, Dennis J.},
  title   = {Can Deep Reinforcement Learning Improve Inventory
             Management? {P}erformance on Lost Sales,
             Dual-Sourcing, and Multi-Echelon Problems},
  journal = {Manufacturing \& Service Operations Management},
  volume  = {24},
  number  = {3},
  pages   = {1349--1368},
  year    = {2022},
  doi     = {10.1287/msom.2021.1064}
}

@book{goodfellow2016deep,
  author    = {Ian Goodfellow and Yoshua Bengio and Aaron Courville},
  title     = {Deep Learning},
  year      = {2016},
  publisher = {MIT Press},
  address   = {Cambridge, MA},
  isbn      = {9780262035613},
  url       = {https://www.deeplearningbook.org/}
}

@inproceedings{hamilton2017inductive,
  author    = {Hamilton, William L. and Ying, Rex and Leskovec, Jure},
  title     = {Inductive Representation Learning on Large Graphs},
  booktitle = {Advances in Neural Information Processing Systems},
  volume    = {30},
  year      = {2017},
  publisher = {Curran Associates, Inc.},
  doi       = {10.48550/arXiv.1706.02216},
  url       = {https://arxiv.org/abs/1706.02216}
}

@article{hochreiter1997long,
  author    = {Hochreiter, Sepp and Schmidhuber, J{\"u}rgen},
  title     = {Long Short-Term Memory},
  journal   = {Neural Computation},
  volume    = {9},
  number    = {8},
  pages     = {1735--1780},
  year      = {1997},
  publisher = {MIT Press},
  doi       = {10.1162/neco.1997.9.8.1735}
}

@article{hornik1989multilayer,
  author  = {Hornik, Kurt and Stinchcombe, Maxwell and White, Halbert},
  title   = {Multilayer Feedforward Networks are Universal
             Approximators},
  journal = {Neural Networks},
  volume  = {2},
  number  = {5},
  pages   = {359--366},
  year    = {1989},
  doi     = {10.1016/0893-6080(89)90020-8}
}

@inproceedings{kendall2017uncertainties,
  author    = {Alex Kendall and Yarin Gal},
  title     = {What Uncertainties Do We Need in {Bayesian} Deep
               Learning for Computer Vision?},
  booktitle = {Advances in Neural Information Processing Systems},
  volume    = {30},
  year      = {2017},
  publisher = {Curran Associates, Inc.},
  doi       = {10.48550/arXiv.1703.04977},
  url       = {https://arxiv.org/abs/1703.04977}
}

@article{KerrEtAl2021,
  author  = {Cliff C. Kerr and Robyn M. Stuart and Dina Mistry and
             Romesh G. Abeysuriya and Katherine Rosenfeld and
             Gregory R. Hart and Rafael C. Nu{\~n}ez and
             Jamie A. Cohen and Prashanth Selvaraj and
             Brittany Hagedorn and Lauren George and
             Michaela Jastrz{\k{e}}bska and Alastair Izzo and
             Greer Fowler and Anna Palmer and Dominic Delport and
             Nicholas Scott and Sherrie Kelly and
             Cristina S. Bennette and Bilal Wagner and
             Stewart T. Chang and Assaf Vassall and
             Bradley J. Pearson and Peter H. Winskill and
             Amanda Panovska-Griffiths and Martin Famulare and
             Daniel J. Klein},
  title   = {Covasim: An Agent-Based Model of {COVID}-19 Dynamics
             and Interventions},
  journal = {PLoS Computational Biology},
  volume  = {17},
  number  = {7},
  pages   = {e1009149},
  year    = {2021},
  doi     = {10.1371/journal.pcbi.1009149}
}

@inproceedings{khalil2017learning,
  author    = {Khalil, Elias B. and Dai, Hanjun and Zhang, Yuyu and
               Dilkina, Bistra and Song, Le},
  title     = {Learning Combinatorial Optimization Algorithms over
               Graphs},
  booktitle = {Advances in Neural Information Processing Systems},
  volume    = {30},
  year      = {2017},
  publisher = {Curran Associates, Inc.},
  doi       = {10.48550/arXiv.1704.01665},
  url       = {https://arxiv.org/abs/1704.01665}
}

@inproceedings{khalil2016branch,
  author    = {Khalil, Elias B. and {Le Bodic}, Pierre and Song, Le
               and Nemhauser, George and Dilkina, Bistra},
  title     = {Learning to Branch in Mixed Integer Programming},
  booktitle = {Proceedings of the {AAAI} Conference on Artificial
               Intelligence},
  volume    = {30},
  number    = {1},
  pages     = {724--731},
  year      = {2016},
  doi       = {10.1609/aaai.v30i1.10080}
}

@inproceedings{kipf2017semi,
  author    = {Kipf, Thomas N. and Welling, Max},
  title     = {Semi-Supervised Classification with Graph
               Convolutional Networks},
  booktitle = {International Conference on Learning Representations
               ({ICLR} 2017)},
  year      = {2017},
  doi       = {10.48550/arXiv.1609.02907},
  url       = {https://openreview.net/forum?id=SJU4ayYgl}
}

@article{komorowski2018ai,
  author  = {Komorowski, Matthieu and Celi, Leo Anthony and
             Badawi, Omar and Gordon, Anthony C. and
             Faisal, Aldo A.},
  title   = {The Artificial Intelligence Clinician Learns Optimal
             Treatment Strategies for Sepsis in Intensive Care},
  journal = {Nature Medicine},
  volume  = {24},
  pages   = {1716--1720},
  year    = {2018},
  doi     = {10.1038/s41591-018-0213-5}
}

@inproceedings{konda2000actor,
  author    = {Vijay R. Konda and John N. Tsitsiklis},
  title     = {Actor-Critic Algorithms},
  booktitle = {Advances in Neural Information Processing Systems},
  volume    = {12},
  pages     = {1008--1014},
  year      = {2000},
  publisher = {MIT Press},
url       = {https://proceedings.neurips.cc/paper_files/paper/1999/hash/6449f44a102fde848669bdd9eb6b76fa-Abstract.html}
}

@inproceedings{kool2019attention,
  author    = {Kool, Wouter and van Hoof, Herke and Welling, Max},
  title     = {Attention, Learn to Solve Routing Problems!},
  booktitle = {International Conference on Learning Representations
               ({ICLR} 2019)},
  year      = {2019},
  doi       = {10.48550/arXiv.1803.08475},
  url       = {https://openreview.net/forum?id=ByxBFsRqYm}
}

@inproceedings{lakshminarayanan2017simple,
  author    = {Balaji Lakshminarayanan and Alexander Pritzel and
               Charles Blundell},
  title     = {Simple and Scalable Predictive Uncertainty Estimation
               Using Deep Ensembles},
  booktitle = {Advances in Neural Information Processing Systems},
  volume    = {30},
  pages     = {6402--6413},
  year      = {2017},
  publisher = {Curran Associates, Inc.},
  doi       = {10.48550/arXiv.1612.01474},
  url       = {https://arxiv.org/abs/1612.01474}
}

@article{lecun2015deep,
  author  = {LeCun, Yann and Bengio, Yoshua and Hinton, Geoffrey},
  title   = {Deep Learning},
  journal = {Nature},
  volume  = {521},
  number  = {7553},
  pages   = {436--444},
  year    = {2015},
  doi     = {10.1038/nature14539}
}

@inproceedings{LibinEtAl2021,
  author    = {Libin, Pieter and Moonens, Arno and Verstraeten, Timothy and {Perez-Sanjines}, Fabian Ramiro and Hens, Niel and Lemey, Philippe and Nowe, Ann},
  editor    = {Dong, Yuxiao and Ifrim, Georgiana and Mladeni{\'c}, Dunja and Saunders, Craig and {Van Hoecke}, Sofie},
  title     = {Deep Reinforcement Learning for Large-Scale Epidemic Control},
  booktitle = {Machine Learning and Knowledge Discovery in Databases. Applied Data Science and Demo Track},
  series    = {Lecture Notes in Computer Science},
  volume    = {12461},
  pages     = {155--170},
  publisher = {Springer},
  year      = {2021},
  doi       = {10.1007/978-3-030-67670-4_10},
  url       = {https://doi.org/10.1007/978-3-030-67670-4_10}
}

@article{lodi2017branching,
  author  = {Lodi, Andrea and Zarpellon, Giulia},
  title   = {On Learning and Branching: A Survey},
  journal = {TOP},
  volume  = {25},
  number  = {2},
  pages   = {207--236},
  year    = {2017},
  doi     = {10.1007/s11750-017-0451-6}
}

@inproceedings{lundberg2017unified,
  author    = {Scott M. Lundberg and Su-In Lee},
  title     = {A Unified Approach to Interpreting Model Predictions},
  booktitle = {Advances in Neural Information Processing Systems},
  volume    = {30},
  pages     = {4765--4774},
  year      = {2017},
  publisher = {Curran Associates, Inc.},
  doi       = {10.48550/arXiv.1705.07874},
  url       = {https://arxiv.org/abs/1705.07874}
}

@inproceedings{luong2015effective,
  author    = {Luong, Minh-Thang and Pham, Hieu and
               Manning, Christopher D.},
  title     = {Effective Approaches to Attention-based Neural
               Machine Translation},
  booktitle = {Proceedings of the 2015 Conference on Empirical
               Methods in Natural Language Processing},
  pages     = {1412--1421},
  year      = {2015},
  publisher = {Association for Computational Linguistics},
  address   = {Lisbon, Portugal},
  doi       = {10.18653/v1/D15-1166}
}

@article{mandi2024decision,
  author    = {Mandi, Jayanta and Kotary, James and Berden, Senne
               and Mulamba, Maxime and Bucarey, V{\'i}ctor
               and Guns, Tias and Fioretto, Ferdinando},
  title     = {Decision-focused learning: Foundations, state of the
               art, benchmark and future opportunities},
  journal   = {Journal of Artificial Intelligence Research},
  volume    = {80},
  pages     = {1623--1701},
  year      = {2024},
  doi       = {10.1613/jair.1.15320},
  url       = {http://dx.doi.org/10.1613/jair.1.15320}
}

@article{MazyavkinaEtAl2021,
  author  = {Nina Mazyavkina and Sergey Sviridov and
             Sergei Ivanov and Evgeny Burnaev},
  title   = {Reinforcement Learning for Combinatorial
             Optimization: A Survey},
  journal = {Computers \& Operations Research},
  volume  = {134},
  pages   = {105400},
  year    = {2021},
  doi     = {10.1016/j.cor.2021.105400}
}

@article{mnih2015human,
  author  = {Volodymyr Mnih and Koray Kavukcuoglu and David Silver
             and Andrei A. Rusu and Joel Veness and
             Marc G. Bellemare and Alex Graves and
             Martin Riedmiller and Andreas K. Fidjeland and
             Georg Ostrovski and Stig Petersen and
             Charles Beattie and Amir Sadik and
             Ioannis Antonoglou and Helen King and
             Dharshan Kumaran and Daan Wierstra and
             Shane Legg and Demis Hassabis},
  title   = {Human-Level Control through Deep Reinforcement
             Learning},
  journal = {Nature},
  volume  = {518},
  number  = {7540},
  pages   = {529--533},
  year    = {2015},
  doi     = {10.1038/nature14236}
}

@article{nakabi2021microgrid,
  author  = {Nakabi, Taha Abdelhalim and Toivanen, Pekka},
  title   = {Deep Reinforcement Learning for Energy Management in
             a Microgrid with Flexible Demand},
  journal = {Sustainable Energy, Grids and Networks},
  volume  = {25},
  pages   = {100413},
  year    = {2021},
  doi     = {10.1016/j.segan.2020.100413}
}

@inproceedings{NazariEtAl2018,
  author    = {Mohammadreza Nazari and Afshin Oroojlooy and
               Lawrence V. Snyder and Martin Tak{\'a}{\v{c}}},
  title     = {Reinforcement Learning for Solving the Vehicle
               Routing Problem},
  booktitle = {Advances in Neural Information Processing Systems},
  volume    = {31},
  pages     = {9839--9849},
  year      = {2018},
  publisher = {Curran Associates, Inc.},
  doi       = {10.48550/arXiv.1802.04240},
  url       = {https://arxiv.org/abs/1802.04240}
}

@book{powell2007approximate,
  author    = {Warren B. Powell},
  title     = {Approximate Dynamic Programming: Solving the
               Curses of Dimensionality},
  year      = {2011},
  publisher = {John Wiley \& Sons},
  address   = {Hoboken, NJ},
  edition   = {2},
  isbn      = {9780470171554},
  doi       = {10.1002/9781118029176},
  url       = {https://onlinelibrary.wiley.com/doi/book/10.1002/9781118029176}
}

@article{Powell2019,
  author  = {Warren B. Powell},
  title   = {A Unified Framework for Stochastic Optimization},
  journal = {European Journal of Operational Research},
  volume  = {275},
  number  = {3},
  pages   = {795--821},
  year    = {2019},
  doi     = {10.1016/j.ejor.2018.07.014}
}

@book{Powell2022,
  author    = {Warren B. Powell},
  title     = {Reinforcement Learning and Stochastic Optimization:
               A Unified Framework for Sequential Decisions},
  year      = {2022},
  publisher = {Wiley},
  address   = {Hoboken, NJ},
  isbn      = {9781119815068},
  url       = {https://www.wiley.com/en-us/Reinforcement+Learning+
               and+Stochastic+Optimization-p-9781119815037}
}

@book{puterman1994markov,
  author    = {Martin L. Puterman},
  title     = {Markov Decision Processes: Discrete Stochastic
               Dynamic Programming},
  year      = {1994},
  publisher = {John Wiley \& Sons},
  address   = {New York},
  doi       = {10.1002/9780470316887}
}

@inproceedings{ribeiro2016why,
  author    = {Marco T{\'u}lio Ribeiro and Sameer Singh and
               Carlos Guestrin},
  title     = {``{W}hy Should {I} Trust You?'': Explaining the
               Predictions of Any Classifier},
  booktitle = {Proceedings of the 22nd {ACM} {SIGKDD} International
               Conference on Knowledge Discovery and Data Mining},
  pages     = {1135--1144},
  year      = {2016},
  doi       = {10.1145/2939672.2939778}
}

@book{rockafellar2009variational,
  author    = {Rockafellar, R. Tyrrell and Wets, Roger J.-B.},
  title     = {Variational Analysis},
  year      = {1998},
  publisher = {Springer},
  address   = {Berlin, Heidelberg},
  doi       = {10.1007/978-3-642-02431-3}
}

@article{rumelhart1986learning,
  author  = {Rumelhart, David E. and Hinton, Geoffrey E. and
             Williams, Ronald J.},
  title   = {Learning Representations by Back-Propagating Errors},
  journal = {Nature},
  volume  = {323},
  number  = {6088},
  pages   = {533--536},
  year    = {1986},
  doi     = {10.1038/323533a0}
}

@article{scarselli2009graph,
  author  = {Scarselli, Franco and Gori, Marco and
             Tsoi, Ah Chung and Hagenbuchner, Markus and
             Monfardini, Gabriele},
  title   = {The Graph Neural Network Model},
  journal = {IEEE Transactions on Neural Networks},
  volume  = {20},
  number  = {1},
  pages   = {61--80},
  year    = {2009},
  doi     = {10.1109/TNN.2008.2005605}
}

@book{shapiro2009lectures,
  author    = {Alexander Shapiro and Darinka Dentcheva and
               Andrzej Ruszczy{\'n}ski},
  title     = {Lectures on Stochastic Programming: Modeling
               and Theory},
  year      = {2009},
  publisher = {Society for Industrial and Applied Mathematics
               and Mathematical Programming Society},
  address   = {Philadelphia, PA},
  doi       = {10.1137/1.9780898718751}
}

@misc{simchilevi2025llm,
  author       = {David Simchi-Levi and Konstantina Mellou and
                  Ishai Menache and Jeevan Pathuri},
  title        = {Large Language Models for Supply Chain Decisions},
  howpublished = {{SSRN} Electronic Journal},
  year         = {2025},
  doi          = {10.2139/ssrn.5370043},
  url          = {https://ssrn.com/abstract=5370043}
}

@inproceedings{sutskever2014sequence,
  author    = {Sutskever, Ilya and Vinyals, Oriol and Le, Quoc V.},
  title     = {Sequence to Sequence Learning with Neural Networks},
  booktitle = {Advances in Neural Information Processing Systems},
  volume    = {27},
  pages     = {3104--3112},
  year      = {2014},
  publisher = {Curran Associates, Inc.},
  doi       = {10.48550/arXiv.1409.3215},
  url       = {https://arxiv.org/abs/1409.3215}
}

@book{sutton2018reinforcement,
  author    = {Richard S. Sutton and Andrew G. Barto},
  title     = {Reinforcement Learning: An Introduction},
  year      = {2018},
  publisher = {MIT Press},
  address   = {Cambridge, MA},
  edition   = {2},
  isbn      = {9780262039246},
  url       = {http://incompleteideas.net/book/the-book-2nd.html}
}

@inproceedings{vaswani2017attention,
  author    = {Ashish Vaswani and Noam Shazeer and Niki Parmar and
               Jakob Uszkoreit and Llion Jones and
               Aidan N. Gomez and {\L}ukasz Kaiser and
               Illia Polosukhin},
  title     = {Attention Is All You Need},
  booktitle = {Advances in Neural Information Processing Systems},
  volume    = {30},
  pages     = {5998--6008},
  year      = {2017},
  publisher = {Curran Associates, Inc.},
  doi       = {10.48550/arXiv.1706.03762},
  url       = {https://arxiv.org/abs/1706.03762}
}

@misc{venkatachalam2025llmnetwork,
  author       = {Saravanan Venkatachalam},
  title        = {Integrating Large Language Models with Network
                  Optimization for Interactive and Explainable
                  Supply Chain Planning: A Real-World Case Study},
  howpublished = {arXiv preprint arXiv:2508.21622},
  year         = {2025},
  doi          = {10.48550/arXiv.2508.21622},
  url          = {https://arxiv.org/abs/2508.21622}
}

@inproceedings{vinyals2015pointer,
  author    = {Vinyals, Oriol and Fortunato, Meire and
               Jaitly, Navdeep},
  title     = {Pointer Networks},
  booktitle = {Advances in Neural Information Processing Systems},
  volume    = {28},
  pages     = {2674--2682},
  year      = {2015},
  publisher = {Curran Associates, Inc.},
  doi       = {10.48550/arXiv.1506.03134},
  url       = {https://arxiv.org/abs/1506.03134}
}

@phdthesis{watkins1989learning,
  author = {Watkins, Christopher J. C. H.},
  title  = {Learning from Delayed Rewards},
  year   = {1989},
  school = {University of Cambridge},
   url    = {http://www.cs.rhul.ac.uk/~chrisw/new_thesis.pdf}
}

@article{watkins1992q,
  author  = {Christopher J. C. H. Watkins and Peter Dayan},
  title   = {Q-learning},
  journal = {Machine Learning},
  volume  = {8},
  number  = {3--4},
  pages   = {279--292},
  year    = {1992},
  doi     = {10.1007/BF00992698}
}

@article{watson2011progressive,
  author  = {Watson, Jean-Paul and Woodruff, David L.},
  title   = {Progressive Hedging Innovations for a Class of
             Stochastic Mixed-Integer Resource Allocation Problems},
  journal = {Computational Management Science},
  volume  = {8},
  number  = {4},
  pages   = {355--370},
  year    = {2011},
  doi     = {10.1007/s10287-010-0125-4}
}

@inproceedings{wilder2019decision,
  author    = {Bryan Wilder and Bistra Dilkina and Milind Tambe},
  title     = {Melding the Data-Decisions Pipeline: Decision-Focused
               Learning for Combinatorial Optimization},
  booktitle = {Proceedings of the {AAAI} Conference on Artificial
               Intelligence},
  volume    = {33},
  number    = {1},
  pages     = {1658--1665},
  year      = {2019},
  doi       = {10.1609/aaai.v33i01.33011658}
}

@article{williams1992simple,
  author  = {Williams, Ronald J.},
  title   = {Simple Statistical Gradient-Following Algorithms for
             Connectionist Reinforcement Learning},
  journal = {Machine Learning},
  volume  = {8},
  number  = {3--4},
  pages   = {229--256},
  year    = {1992},
  doi     = {10.1007/BF00992696}
}

@article{wu2021comprehensive,
  author  = {Wu, Zonghan and Pan, Shirui and Chen, Fengwen and
             Long, Guodong and Zhang, Chengqi and Yu, Philip S.},
  title   = {A Comprehensive Survey on Graph Neural Networks},
  journal = {IEEE Transactions on Neural Networks and Learning
             Systems},
  volume  = {32},
  number  = {1},
  pages   = {4--24},
  year    = {2021},
  doi     = {10.1109/TNNLS.2020.2978386}
}

@article{wu2023truckdrone,
  author  = {Wu, Guohua and Fan, Mingfeng and Shi, Jianmai and
             Feng, Yanghe},
  title   = {Reinforcement Learning Based Truck-and-Drone
             Coordinated Delivery},
  journal = {IEEE Transactions on Artificial Intelligence},
  volume  = {4},
  number  = {4},
  pages   = {754--763},
  year    = {2023},
  doi     = {10.1109/TAI.2021.3087666}
}

@article{yilmaz2023lstm,
  author  = {Yilmaz, Do{\u{g}}a{\c{c}}an and
             B{\"u}y{\"u}ktahtak{\i}n, {\.{I}}. Esra},
  title   = {Learning Optimal Solutions via an {LSTM}-Optimization
             Framework},
  journal = {Operations Research Forum},
  volume  = {4},
  number  = {2},
  pages   = {48},
  year    = {2023},
  doi     = {10.1007/s43069-023-00224-5}
}

@article{yilmaz2024expandable,
  author  = {Yilmaz, Do{\u{g}}a{\c{c}}an and
             B{\"u}y{\"u}ktahtak{\i}n, {\.{I}}. Esra},
  title   = {An Expandable Machine Learning-Optimization Framework
             for Sequential Decision-Making},
  journal = {European Journal of Operational Research},
  volume  = {314},
  number  = {1},
  pages   = {280--296},
  year    = {2024},
  doi     = {10.1016/j.ejor.2023.10.045}
}

@article{yilmaz2025nonanticipative,
  author  = {Yilmaz, Do{\u{g}}a{\c{c}}an and
             B{\"u}y{\"u}ktahtak{\i}n, {\.{I}}. Esra},
  title   = {A Non-Anticipative Learning-Optimization Framework for
             Solving Multi-Stage Stochastic Programs},
  journal = {Annals of Operations Research},
  volume  = {355},
  number  = {3},
  pages   = {2859--2899},
  year    = {2025},
  doi     = {10.1007/s10479-024-06100-7}
}

@article{YilmazBuyuktahtakin2024,
  author  = {Yilmaz, Do{\u{g}}a{\c{c}}an and
             B{\"u}y{\"u}ktahtak{\i}n, {\.{I}}. Esra},
  title   = {A Deep Reinforcement Learning Framework for Solving
             Two-Stage Stochastic Programs},
  journal = {Optimization Letters},
  volume  = {18},
  pages   = {1993--2020},
  year    = {2024},
  doi     = {10.1007/s11590-023-02009-5}
}

@article{zou2019sddip,
  author  = {Zou, Jikai and Ahmed, Shabbir and Sun, Xu Andy},
  title   = {Stochastic Dual Dynamic Integer Programming},
  journal = {Mathematical Programming},
  volume  = {175},
  number  = {1},
  pages   = {461--502},
  year    = {2019},
  doi     = {10.1007/s10107-018-1249-5}
}

@article{elbalghiti2023generalization,
  author  = {El Balghiti, Othman and Elmachtoub, Adam N. and
             Grigas, Paul and Tewari, Ambuj},
  title   = {Generalization Bounds in the Predict-Then-Optimize
             Framework},
  journal = {Mathematics of Operations Research},
  volume  = {48},
  number  = {4},
  pages   = {2043--2065},
  year    = {2023},
  doi     = {10.1287/moor.2022.1330}
}

@article{kallus2023stochastic,
  author  = {Kallus, Nathan and Mao, Xiaojie},
  title   = {Stochastic Optimization Forests},
  journal = {Management Science},
  volume  = {69},
  number  = {4},
  pages   = {1975--1994},
  year    = {2023},
  doi     = {10.1287/mnsc.2022.4458}
}

@article{misic2020treeensembles,
  author  = {Mi{\v{s}}i{\'c}, Velibor V.},
  title   = {Optimization of Tree Ensembles},
  journal = {Operations Research},
  volume  = {68},
  number  = {5},
  pages   = {1605--1624},
  year    = {2020},
  doi     = {10.1287/opre.2019.1928}
}

@article{oroojlooyjadid2022beergame,
  author  = {Oroojlooyjadid, Afshin and Nazari, MohammadReza and
             Snyder, Lawrence V. and Tak{\'a}{\v{c}}, Martin},
  title   = {A Deep {Q}-Network for the Beer Game: Deep Reinforcement
             Learning for Inventory Optimization},
  journal = {Manufacturing \& Service Operations Management},
  volume  = {24},
  number  = {1},
  pages   = {285--304},
  year    = {2022},
  doi     = {10.1287/msom.2020.0939}
}

@article{harsha2025deeppolicy,
  author  = {Harsha, Pavithra and Jagmohan, Ashish and
             Kalagnanam, Jayant and Quanz, Brian and
             Singhvi, Divya},
  title   = {Deep Policy Iteration with Integer Programming for
             Inventory Management},
  journal = {Manufacturing \& Service Operations Management},
  volume  = {27},
  number  = {2},
  pages   = {369--388},
  year    = {2025},
  doi     = {10.1287/msom.2022.0617}
}

@article{mandl2023commodity,
  author  = {Mandl, Christian and Minner, Stefan},
  title   = {Data-Driven Optimization for Commodity Procurement
             Under Price Uncertainty},
  journal = {Manufacturing \& Service Operations Management},
  volume  = {25},
  number  = {2},
  pages   = {371--390},
  year    = {2023},
  doi     = {10.1287/msom.2020.0890}
}

@article{busch2023commoditydeep,
  author  = {Busch, Nicolas and Cr{\"o}nert, Tobias and
             Minner, Stefan and Rettinger, Moritz and
             Sel, Burakhan},
  title   = {Deep Learning for Commodity Procurement: Nonlinear
             Data-Driven Optimization of Hedging Decisions},
  journal = {{INFORMS} Journal on Optimization},
  volume  = {5},
  number  = {3},
  pages   = {273--294},
  year    = {2023},
  doi     = {10.1287/ijoo.2022.0086}
}

@article{liu2022mlepomdp,
  author  = {Liu, Zeyu and Khojandi, Anahita and Li, Xueping and
             Mohammed, Akram and Davis, Robert L. and
             Kamaleswaran, Rishikesan},
  title   = {A Machine Learning--Enabled Partially Observable
             {Markov} Decision Process Framework for Early Sepsis
             Prediction},
  journal = {{INFORMS} Journal on Computing},
  volume  = {34},
  number  = {4},
  pages   = {2039--2057},
  year    = {2022},
  doi     = {10.1287/ijoc.2022.1176}
}

@article{saghafian2024adtr,
  author  = {Saghafian, Soroush},
  title   = {Ambiguous Dynamic Treatment Regimes: A Reinforcement
             Learning Approach},
  journal = {Management Science},
  volume  = {70},
  number  = {9},
  pages   = {5667--5690},
  year    = {2024},
  doi     = {10.1287/mnsc.2022.00883}
}

@article{gautron2022crop,
  author  = {Gautron, Romain and Maillard, Odalric-Ambrym and
             Preux, Philippe and Corbeels, Marc and
             Sabbadin, R{\'e}gis},
  title   = {Reinforcement Learning for Crop Management Support:
             Review, Prospects and Challenges},
  journal = {Computers and Electronics in Agriculture},
  volume  = {200},
  pages   = {107182},
  year    = {2022},
  doi     = {10.1016/j.compag.2022.107182}
}

@article{pan2021deepopf,
  author  = {Pan, Xiang and Zhao, Tianyu and Chen, Minghua and
             Zhang, Shengyu},
  title   = {{DeepOPF}: A Deep Neural Network Approach for
             Security-Constrained {DC} Optimal Power Flow},
  journal = {IEEE Transactions on Power Systems},
  volume  = {36},
  number  = {3},
  pages   = {1725--1735},
  year    = {2021},
  doi     = {10.1109/TPWRS.2020.3026379}
}

@article{delage2010dro,
  author  = {Delage, Erick and Ye, Yinyu},
  title   = {Distributionally Robust Optimization Under Moment
             Uncertainty with Application to Data-Driven Problems},
  journal = {Operations Research},
  volume  = {58},
  number  = {3},
  pages   = {595--612},
  year    = {2010},
  doi     = {10.1287/opre.1090.0741}
}

@article{wiesemann2014dro,
  author  = {Wiesemann, Wolfram and Kuhn, Daniel and Sim, Melvyn},
  title   = {Distributionally Robust Convex Optimization},
  journal = {Operations Research},
  volume  = {62},
  number  = {6},
  pages   = {1358--1376},
  year    = {2014},
  doi     = {10.1287/opre.2014.1314}
}

@article{esfahani2018wasserstein,
  author  = {Mohajerin Esfahani, Peyman and Kuhn, Daniel},
  title   = {Data-Driven Distributionally Robust Optimization Using
             the {Wasserstein} Metric: Performance Guarantees and
             Tractable Reformulations},
  journal = {Mathematical Programming},
  volume  = {171},
  number  = {1--2},
  pages   = {115--166},
  year    = {2018},
  doi     = {10.1007/s10107-017-1172-1}
}

@article{xu2012drmdp,
  author  = {Xu, Huan and Mannor, Shie},
  title   = {Distributionally Robust {Markov} Decision Processes},
  journal = {Mathematics of Operations Research},
  volume  = {37},
  number  = {2},
  pages   = {288--300},
  year    = {2012},
  doi     = {10.1287/moor.1120.0540}
}

@inproceedings{sinha2018certifying,
  author    = {Sinha, Aman and Namkoong, Hongseok and Duchi, John},
  title     = {Certifying Some Distributional Robustness with
               Principled Adversarial Training},
  booktitle = {International Conference on Learning Representations
               ({ICLR} 2018)},
  year      = {2018},
  doi       = {10.48550/arXiv.1710.10571},
  url       = {https://arxiv.org/abs/1710.10571}
}

@article{duchi2021uniform,
  author  = {Duchi, John C. and Namkoong, Hongseok},
  title   = {Learning Models with Uniform Performance via
             Distributionally Robust Optimization},
  journal = {The Annals of Statistics},
  volume  = {49},
  number  = {3},
  pages   = {1378--1406},
  year    = {2021},
  doi     = {10.1214/20-AOS2004}
}

@inproceedings{sagawa2020groupdro,
  author    = {Sagawa, Shiori and Koh, Pang Wei and
               Hashimoto, Tatsunori B. and Liang, Percy},
  title     = {Distributionally Robust Neural Networks for Group
               Shifts: On the Importance of Regularization for
               Worst-Case Generalization},
  booktitle = {International Conference on Learning Representations
               ({ICLR} 2020)},
  year      = {2020},
  doi       = {10.48550/arXiv.1911.08731},
  url       = {https://arxiv.org/abs/1911.08731}
}

@article{duchi2023latentdro,
  author  = {Duchi, John and Hashimoto, Tatsunori and
             Namkoong, Hongseok},
  title   = {Distributionally Robust Losses for Latent Covariate
             Mixtures},
  journal = {Operations Research},
  volume  = {71},
  number  = {2},
  pages   = {649--664},
  year    = {2023},
  doi     = {10.1287/opre.2022.2363}
}

@article{romeraparedes2024funsearch,
  author  = {Romera-Paredes, Bernardino and Barekatain,
             Mohammadamin and Novikov, Alexander and Balog, Matej
             and Kumar, M. Pawan and Dupont, Emilien and
             Ruiz, Francisco J. R. and Ellenberg, Jordan S. and
             Wang, Pengming and Fawzi, Omar and Kohli, Pushmeet
             and Fawzi, Alhussein},
  title   = {Mathematical Discoveries from Program Search with
             Large Language Models},
  journal = {Nature},
  volume  = {625},
  pages   = {468--475},
  year    = {2024},
  doi     = {10.1038/s41586-023-06924-6}
}

@article{blekos2024qaoa,
  author  = {Blekos, Kostas and Brand, Dean and Ceschini, Andrea
             and Chou, Chiao-Hui and Li, Rui-Hao and
             Pandya, Komal and Summer, Alessandro},
  title   = {A Review on Quantum Approximate Optimization Algorithm
             and its Variants},
  journal = {Physics Reports},
  volume  = {1068},
  pages   = {1--66},
  year    = {2024},
  doi     = {10.1016/j.physrep.2024.03.002}
}

@article{yang2025advisor,
  author  = {Yang, Cathy (Liu) and Bauer, Kevin and Li, Xitong
             and Hinz, Oliver},
  title   = {My Advisor, Her {AI}, and Me: Evidence from a Field
             Experiment on Human--{AI} Collaboration and Investment
             Decisions},
  journal = {Management Science},
  volume  = {72},
  number  = {1},
  pages   = {242--264},
  year    = {2025},
  doi     = {10.1287/mnsc.2022.03918}
}

@article{vaccaro2024useful,
  author  = {Vaccaro, Michelle and Almaatouq, Abdullah and
             Malone, Thomas},
  title   = {When Combinations of Humans and {AI} Are Useful:
             A Systematic Review and Meta-Analysis},
  journal = {Nature Human Behaviour},
  volume  = {8},
  pages   = {2293--2303},
  year    = {2024},
  doi     = {10.1038/s41562-024-02024-1}
}

@misc{vandenoord2016wavenet,
  author       = {van den Oord, Aaron and Dieleman, Sander and
                  Zen, Heiga and Simonyan, Karen and Vinyals, Oriol
                  and Graves, Alex and Kalchbrenner, Nal and
                  Senior, Andrew and Kavukcuoglu, Koray},
  title        = {{WaveNet}: A Generative Model for Raw Audio},
  howpublished = {arXiv preprint arXiv:1609.03499},
  year         = {2016},
  doi          = {10.48550/arXiv.1609.03499},
  url          = {https://arxiv.org/abs/1609.03499}
}

@misc{bai2018empirical,
  author       = {Bai, Shaojie and Kolter, J. Zico and
                  Koltun, Vladlen},
  title        = {An Empirical Evaluation of Generic Convolutional
                  and Recurrent Networks for Sequence Modeling},
  howpublished = {arXiv preprint arXiv:1803.01271},
  year         = {2018},
  doi          = {10.48550/arXiv.1803.01271},
  url          = {https://arxiv.org/abs/1803.01271}
}

@article{buyuktahtakin2018ebola,
  author  = {B{\"u}y{\"u}ktahtak{\i}n, {\.I}. Esra and
             des-Bordes, Emmanuel and
             K{\i}b{\i}{\c{s}}, Eyy{\"u}b Y.},
  title   = {A New Epidemics--Logistics Model: Insights into
             Controlling the {Ebola} Virus Disease in
             {West Africa}},
  journal = {European Journal of Operational Research},
  volume  = {265},
  number  = {3},
  pages   = {1046--1063},
  year    = {2018},
  doi     = {10.1016/j.ejor.2017.08.037}
}

@article{kibis2021eab,
  author  = {K{\i}b{\i}{\c{s}}, Eyy{\"u}b Y. and
             B{\"u}y{\"u}ktahtak{\i}n, {\.I}. Esra and
             Haight, Robert G. and Akhundov, Najmaddin and
             Knight, Kathleen and Flower, Charles E.},
  title   = {A Multistage Stochastic Programming Approach to the
             Optimal Surveillance and Control of the Emerald Ash
             Borer in Cities},
  journal = {{INFORMS} Journal on Computing},
  volume  = {33},
  number  = {2},
  pages   = {808--834},
  year    = {2021},
  doi     = {10.1287/ijoc.2020.0963}
}

@inproceedings{lee2025transformer,
  author    = {Lee, Minseo and B{\"u}y{\"u}ktahtak{\i}n,
               {\.{I}}. Esra},
  title     = {Mathematical Formulation of Transformer
               Architecture},
  booktitle = {Proceedings of the 20th {INFORMS} Data Mining and
               Decision Analytics ({DMDA}) Workshop},
  year      = {2025},
  note      = {Finalist paper in the Workshop Best Paper
               Competition}
}

@article{cosgun2018hiv,
  author  = {Co{\c{s}}gun, {\"O}zlem and
             B{\"u}y{\"u}ktahtak{\i}n, {\.I}. Esra},
  title   = {Stochastic Dynamic Resource Allocation for {HIV}
             Prevention and Treatment: An Approximate Dynamic
             Programming Approach},
  journal = {Computers \& Industrial Engineering},
  volume  = {118},
  pages   = {423--439},
  year    = {2018},
  doi     = {10.1016/j.cie.2018.01.018}
}

@article{yin2023ventilator,
  author  = {Yin, Xuecheng and B{\"u}y{\"u}ktahtak{\i}n,
             {\.I}. Esra and Patel, Bhumi P.},
  title   = {{COVID-19}: Data-Driven Optimal Allocation of
             Ventilator Supply Under Uncertainty and Risk},
  journal = {European Journal of Operational Research},
  volume  = {304},
  number  = {1},
  pages   = {255--275},
  year    = {2023},
  doi     = {10.1016/j.ejor.2021.11.052}
}

@article{cibaku2026vaccine,
  author  = {Cibaku, Elson and B{\"u}y{\"u}ktahtak{\i}n,
             {\.I}. Esra},
  title   = {An Adaptive {K}-means and Reinforcement Learning
             ({RL}) Algorithm to Effective Vaccine Distribution},
  journal = {Computers \& Operations Research},
  volume  = {185},
  pages   = {107275},
  year    = {2026},
  doi     = {10.1016/j.cor.2025.107275}
}

@article{ivanov2020viability,
  author  = {Ivanov, Dmitry and Dolgui, Alexandre},
  title   = {Viability of Intertwined Supply Networks: Extending
             the Supply Chain Resilience Angles Towards
             Survivability. {A} Position Paper Motivated by
             {COVID-19} Outbreak},
  journal = {International Journal of Production Research},
  volume  = {58},
  number  = {10},
  pages   = {2904--2915},
  year    = {2020},
  doi     = {10.1080/00207543.2020.1750727}
}

@article{ivanov2023idt,
  author  = {Ivanov, Dmitry},
  title   = {Intelligent Digital Twin ({iDT}) for Supply Chain
             Stress-Testing, Resilience, and Viability},
  journal = {International Journal of Production Economics},
  volume  = {263},
  pages   = {108938},
  year    = {2023},
  doi     = {10.1016/j.ijpe.2023.108938}
}

@article{jozani2025proactive,
  author       = {Jozani, Kimiya and Tunc, Sait and Eldardiry, Hoda and Sageer, Nihal A. and
                  B{\"u}y{\"u}ktahtak{\i}n, {\.{I}}. Esra},
  title        = {Adaptive Multi-Echelon Epidemic--Supply Chain Optimization
                  at National Scale},
  journal      = {IISE Transactions},
  year         = {2026},
  note         = {Accepted for publication},
  doi          = {https://doi.org/10.1080/24725854.2026.2698692},
  url          = {https://doi.org/10.1080/24725854.2026.2698692}
}

@article{kibis2019cancer,
  author  = {K{\i}b{\i}{\c{s}}, Eyy{\"u}b Y. and
             B{\"u}y{\"u}ktahtak{\i}n, {\.I}. Esra},
  title   = {Optimizing Multi-Modal Cancer Treatment Under {3D}
             Spatio-Temporal Tumor Growth},
  journal = {Mathematical Biosciences},
  volume  = {307},
  pages   = {53--69},
  year    = {2019},
  doi     = {10.1016/j.mbs.2018.10.010},
  url     = {https://doi.org/10.1016/j.mbs.2018.10.010}
}

@inproceedings{buyuktahtakin2024transformers,
  author    = {Cooper, Joshua F. and Choi, Seung Jin and
               B{\"u}y{\"u}ktahtak{\i}n, {\.{I}}. Esra},
  title     = {Toward Transf{OR}mers: Revolutionizing the Solution
               of Mixed Integer Programs with Transformers},
  booktitle = {Proceedings of the 2024 Industrial and Systems
               Engineering Research Conference ({ISERC})},
  address   = {Montreal, Canada},
  year      = {2024},
  doi       = {10.48550/arXiv.2402.13380},
  url       = {https://arxiv.org/abs/2402.13380}
}

@inproceedings{choi2024tcnn,
  author    = {Choi, Seung Jin and Cooper, Joshua and
              B{\"u}y{\"u}ktahtak{\i}n, {\.{I}}. Esra},
  title     = {A Temporal Convolutional Neural Network ({TCNN})
               Approach to Predicting Capacitated Lot-Sizing
               Solutions},
  booktitle = {Proceedings of the 2024 {IISE} Annual Conference
               \& Expo},
  pages     = {1--6},
  year      = {2024},
  publisher = {Institute of Industrial and Systems Engineers
               ({IISE})},
  doi       = {10.21872/2024IISE_7151}
}

@inproceedings{choi2025benderstransformers,
  author    = {Choi, Seung Jin and Jozani, Kimiya and
               Cooper, Joshua F. and
               B{\"u}y{\"u}ktahtak{\i}n, {\.{I}}. Esra},
  title     = {Learning to Optimize at Scale: A Benders
               Decomposition-Transf{OR}mers Framework for
               Stochastic Combinatorial Optimization},
  booktitle = {{NeurIPS} 2025 Workshop {MLxOR}: Mathematical
               Foundations and Operational Integration of Machine
               Learning for Uncertainty-Aware Decision-Making},
  year      = {2025},
  note      = {Poster paper, published on OpenReview},
  url       = {https://openreview.net/forum?id=jVcPvWjrQ5}
}

@article{benders1962partitioning,
  author  = {Benders, J. F.},
  title   = {Partitioning Procedures for Solving Mixed-Variables
             Programming Problems},
  journal = {Numerische Mathematik},
  volume  = {4},
  number  = {1},
  pages   = {238--252},
  year    = {1962},
  doi     = {10.1007/BF01386316}
}

@inproceedings{tong2024ottnn,
  author    = {Tong, Jiatai and Cai, Junyang and Serra, Thiago},
  title     = {Optimization over Trained Neural Networks: Taking
               a Relaxing Walk},
  booktitle = {Integration of Constraint Programming, Artificial
               Intelligence, and Operations Research},
  series    = {Lecture Notes in Computer Science},
  pages     = {221--233},
  year      = {2024},
  publisher = {Springer},
  doi       = {10.1007/978-3-031-60599-4_14}
}

@inproceedings{kotary2021survey,
  author    = {Kotary, James and Fioretto, Ferdinando and
               Van Hentenryck, Pascal and Wilder, Bryan},
  title     = {End-to-End Constrained Optimization Learning:
               A Survey},
  booktitle = {Proceedings of the Thirtieth International Joint
               Conference on Artificial Intelligence
               ({IJCAI-21})},
  pages     = {4475--4482},
  year      = {2021},
  doi       = {10.24963/ijcai.2021/610}
}

@article{sadana2025survey,
  author  = {Sadana, Utsav and Chenreddy, Abhilash and
             Delage, Erick and Forel, Alexandre and
             Frejinger, Emma and Vidal, Thibaut},
  title   = {A Survey of Contextual Optimization Methods for
             Decision-Making under Uncertainty},
  journal = {European Journal of Operational Research},
  volume  = {320},
  number  = {2},
  pages   = {271--289},
  year    = {2025},
  doi     = {10.1016/j.ejor.2024.03.020}
}

@article{kaelbling1998pomdp,
  author  = {Kaelbling, Leslie Pack and Littman, Michael L. and
             Cassandra, Anthony R.},
  title   = {Planning and Acting in Partially Observable
             Stochastic Domains},
  journal = {Artificial Intelligence},
  volume  = {101},
  number  = {1--2},
  pages   = {99--134},
  year    = {1998},
  doi     = {10.1016/S0004-3702(98)00023-X}
}

@misc{hausknecht2015drqn,
  author       = {Hausknecht, Matthew and Stone, Peter},
  title        = {Deep Recurrent Q-Learning for Partially
                  Observable {MDPs}},
  howpublished = {arXiv preprint arXiv:1507.06527},
  year         = {2015},
  doi          = {10.48550/arXiv.1507.06527},
  url          = {https://arxiv.org/abs/1507.06527}
}

@article{yin2021equity,
  author  = {Yin, Xuecheng and B{\"u}y{\"u}ktahtak{\i}n,
             {\.I}. Esra},
  title   = {A Multi-Stage Stochastic Programming Approach to
             Epidemic Resource Allocation with Equity
             Considerations},
  journal = {Health Care Management Science},
  volume  = {24},
  number  = {3},
  pages   = {597--622},
  year    = {2021},
  doi     = {10.1007/s10729-021-09559-z}
}

@inproceedings{gasse2019learn2branch,
  author    = {Gasse, Maxime and Ch{\'e}telat, Didier and
               Ferroni, Nicola and Charlin, Laurent and
               Lodi, Andrea},
  title     = {Exact Combinatorial Optimization with Graph
               Convolutional Neural Networks},
  booktitle = {Advances in Neural Information Processing Systems},
  volume    = {32},
  pages     = {15554--15566},
  year      = {2019},
  publisher = {Curran Associates, Inc.},
  doi       = {10.48550/arXiv.1906.01629}
}

@inproceedings{khalil2022mipgnn,
  author    = {Khalil, Elias B. and Morris, Christopher and
               Lodi, Andrea},
  title     = {{MIP-GNN}: A Data-Driven Framework for Guiding
               Combinatorial Solvers},
  booktitle = {Proceedings of the {AAAI} Conference on Artificial
               Intelligence},
  volume    = {36},
  number    = {9},
  pages     = {10219--10227},
  year      = {2022},
  doi       = {10.1609/aaai.v36i9.21262}
}

@inproceedings{patel2022neur2sp,
  author    = {Dumouchelle, Justin and Patel, Rahul
               and Khalil, Elias B. and Bodur, Merve},
  title     = {{Neur2SP}: Neural Two-Stage Stochastic
               Programming},
  booktitle = {Advances in Neural Information Processing Systems},
  volume    = {35},
  pages     = {23992--24005},
  year      = {2022},
  publisher = {Curran Associates, Inc.},
  doi       = {10.48550/arXiv.2205.12006}
}

@article{buyuktahtakin2011dynamic,
  author  = {B{\"u}y{\"u}ktahtak{\i}n, {\.{I}}. Esra and
             Feng, Zhuo and Frisvold, George and
             Szidarovszky, Ferenc and Olsson, Aaryn},
  title   = {A Dynamic Model of Controlling Invasive Species},
  journal = {Computers \& Mathematics with Applications},
  volume  = {62},
  number  = {9},
  pages   = {3326--3333},
  year    = {2011},
  doi     = {10.1016/j.camwa.2011.08.037}
}

@inproceedings{fioretto2020predicting,
  author    = {Fioretto, Ferdinando and Mak, Terrence W. K. and
               Van Hentenryck, Pascal},
  title     = {Predicting {AC} Optimal Power Flows: Combining
               Deep Learning and {Lagrangian} Dual Methods},
  booktitle = {Proceedings of the {AAAI} Conference on Artificial
               Intelligence},
  volume    = {34},
  number    = {01},
  pages     = {630--637},
  year      = {2020},
  doi       = {10.1609/aaai.v34i01.5403}
}

@inproceedings{choi2026safety,
  author    = {Choi, Seung Jin and Cibaku, Elson and
               Svirsko, Anna and Skipper, Daphne and
               B{\"u}y{\"u}ktahtak{\i}n, {\.{I}}. Esra},
  title     = {Safety-Constrained Reinforcement Learning for
               Naval Warfare Searching with an Intelligent
               Target},
  booktitle = {Refereed Proceedings of the 2026 {INFORMS}
               Optimization Society Conference ({IOS} 2026)},
  address   = {Atlanta, GA},
  year      = {2026}
}

@misc{baswapuram2026principalagent,
  author       = {Baswapuram, Avinash Kumar and Chen, Chen and
                  Cai, Wenbo and
                  B{\"u}y{\"u}ktahtak{\i}n, {\.I}. Esra},
  title        = {An Interpretable Ensemble Heuristic for
                  Principal-Agent Games with Machine Learning},
  howpublished = {Optimization Online},
  pages        = {1--33},
  year         = {2026},
  month        = apr,
  note         = {Published online April 12, 2026. Working paper,
                  Virginia Tech. Also available at
                  \url{https://hdl.handle.net/10919/143025}},
  url          = {https://optimization-online.org/?p=34384}
}

@inproceedings{lanctot2017unified,
  author    = {Lanctot, Marc and Zambaldi, Vinicius and
               Gruslys, Audrunas and Lazaridou, Angeliki and
               Tuyls, Karl and Perolat, Julien and
               Silver, David and Graepel, Thore},
  title     = {A Unified Game-Theoretic Approach to Multiagent
               Reinforcement Learning},
  booktitle = {Advances in Neural Information Processing Systems},
  volume    = {30},
  year      = {2017},
  editor    = {Guyon, Isabelle and von Luxburg, Ulrike and
               Bengio, Samy and Wallach, Hanna and Fergus, Rob
               and Vishwanathan, S. V. N. and Garnett, Roman},
  publisher = {Curran Associates, Inc.},
url       = {https://proceedings.neurips.cc/paper_files/paper/2017/hash/3323fe11e9595c09af38fe67567a9394-Abstract.html}
}

@article{moerland2023model,
  author  = {Moerland, Thomas M. and Broekens, Joost and
             Plaat, Aske and Jonker, Catholijn M.},
  title   = {Model-Based Reinforcement Learning: A Survey},
  journal = {Foundations and Trends in Machine Learning},
  volume  = {16},
  number  = {1},
  pages   = {1--118},
  year    = {2023},
  doi     = {10.1561/2200000086}
}

@misc{levine2020offline,
  author       = {Levine, Sergey and Kumar, Aviral and
                  Tucker, George and Fu, Justin},
  title        = {Offline Reinforcement Learning: Tutorial,
                  Review, and Perspectives on Open Problems},
  howpublished = {arXiv preprint arXiv:2005.01643},
  year         = {2020},
  url          = {https://arxiv.org/abs/2005.01643}
}

@article{garcia2015comprehensive,
  author  = {Garc{\'i}a, Javier and Fern{\'a}ndez, Fernando},
  title   = {A Comprehensive Survey on Safe Reinforcement
             Learning},
  journal = {Journal of Machine Learning Research},
  volume  = {16},
  number  = {1},
  pages   = {1437--1480},
  year    = {2015},
  url     = {https://jmlr.org/papers/v16/garcia15a.html}
}

@inproceedings{achiam2017constrained,
  author    = {Achiam, Joshua and Held, David and
               Tamar, Aviv and Abbeel, Pieter},
  title     = {Constrained Policy Optimization},
  booktitle = {Proceedings of the 34th International Conference
               on Machine Learning},
  series    = {Proceedings of Machine Learning Research},
  volume    = {70},
  pages     = {22--31},
  year      = {2017},
  publisher = {PMLR},
  url       = {https://proceedings.mlr.press/v70/achiam17a.html}
}

@incollection{zhang2021multi,
  author    = {Zhang, Kaiqing and Yang, Zhuoran and
               Ba{\c{s}}ar, Tamer},
  title     = {Multi-Agent Reinforcement Learning: A Selective
               Overview of Theories and Algorithms},
  booktitle = {Handbook of Reinforcement Learning and Control},
  editor    = {Vamvoudakis, Kyriakos G. and Wan, Yan and
               Lewis, Frank L. and Cansever, Derya},
  series    = {Studies in Systems, Decision and Control},
  volume    = {325},
  pages     = {321--384},
  publisher = {Springer},
  year      = {2021},
  doi       = {10.1007/978-3-030-60990-0_12}
}

@article{jones1998efficient,
  author  = {Jones, Donald R. and Schonlau, Matthias and
             Welch, William J.},
  title   = {Efficient Global Optimization of Expensive
             Black-Box Functions},
  journal = {Journal of Global Optimization},
  volume  = {13},
  number  = {4},
  pages   = {455--492},
  year    = {1998},
  doi     = {10.1023/A:1008306431147}
}

@misc{frazier2018tutorial,
  author       = {Frazier, Peter I.},
  title        = {A Tutorial on {Bayesian} Optimization},
  howpublished = {arXiv preprint arXiv:1807.02811},
  year         = {2018},
  url          = {https://arxiv.org/abs/1807.02811}
}

@article{caballero2025bayesian,
  author  = {Caballero, William N. and Jenkins, Phillip R. and
             Banks, David and Robbins, Matthew},
  title   = {Bayesian Graph Traversal},
  journal = {Decision Analysis},
  volume  = {23},
  number  = {1},
  pages   = {11--27},
  year    = {2025},
  doi     = {10.1287/deca.2024.0239}
}

@inproceedings{chen2021decision,
  author    = {Chen, Lili and Lu, Kevin and Rajeswaran, Aravind and 
               Lee, Kimin and Grover, Aditya and Laskin, Michael and 
               Abbeel, Pieter and Srinivas, Aravind and Mordatch, Igor},
  title     = {Decision Transformer: Reinforcement Learning via 
               Sequence Modeling},
  booktitle = {Advances in Neural Information Processing Systems},
  volume    = {34},
  pages     = {15084--15097},
  year      = {2021},
  url       = {https://arxiv.org/abs/2106.01345}
}

@article{basciftci2024adaptive,
  author  = {Basciftci, Beste and Ahmed, Shabbir and Gebraeel, Nagi},
  title   = {Adaptive Two-Stage Stochastic Programming with an Analysis on Capacity Expansion Planning Problem},
  journal = {Manufacturing \& Service Operations Management},
  year    = {2024},
  volume  = {26},
  number  = {6},
  doi     = {10.1287/msom.2023.0157},
  url     = {https://doi.org/10.1287/msom.2023.0157}
}

@article{siddig2025multistage,
  author  = {Siddig, Murwan and Song, Yongjia},
  title   = {Multistage Stochastic Programming with a Random Number of Stages: Applications in Hurricane Disaster Relief Logistics Planning},
  journal = {European Journal of Operational Research},
  year    = {2025},
  volume  = {321},
  number  = {3},
  pages   = {925--941},
  publisher = {North-Holland},
  doi       = {10.1016/j.ejor.2024.10.004},
  url       = {https://doi.org/10.1016/j.ejor.2024.10.004}
}

@article{daryalal2024lagrangian,
  author  = {Daryalal, Maryam and Bodur, Merve and Luedtke, James R.},
  title   = {Lagrangian Dual Decision Rules for Multistage Stochastic 
             Mixed-Integer Programming},
  journal = {Operations Research},
  volume  = {72},
  number  = {2},
  pages   = {717--737},
  year    = {2024},
  doi     = {10.1287/opre.2022.2366}
}

@article{boland2018combining,
  author  = {Boland, Natashia and Christiansen, Jeffrey and Dandurand, 
             Brian and Eberhard, Andrew and Linderoth, Jeff and 
             Luedtke, James and Oliveira, Fabricio},
  title   = {Combining Progressive Hedging with a Frank-Wolfe Method 
             to Compute Lagrangian Dual Bounds in Stochastic 
             Mixed-Integer Programming},
  journal = {SIAM Journal on Optimization},
  volume  = {28},
  number  = {2},
  pages   = {1312--1336},
  year    = {2018},
  doi     = {10.1137/16M1076290}
}

@article{lulli2004branch,
  author  = {Lulli, Guglielmo and Sen, Suvrajeet},
  title   = {A Branch-and-Price Algorithm for Multistage Stochastic 
             Integer Programming with Application to Stochastic 
             Batch-Sizing Problems},
  journal = {Management Science},
  volume  = {50},
  number  = {6},
  pages   = {786--796},
  year    = {2004},
  doi     = {10.1287/mnsc.1030.0164}
}

@article{liu2021adaptive,
  author  = {Li, Cong and Wang, Yongchao and Liu, Fangzhou and Liu, Qingchen and Buss, Martin},
  title   = {Model-Free Incremental Adaptive Dynamic Programming Based Approximate Robust Optimal Regulation},
  journal = {arXiv preprint arXiv:2105.01698},
  year    = {2021},
  doi     = {10.48550/arXiv.2105.01698},
  url     = {https://arxiv.org/abs/2105.01698}
}

@article{jiang2022chance,
  author  = {Jiang, Ruiwei and K{\"u}{\c{c}}{\"u}kyavuz, Simge},
  title   = {Chance-Constrained Optimization under Limited Distributional Information: A Review of Reformulations Based on Sampling and Distributional Robustness},
  journal = {EURO Journal on Computational Optimization},
  volume  = {10},
  pages   = {100030},
  year    = {2022},
  doi     = {10.1016/j.ejco.2022.100030},
  url     = {https://doi.org/10.1016/j.ejco.2022.100030}
}

@article{homemdemello2014monte,
  author  = {Homem-de-Mello, Tito and Bayraksan, G{\"u}zin},
  title   = {Monte Carlo Sampling-Based Methods for Stochastic Optimization},
  journal = {Surveys in Operations Research and Management Science},
  volume  = {19},
  number  = {1},
  pages   = {56--85},
  year    = {2014},
  doi     = {10.1016/j.sorms.2014.05.001},
  url     = {https://doi.org/10.1016/j.sorms.2014.05.001}
}

@article{morton2026mdp,
  author  = {Morton, David P. and Dowson, Oscar and Pagnoncelli, Bernardo K.},
  title   = {{MDP} Modeling for Multi-Stage Stochastic Programs},
  journal = {Mathematical Programming},
  year    = {2026},
  pages   = {1--36},
  doi     = {10.1007/s10107-026-02368-8},
  url     = {https://doi.org/10.1007/s10107-026-02368-8}
}

@article{rockafellar2000optimization,
  author  = {Rockafellar, R. Tyrrell and Uryasev, Stanislav},
  title   = {Optimization of Conditional Value-at-Risk},
  journal = {Journal of Risk},
  volume  = {2},
  number  = {3},
  pages   = {21--41},
  year    = {2000},
  doi     = {10.21314/JOR.2000.038},
  url     = {https://doi.org/10.21314/JOR.2000.038}
}

@book{dentcheva2024risk,
  author    = {Dentcheva, Darinka and Ruszczy{\'n}ski, Andrzej},
  title     = {Risk-Averse Optimization and Control: Theory and Methods},
  series    = {Springer Series in Operations Research and Financial Engineering},
  publisher = {Springer Cham},
  year      = {2024},
  pages     = {i--xv, 1--451},
  doi       = {10.1007/978-3-031-57988-2},
  isbn      = {978-3-031-57988-2},
  url       = {https://doi.org/10.1007/978-3-031-57988-2}
}

@article{ntaimo2013simulation,
  author  = {Ntaimo, Lewis and Gallego-Arrubla, Julian A. and Gan, Jianbang and Stripling, Curt and Young, Joshua and Spencer, Thomas},
  title   = {A Simulation and Stochastic Integer Programming Approach to Wildfire Initial Attack Planning},
  journal = {Forest Science},
  volume  = {59},
  number  = {1},
  pages   = {105--117},
  year    = {2013},
  doi     = {10.5849/forsci.11-022},
  url     = {https://doi.org/10.5849/forsci.11-022}
}

@article{hijazi2026improving,
  author  = {Hijazi, Amira and {\"O}zalt{\i}n, Osman Y. and Uzsoy, Reha},
  title   = {An Improving Column Is All You Need: Enhancing Column Generation for Parallel Machine Scheduling via Transformers},
  journal = {INFORMS Journal on Computing},
  volume  = {0},
  number  = {0},
  pages   = {1--36},
  year    = {2026},
  doi     = {10.1287/ijoc.2024.1005},
  url     = {https://doi.org/10.1287/ijoc.2024.1005}
}

\clearpage
\section*{Author Biography}

\noindent {\textbf{{\.I.} Esra B{\"u}y{\"u}ktahtak{\i}n}} (``Deep Learning for Sequential
Decision Making under Uncertainty: Foundations, Frameworks, and Frontiers'') is a Full
Professor of Operations Research in the Grado Department of Industrial and Systems
Engineering at Virginia Tech, where she directs the
\href{https://soml.ise.vt.edu}{Systems Optimization and Machine Learning Lab (SysOptiMaL)}.
She earned her Ph.D.\ in Operations Research from the University of Florida. Her research advances multi-stage stochastic mixed-integer programming through optimization theory, algorithmic innovation, and learning-enhanced optimization.
She is recognized for pioneering contributions at the interface of
deep learning and optimization for sequential decision making under uncertainty.
Her work bridges operations research, machine learning, and artificial intelligence, with applications in health systems, epidemic supply chains, biosecurity, defense, ecological conservation, agriculture, forestry, and environmental
sustainability.

Dr.~B{\"u}y{\"u}ktahtak{\i}n is the recipient of the National Science Foundation (NSF)
CAREER Award (2016) and the INFORMS Minority Issues Forum (MIF) Early Career Award (2016),
and her research has been supported by NSF, the U.S. Department of
Agriculture (USDA), the Office of Naval Research (ONR), the U.S. Forest Service, the
Virginia Department of Forestry, the Minnesota Aquatic Invasive Species Research Center
(MAISRC), and the 4-VA Collaborative Research Program, with more than \$3 million in
external funding. She has authored 49 peer-reviewed journal publications in leading scholarly outlets,
including 27 papers in A* and A ranked journals.
Her work has received six INFORMS Best Publication Awards and has been featured three
times in \textit{ISE Magazine}, as well as in the INFORMS Computing Society and U.S. Forest
Service newsletters. Her professional service includes leadership in the INFORMS community, where she served as President of the INFORMS Junior Faculty Interest Group (JFIG, 2014--2015) and received two INFORMS Service Awards. She currently serves as an Associate Editor of the \textit{INFORMS Journal on Computing}.

\end{document}